\documentclass[11pt]{article}
\usepackage{geometry}  
\usepackage[pdftex]{hyperref}
\usepackage{amssymb,amsmath} 
\usepackage{amsthm}
\usepackage{url}
\usepackage{verbatim}
\usepackage{algorithm}
 
\numberwithin{equation}{section}
\usepackage[british]{babel}  
\usepackage[T1]{fontenc}
\usepackage[utf8]{inputenc}
\usepackage[pdftex]{graphicx} 
\usepackage[british]{babel} 
\usepackage[round]{natbib}
\usepackage{fourier,charter} 
\usepackage{color}
\usepackage{placeins}
\let\mathds=\mathbb
\let\mathbbm=\mathbb 
\let\phi=\varphi
\newcommand{\indep}{\perp\kern-0.5em\perp}
\newcommand{\Real}{\mathbb{R}}
\newcommand{\Rlang}{{\normalfont\textsf{R}}{}}
\newcommand{\N}{\mathrm{N}{}}
\newcommand{\SN}{\mathrm{SN}{}}
\newcommand{\ESN}{\mathrm{ESN}{}}
\renewcommand{\d}{\,\mathrm{d}}
\newcommand{\inv}{^{-1}}
\newcommand\NB[1]{\kern-0.4em\raisebox{-1.5ex}{$\stackrel{\big|}{\hbox{%
  \tiny\sc NB}}$}\marginpar{\sl\footnotesize\raggedright\color{red}#1\hfill}\kern-0.2em}
\newcommand{\rv}{r.\,v.{}}  
\title{\huge\bf
      The information matrix of the bivariate \\
      extended skew-normal distribution }
\author{\large
   {\Large\bf Stefano Franco} \\  
   \url{<franco.ste@gmail.com>}\\
   \url{https://www.linkedin.com/in/stefanofranco}
  \and
  {\Large\bf Adelchi Azzalini} \\   
  Dipartimento di Scienze Statistiche \\
  Università degli Studi di Padova, Italia
 }
\date{\today}

\newtheorem{theo}{Teorema} 

\newtheorem{lemma}[theo]{Lemma}
\newtheorem{prop}[theo]{Proposizione}
\newcommand{\ud}{\mathrm{d}}
\newcommand{\var}{\mathbbm{V}\mathrm{ar}}
\newcommand{\mean}{\mathbbm{E}}
\allowdisplaybreaks[1]

\begin{document}
\maketitle
\begin{abstract}\noindent
For the  extended skew-normal distribution, which represents an extension of
the normal (or Gaussian) distribution, we focus on the properties of
the log-likelihood function and derived quantities in the the bivariate case.
Specifically, we derive explicit expressions for the score function and 
the information matrix, in the observed and the expected form;
these do not appear to have been examined before in the literature.
Corresponding computing code in \Rlang\ language is provided, 
which implements the formal expressions.
\end{abstract}
\section{Background and aim} \label{s:intro}

\subsection{The skew-normal and some related distributions}

The skew-normal (SN) family of continuous probability distributions 
represents an extension of the normal, or Gaussian, family of distributions. 
In the last 20--25 years, the SN distribution and a number of related ones, 
which arise from a similar type of construction, 
have received substantial attention in the statistical literature.
Within this context, the present contribution refers to the extended 
skew-normal (ESN) distribution, specifically to its bivariate case.

To start from the most basic case, the expression of the univariate 
SN density  function featuring only a shape parameter $\alpha\in\Real$ is
\begin{equation}
        \label{eq:skewnormZ}
  f_{\SN}(z;\alpha) = 2\: \phi(z)\:\Phi(\alpha z), \qquad -\infty < z < \infty,
\end{equation}
where $\phi$ and $\Phi$ denote the $\N(0,1)$ density and distribution function.
It is clear that, if we set $\alpha=0$, (\ref{eq:skewnormZ}) reproduces the $\N(0,1)$ density;
otherwise the density is asymmetric, and the direction of the asymmetry 
is in agreement with the sign of $\alpha$. For practical work, 
additional parameters for regulating location and scale are 
introduced, as follows. If $Z$ is a random variable (\rv) with 
density (\ref{eq:skewnormZ}), the distribution of $Y=\xi+\omega\,Z$
where $\omega>0$ involves also a location $\xi$ and a scale parameter $\omega$;
in this case, we write  $Y\sim\SN(\xi, \omega^2, \alpha)$.

The extended skew-normal (ESN) density introduces another parameter, 
$\tau$, which allows further regulation of the shape;
see \citet[Section 3.3]{art:azza85} up to a minor change of parameterization.
Accordingly, density (\ref{eq:skewnormZ})  evolves into
\begin{equation}
  \label{eq:esnZ}
   f_{\ESN}(z;\alpha,\tau) = 
        \frac{1}{\Phi(\tau)}\: \phi(z) \:
        \Phi\left(\tau \sqrt{1+\alpha^2} + \alpha z\right),
        \qquad -\infty < z < \infty,
\end{equation}
for any $\tau\in\Real$. If $\tau=0$, we clearly are back to (\ref{eq:skewnormZ}).

The univariate SN and ESN distribution are simple instances of a wide set
of constructions, which we hinted at in the opening paragraph. 
This broader theme is developed in the  monograph  of \cite{azza:capi:2014},
which we shall extensively use as our primary reference source; 
for brevity, we shall refer to this publication as `the SN book'.

\subsection{Aim and origin of the present note}

The central theme of the present note is the bivariate version of
the ESN distribution.
Our specific target is the development of explicit expressions for its
log-likelihood function and some derived quantities, notably the score 
function and the information matrix, in the observed and the expected form.
Computer code which implements these expressions in the \Rlang\ 
programming language is available in the appendix.

Towards this program, Section~\ref{cap:2}  recalls a number of formal properties  
of the univariate and multivariate ESN distributions, with special focus on
the bivariate case.
Proofs of these properties are available in the SN book and will not replicated here; 
similarly, we do not recall a range of results not required for our development.
As a general indication, Sections~2.2 and 3.3.2 of the SN book deal with
the univariate ESN distribution; Section~5.3 covers the 
multivariate version. 

After recalling briefly a number of known facts in Section~\ref{cap:2}, 
our main development takes place in Section~\ref{cap:t3}; 
this delivers the expressions of the functions of interest, 
followed by the \Rlang\ code in the appendix.
Apart from the language translation and some minor changes, the content of 
Section~\ref{cap:t3} and  the appendix is the nearly intact transposition 
of the corresponding units in the 2011 dissertation presented by the 
first author for his  master-level degree at the University of Padua; 
see \citet{franco:2011}.

\section{The extended skew-normal distribution}
\label{cap:2}

\subsection{The univariate ESN distribution}

The basic expression of the univariate ESN density has already been introduced
in equation (\ref{eq:esnZ}). 
The proof that this expression actually integrates to $1$ is given,
along with other results to be recalled next, in Section~2.2 of the SN book.

Similarly to the SN distribution, location and scale parameters are
usually introduced via the transformation $Y=\xi+\omega\,Z$ where $Z$ is
taken to have density  (\ref{eq:esnZ}). 
The density function of $Y$ is
\begin{equation}
\label{eq:esnY}
    f_{ESN}(y;\xi,\omega,\alpha,\tau) = \frac{1}{\omega \Phi(\tau)} \:
        \phi\left(\frac{y- \xi}{\omega}\right) \:
        \Phi\left(\alpha_0(\tau) + \alpha \frac{y-\xi}{\omega}\right),
        \qquad -\infty < y < \infty,
\end{equation}
where $\alpha_0(\tau)=\tau \sqrt{1+\alpha^2}$. 
In this case we write $Y \sim \ESN(\xi, \omega^2, \alpha, \tau)$.

The moment generating function of $Y$ has a simple explicit expression. 
For our purposes, it is sufficient to recall the expression for the normalized 
case  with $\xi=0$, $\omega=1$, which is
\begin{equation}
 \label{eq:ESNfgm}
  M_Z(t) = \frac{1}{\Phi(\tau)} \exp\left(\frac{t^2}{2}\right)\:\Phi(\tau +\delta t );
\end{equation}
where
\[
  \delta = \frac{\alpha}{\sqrt{1+\alpha^{2}}}
\] 
lies in $(-1, 1)$.
From $M_Z(t)$, we obtain the cumulant generating function, that is,
\begin{equation}
  \label{eq:ESNfgk}
  K_Z(t) = \log M_Z(t) = -\log\Phi(\tau) + \frac{t^2}{2} + \zeta_0(\tau + \delta t)
\end{equation}
which has been written using the function $\zeta_0(x)$, whose definition is given, 
along with expressions of its derivatives,  as follows:
\begin{eqnarray}
  \zeta_0(x) &=& \log(\Phi(x)), \label{eq:zeta_0cap2} \\
  \zeta_1(x) &=& \frac{\d~}{\d x} \zeta_0(x)\, =\, \frac{\phi(x)}{\Phi(x)}, \label{eq:zeta1} \\
  \zeta_2(x) &=&\frac{\d~}{\d x} \zeta_1(x)\, =\, -\big(\zeta_1(x)x + [\zeta_1(x)]^2\big)  \label{eq:zeta2} \\
  \zeta_m(x) &=& \frac{\d^m}{\d x^m} \zeta_0(x), \qquad (m=1,2,\ldots).
\end{eqnarray}
This set of functions will occur repeatedly in our development.

Successive differentiation of (\ref{eq:ESNfgk}) and evaluation at $t=0$ delivers
the cumulants of $Z$. In particular, after some simple algebra, we arrive at
\begin{eqnarray}
  \mean \left[Z\right] &=& \zeta_1(\tau) \delta, \label{eq:media} \\
  \var \left[Z\right] &=& 1 + \zeta_2(\tau)\,\delta^2\,, \label{eq:varianza}  
\end{eqnarray}
where $-1<\zeta_{2}(\tau)<0$.



\subsection{The multivariate ESN distribution}
\label{caso_multi}

The $k$-dimensional version of the ESN density function is introduced as follows. 
Denote by $\phi_k(y;\Omega)$ the $k$-dimensional density function of a $\N_{k}(0,\Omega)$
variable, where $\Omega=(\Omega_{ij})$ is a $k\times k$ symmetric positive definite matrix.
Also denote by $\bar\Omega=\omega\inv\Omega\omega\inv$ the corresponding correlation matrix, where 
\begin{equation} 
  \label{eq:omega_cap2}
  \omega = \text{diag} (\Omega_{11},\Omega_{22},\dots,\Omega_{kk})^\frac{1}{2},
\end{equation}
and let $\xi, \alpha$ denote two arbitrary $k$-dimensional vectors.
Then 
\begin{equation}
\label{eq:esnd}
f(y;\xi,\Omega, \alpha, \tau) = \phi_k(y-\xi;\Omega) \:
   \frac{\Phi(\alpha_0 + \alpha^T \omega^{-1}(y-\xi))}{\Phi(\tau)},
    \qquad y\in \mathds{R}^k,
\end{equation}   
where 
\begin{eqnarray}
\label{eq:alpha0} 
\alpha_{0} &=& \tau\:\left(1+\alpha^T \bar{\Omega}\alpha\right)^{1/2} = \tau\sqrt{1+\alpha^2_*}, \\
\label{eq:alphaasterisco}
    \alpha^2_{*}&=& \alpha^T\bar{\Omega}\alpha,
\end{eqnarray}
is the density function of an ESN distribution with location $\xi$, scale matrix $\Omega$,
and shape parameter $\alpha$. 
For a random variable with this density function, we use the notation
$\ESN_{k}(\xi, \Omega, \alpha, \tau)$.

The cumulant generating function of $Z\sim \ESN_k(0, \bar{\Omega},\alpha,\tau)$ is
\begin{eqnarray}
\label{eq:fngenesn} 
K(t) = \frac{1}{2} t^T\Omega t + \zeta_0(\tau+\delta^T t) - \zeta_0(\tau),
\qquad t\in\Real^k,
\end{eqnarray}
where $\zeta_0(x)$ is defined by (\ref{eq:zeta_0cap2}). 
Differentiation of $K(t)$ and evaluation at $t=0$ leads to
\begin{eqnarray}
  \mean\left[ Z\right] &=& K'(0) =\zeta_1(\tau)\delta,\\ 
  \var\left[ Z\right] &=& K''(0) = \bar{\Omega} + \zeta_2(\tau)\delta\delta^T.
\end{eqnarray}
where $\zeta_1(x)$ and $\zeta_2(x)$ are defined by (\ref{eq:zeta1}) and (\ref{eq:zeta2}).
The $m$th order derivatives of (\ref{eq:fngenesn}) are
\begin{eqnarray}
\frac{\d^m K(t)}{\d t_i \d t_j\cdots \d t_r} = 
    \zeta_m(\tau+\delta^T t)\: \delta_i \delta_j \dots \delta_r, \qquad  (m>2);
\end{eqnarray}
here $\delta_{j}$ denotes the $j$-th entry of $\delta$ and $t_{j}$ denotes 
the $j$-th entry of $t$, for $j=1,\ldots,k$. \emph{The same convention will be used 
throughout the rest of the paper without further specification.} 

The cumulant generating function of $Y\sim\ESN_{k}(\xi, \Omega, \alpha, \tau)$,
having density function (\ref{eq:esnd}), is:
\begin{eqnarray}
\label{eq:fngen}        
K_Y(t) = \xi^T t + \frac{1}{2} t^T\Omega t + \zeta_0(\tau+\delta^T\omega t) - \zeta_0(\tau);
\end{eqnarray}
see \citet{art:capiteal03}.
By differentiation of this function, we obtain:
\begin{eqnarray}
&& \mean\left[ Y\right] = K'_Y(0) = \xi + \zeta_1(\tau)\omega\delta,\\ 
&& \var[Y] = K''_Y(0) = \Omega + \zeta_2(\tau)\omega\delta\delta^T\omega,
\end{eqnarray}
 
\section{Statistical aspects of the $\ESN_{2}$ distribution}
\label{cap:t3}
In this section, we  deal with some inferential aspects of the bivariate 
extended skew-normal distribution ($\ESN_2$).
First, in Subsection~\ref{sec:esntony}, we recall briefly some results for the
univariate case contained in the dissertation of \citet{tesi:canale}%
\footnote{Subsequently published in \cite{canale:2011}.},
and then, in the remaining subsections, we focus on new developments 
for the multivariate case, as summarized next.

In Subection~\ref{sec:versim},  the $\ESN_2$ log-likelihood function 
is introduced, followed by its score function.
In Subsection~\ref{sec:info-obs}, we derive the observed  Fisher 
information matrix.
The computation of the expected Fisher  information matrix
takes place in  Subsection~\ref{sec:info-exp}.
This computation could be achieved
thanks to a new lemma, stated originally in the dissertation of
\cite{franco:2011}, to be recalled in Section~\ref{sec:info-exp}.

These matrices will be compared with the analogous expressions
of the $\SN_2$ case,  obtained by \citet{art:azzarei2008}, and the two sets of
matrices coincide when $\tau=0$; obviously, the row and column corresponding
to $\tau$ do not exist in the $\SN_2$ matrix.

Finally, in Subsection~\ref{s:peculiar}, 
we examine some peculiar aspects of the expected information matrix,
by considering certain block partitions of the parameters, 
and we examine numerically the singularity of this matrix 
as $\alpha \to 0 $ and as $\tau \to \pm \infty$.

%
\subsection{The univariate case}
\label{sec:esntony}

From the work of \citet{tesi:canale}, we see that the univariate $\ESN$ distribution
exhibits a range of issues when it comes to inferential work.
The expected information matrix  turns out to be singular, not only when  $\alpha=0$,
but also when $\tau\rightarrow\pm\infty$. 
In some cases, the determinant of the matrix can approach zero, even 
for finite values of $\tau$.
The profile log-likelihood function for $\tau$ is often rather flat,
which makes it difficult for the numerical optimization algorithms 
to locate the maximum likelihood estimate (MLE).
If it can happen that, when a parameter is moved away from the MLE,
the likelihood function can `adjust' the estimates of the other parameters,
leading in practice to equivalent densities for different parameter sets. 
These aspects indicate a situation of `near non-identifiability' of this model. 

The next lemma, introduced by \citet{tesi:canale}, is useful for computing
the expected information matrix in the scalar case.
\begin{lemma}[Canale, 2008]
  \label{lemmatony}
  Consider a \rv\ $Z \sim \ESN(0,1,\alpha,\tau)$,  denote
  by $\zeta_1(\cdot)$  the function defined in (\ref{eq:zeta1}),and 
  let $g(x) = \alpha \tau + \sqrt{\alpha^2+1} x$; then
  \begin{eqnarray*}
  \mean [ f(Z) \zeta_1(\alpha_0 +\alpha Z )] =  
    \frac{\zeta_1(\tau)}{\sqrt{\alpha^2+1}}~ \mean [f(g^{-1}(V))]
   \end{eqnarray*}
where $V \sim N(0,1)$, \, $\alpha_0=\tau\sqrt{1+\alpha^2}$ and  $f(\cdot)$  
any function from $\Real$ to $\Real$ such that the involved integrals exist. 
\end{lemma}

This lemma will later be generalized by a new lemma for the multidimensional case.

\subsection{Log-likelihood and score functions}
\label{sec:versim}
Consider  $Y \sim ESN_{2}(\xi,\Omega,\alpha,\tau)$ and denote by 
$\theta=(\xi_1,\xi_2,\Omega_{11},\Omega_{12},\Omega_{22},\alpha_1,\alpha_2,\tau)$  
the vector of its parameters;
it is required that  $\Omega_{11}, \Omega_{12}, \Omega_{22}$ identify a positive definitte matrix.
The $\theta$ components represent the so-called `direct parameterization' (DP) 
of the bivariate ESN distribution.

In the following development we re-write the log-determinant of  $\Omega$ as 
\[
  \log\Omega_{11} + \log\Omega_{22} - 2\log\Omega_{12}=
   \log\Omega_{11} + \log\Omega_{22} + \log(1- \lambda^2)    
\]
where
  \begin{eqnarray}
    \label{eq:lambda}   
    \lambda &=& \frac{\Omega_{12}}{\sqrt{\Omega_{11}\Omega_{22}}}. 
  \end{eqnarray}
Given an observation $y=(y_1,y_2)^T $ from $Y$, the log-likelihood function can
be written as 
\begin{eqnarray}
  \label{eq:loglikelihood}
    \ell(\theta;y) &=& \mathrm{constant} \, -\,\zeta_0(\tau) - \frac{1}{2}\log\Omega_{11} 
    - \frac{1}{2}\log\Omega_{22}  - \frac{1}{2}\log \left(1- \lambda^2 \right) +      \notag \\ 
      && \phantom{\ell()}     
      -\frac{1}{2 \left(1- \lambda^2 \right)}
    \left[ \frac{(y_1-\xi_1)^2}{\Omega_{11}} + \frac{(y_2-\xi_2)^2}{\Omega_{22}} 
    -2 \lambda
    \left(\frac{y_1-\xi_1}{\sqrt{\Omega_{11}}}\right) 
    \left(\frac{y_2-\xi_2}{\sqrt{\Omega_{22}}}\right)\right] 
    \notag + \notag\\
      && \phantom{\ell()}
    +\zeta_0 \left( \alpha_0 + \alpha_1 \frac{(y_1-\xi_1)}{\sqrt{\Omega_{11}}} 
    + \alpha_2 \frac{(y_2-\xi_2)} {\sqrt{\Omega_{22}}} \right) 
  \end{eqnarray}
having used  $\zeta_{0}(\cdot)$ as defined by (\ref{eq:zeta_0cap2}), 
and $\alpha_0$ as in (\ref{eq:alpha0}).
On setting  
  \begin{eqnarray*}
    z_1 = \frac{y_1-\xi_1}{\sqrt{\Omega_{11}}} \quad ,\quad z_2 = \frac{y_2-\xi_2}{\sqrt{\Omega_{22}}},
  \end{eqnarray*}
  \begin{eqnarray}
    \label{eq:t}
    t
    &=&\tau(1+\alpha^T \bar{\Omega}\alpha)^{1/2}+\alpha_{1}z_1+\alpha_2 z_2\notag\\
    &=&\tau\sqrt{1+\alpha^2_*} +\alpha_{1}z_1+\alpha_2 z_2\\
    &=&\alpha_0+\alpha_{1}z_1+\alpha_2 z_2 \notag
  \end{eqnarray}
and $\zeta_{1}(\cdot)$ as defined by (\ref{eq:zeta1}), the partial derivatives of the
log-likelihood (\ref{eq:loglikelihood}) with respect to the components of $\theta$
provide the elements of the score function, namely
\begin{eqnarray*}
\frac{\partial\ell}{\partial \xi_1} & = &  
    \frac{1}{\sqrt{\Omega_{11}}}\left[ \frac{z_1-\lambda z_2}{1-\lambda^2} -\alpha_1 \zeta_1(t)  \right],\\
\frac{\partial\ell}{\partial \xi_2} & = &  
    \frac{1}{\sqrt{\Omega_{22}}}\left[ \frac{z_2-\lambda z_1}{1-\lambda^2} -\alpha_2 \zeta_1(t)  \right], \\
\frac{\partial\ell}{\partial \Omega_{11}}  &=& 
    \frac{1}{2\Omega_{11}}  \bigg[ \frac{(z^2_1 + z^2_2-2z_1z_2\lambda)\lambda^2}{(1-\lambda^2)^2} +              \frac{z_1^2-2z_1z_2\lambda -1}{1-\lambda^2} \bigg]-
    \frac{1}{2\Omega_{11}}\left( \frac{\alpha_1\alpha_2\lambda\tau}{\sqrt{1+\alpha^2_*}} +\alpha_1 z_1\right)\zeta_1(t),\\
\frac{\partial\ell}{\partial \Omega_{12}}  & = & 
    \frac{1}{\sqrt{\Omega_{11}\Omega_{22}}} \left[\frac{\lambda+z_1z_2}{1-\lambda^2} -                              \frac{(z_1^2+z_2^2-2z_1z_2\lambda)\lambda}{(1-\lambda^2)^2}\right] + \frac{\alpha_1\alpha_2\tau}{\sqrt{\Omega_{11}\Omega_{22}} \sqrt{1+\alpha^2_*}} \zeta_1(t), \\ 
\frac{\partial\ell}{\partial \Omega_{22}}  & = &  
    \frac{1}{2\Omega_{22}} \bigg[ \frac{(z^2_1 + z^2_2 -2z_1z_2\lambda)\lambda^2}{(1-\lambda^2)^2} 
    + \frac{z_2^2 -2z_1z_2\lambda -1}{1-\lambda^2} \bigg] -
    \frac{1}{2\Omega_{22}}\left( \frac{\alpha_1\alpha_2\lambda\tau}{\sqrt{1+\alpha^2_*}} +\alpha_2 z_2\right)\zeta_1(t),\\
%
\frac{\partial\ell}{\partial \alpha_{1}}  & = & 
    \left[\frac{(\alpha_1+\lambda\alpha_2)\tau}{\sqrt{1+\alpha^2_*}}+z_1 \right] 
    \zeta_1(t), \\    
\frac{\partial\ell}{\partial \alpha_{2}}  & = & 
    \left[\frac{(\alpha_2+\lambda\alpha_1)\tau}{\sqrt{1+\alpha^2_*}}+z_2 \right] 
    \zeta_1(t), \\
  \frac{\partial\ell}{\partial \tau}  & = & 
    \sqrt{1+\alpha^2_*} \, \zeta_1(t)  - \zeta_1(\tau).
\end{eqnarray*}

The above algebraic expressions have been compared with the numerical derivatives computed,
within the \Rlang\ environment  \citep{man:rrr}, 
using the function \texttt{grad} of the package \texttt{numDeriv} \citep{man:numderiv}.
These checks have been repeated for numerous combinations of the parameters, and the
vector $y$; in all these comparisons, matching values have been obtained.  
 
The likelihood equations $\ud \ell(\theta)/\ud \theta=0$ do not lend themselves 
to an explicit solution, which must then be searched for by numerical methods.


\subsection{The observed Fisher information matrix}
\label{sec:info-obs}


On defining $\zeta_{2}(\cdot)$ and $t$ as in (\ref{eq:zeta2}) and (\ref{eq:t}),
the second order partial derivatives of the log-likelihood (\ref{eq:loglikelihood})
can be written as follows. 
With some extension of the terminology,
we shall refer to the matrix so constructed as the `observed information matrix',
even if the matrix is not necessarily positive definite, when it is evaluated at a 
parameter point which is not the maximum likelihood estimate.

For notational simplicity, the generic entry $j_{rs}(\theta)$ of the hessian matrix
with reversed sign, $j(\theta)$, is written  as $j_{rs}$. 
\begin{eqnarray*}
-j_{11} & = &\frac{\partial^2\ell}{\partial \xi^2_1} = 
  -\frac{1}{\Omega_{11}}\left[ \frac{1}{1-\lambda^2} - \alpha^2_{1} \zeta_{2}(t)  \right], \\
-j_{12}&=&\frac{\partial^2 \ell}{\partial \xi_{1}  \partial \xi_{2} } = 
  \frac{1}{\sqrt{\Omega_{11}\Omega_{22}}} \left[ \frac{\lambda}{1-\lambda^2} + \alpha_1\alpha_2 \zeta_{2}(t)    \right], \\
-j_{13}&=&\frac{\partial^2 \ell}{\partial \xi_{1}  \partial \Omega_{11} } = 
    \frac{\lambda z_2-z_1}{(1-\lambda^{2})^{2}\Omega^{3/2}_{11}}
    +\frac{\alpha_1}{2\Omega^{3/2}_{11}}    
    \left( \frac{\alpha_1\alpha_2\lambda\tau}{\sqrt{1+\alpha^2_*}} +\alpha_1 z_1\right)\zeta_2(t) + 
    \frac{\alpha_1}{2\Omega^{3/2}_{11}}\zeta_1(t),  \\
-j_{14}&=&\frac{\partial^2 \ell}{\partial \xi_{1}  \partial \Omega_{12} } = 
  -\frac{2\lambda(\lambda z_2-z_1)} {(1-\lambda^2)^2 \Omega_{11}\sqrt{\Omega_{22}}} - 
  \frac{z_2}{(1-\lambda^2) \Omega_{11}\sqrt{\Omega_{22}}} -
  \frac{\alpha^2_1\alpha_2\tau}{\Omega_{11}\sqrt{\Omega_{22}}\sqrt{1+\alpha^2_*}} \zeta_2(t)  ,\\
%
-j_{15}& = &\frac{\partial^2 \ell}{\partial \xi_{1}  \partial \Omega_{22} } =
  \frac{\lambda(z_2 - z_1\lambda)}{(1-\lambda^2)^2\Omega_{22}\sqrt{\Omega_{11}}} +
  \frac{\alpha_1}{2\Omega_{22}\sqrt{\Omega_{11}}}
  \left( \frac{\alpha_1\alpha_2\lambda\tau}{\sqrt{1+\alpha^2_*}} + \alpha_2 z_2\right)\zeta_2(t), \\
-j_{16}& = &\frac{\partial^2 \ell}{\partial \xi_{1}  \partial \alpha_1 } =
  -\frac{\alpha_1}{\sqrt{\Omega_{11}}}\left[\frac{(\alpha_1+\lambda\alpha_2)\tau}{\sqrt{1+\alpha^2_*}}+z_1                  \right]  \zeta_2(t)   -\frac{1}{\sqrt{\Omega_{11}}}       \zeta_1(t),\\
-j_{17}& = &\frac{\partial^2 \ell}{\partial \xi_{1}  \partial \alpha_2 } =
  -\frac{\alpha_1}{\sqrt{\Omega_{11}}}  \left[\frac{(\alpha_2+\lambda\alpha_1)\tau}{\sqrt{1+\alpha^2_*}}+z_2 \right]          \zeta_2(t), \\
-j_{18}& = &\frac{\partial^2 \ell}{\partial \xi_{1}  \partial \tau } =
  -\frac{\alpha_1\sqrt{1+\alpha^2_*}}{\sqrt{\Omega_{11}}}\zeta_2(t),
\end{eqnarray*}
\begin{eqnarray*}
-j_{22}& = &\frac{\partial^2\ell}{\partial \xi^2_2} =
  -\frac{1}{\Omega_{22}}\left[ \frac{1}{1-\lambda^2} - \alpha^2_{2} \zeta_{2}(t)  \right], \\
-j_{23}& = &\frac{\partial^2 \ell}{\partial \xi_{2}  \partial \Omega_{11} } =
  \frac{\lambda(z_1 - z_2\lambda)}{(1-\lambda^2)^2\Omega_{11}\sqrt{\Omega_{22}}} +
  \frac{\alpha_2}{2\Omega_{11}\sqrt{\Omega_{22}}}
  \left( \frac{\alpha_1\alpha_2\lambda\tau}{\sqrt{1+\alpha^2_*}} +\alpha_1 z_1\right)\zeta_2(t), \\
-j_{24}& = &\frac{\partial^2 \ell}{\partial \xi_{2}  \partial \Omega_{12} } =
    -\frac{2\lambda(\lambda z_1-z_2)} {(1-\lambda^2)^2 \Omega_{22}\sqrt{\Omega_{11}}} - 
  \frac{z_1}{(1-\lambda^2) \Omega_{22}\sqrt{\Omega_{11}}} -
  \frac{\alpha^2_2\alpha_1\tau}{\Omega_{22}\sqrt{\Omega_{11}}\sqrt{1+\alpha^2_*}} \zeta_2(t)  ,\\
-j_{25}& = &\frac{\partial^2 \ell}{\partial \xi_{2}  \partial \Omega_{22} } =
  \frac{\lambda z_1-z_2}{(1-\lambda^{2})^{2}\Omega^{3/2}_{22}}  + 
  \frac{\alpha_2}{2\Omega^{3/2}_{22}}   
  \left( \frac{\alpha_1\alpha_2\lambda\tau}{\sqrt{1+\alpha^2_*}} +\alpha_2 z_2\right)\zeta_2(t) + 
  \frac{\alpha_2}{2\Omega^{3/2}_{22}}\zeta_1(t),  \\
-j_{26}& = &\frac{\partial^2 \ell}{\partial \xi_{2}  \partial \alpha_1 } =
  -\frac{\alpha_2}{\sqrt{\Omega_{22}}}  \left[\frac{(\alpha_1+\lambda\alpha_2)\tau}{\sqrt{1+\alpha^2_*}}+z_1 \right]          \zeta_2(t), \\
-j_{27}& = &\frac{\partial^2 \ell}{\partial \xi_{2}  \partial \alpha_2 } =
  -\frac{\alpha_2}{\sqrt{\Omega_{22}}}\left[\frac{(\alpha_2+\lambda\alpha_1)\tau}{\sqrt{1+\alpha^2_*}}+z_2                  \right]  \zeta_2(t)-\frac{1}{\sqrt{\Omega_{22}}}\zeta_1(t),\\
-j_{28}& = &\frac{\partial^2 \ell}{\partial \xi_{2}  \partial \tau } =
  -\frac{\alpha_2\sqrt{1+\alpha^2_*}}{\sqrt{\Omega_{22}}}\zeta_2(t),
\end{eqnarray*}
\begin{eqnarray*}
-j_{33}&=&\frac{\partial^2 \ell}{\partial\Omega^2_{11} }  = 
  \frac{\lambda^2-z^2_1+2z_1z_2\lambda}{\Omega^2_{11}(1-\lambda^{2})} +                                   \frac{4\lambda^3z_1z_2-2\lambda^2z^2_1-\lambda^2z^2_2}{\Omega^2_{11}(1-\lambda^{2})^{2}} -
  \frac{\lambda^4(z^2_1+z^2_2-2z_1z_2\lambda)}{\Omega^2_{11}(1-\lambda^{2})^{3}} +\\
    && \phantom{ \frac{\partial^2 \ell}{\partial \Omega^2_{11} } }    
  +\frac{1}{2\Omega^2_{11}} + \frac{\lambda^4}{2\Omega^2_{11}(1-\lambda^{2})^{2}}+
  \frac{1}{4\Omega^2_{11}} 
  \left( \frac{3\alpha_1\alpha_2\tau\lambda}{\sqrt{1+\alpha^2_*}} 
  - \frac{\alpha^2_1\alpha^2_2\tau\lambda^2}{(1+\alpha^2_*)^{3/2}} 
  + 3\alpha_1 z_1 \right) \zeta_1(t) + \\
    && \phantom{ \frac{\partial^2 \ell}{\partial \Omega^2_{11} } }  
  +\frac{1}{4\Omega^2_{11}} \left( \frac{\alpha_1\alpha_2\lambda\tau}{\sqrt{1+\alpha^2_*}} +\alpha_1 z_1            \right)^2 \zeta_2(t),\\
-j_{34}&=&\frac{\partial^2 \ell}{\partial \Omega_{11} \partial \Omega_{12} } =
  -\frac{\lambda+z_1z_2}{(1-\lambda^{2}) \Omega^{3/2}_{11}\sqrt{\Omega_{22}}} + 
  \frac{2\lambda z^2_1 + \lambda z^2_2 - 5\lambda^2 z_1z_2 - \lambda^3}{(1-\lambda^{2})^2                     \Omega^{3/2}_{11}\sqrt{\Omega_{22}}} +
  \frac{2\lambda^3(z^2_1+ z^2_2 -2z_1z_2\lambda)} {(1-\lambda^{2})^3\Omega^{3/2}_{11}\sqrt{\Omega_{22}}} +\\
    && \phantom{ \frac{\partial^2 \ell}{\partial \Omega^2_{11} } }    
  +\left( \frac{\alpha^2_1\alpha^2_2\tau\lambda}{2\Omega^{3/2}_{11}\sqrt{\Omega_{22}} (1+\alpha^2_*)^{3/2} } -      \frac{\alpha_1\alpha_2\tau}{2\Omega^{3/2}_{11}\sqrt{\Omega_{22}}\sqrt{1+\alpha^2_*}}\right) \zeta_1(t)  +\\ 
    && \phantom{ \frac{\partial^2 \ell}{\partial \Omega^2_{11} } }
  -\frac{\alpha_1\alpha_2\tau}{2\Omega^{3/2}_{11}\sqrt{\Omega_{22}}\sqrt{1+\alpha^2_*}}
   \left( \frac{\alpha_1\alpha_2\lambda\tau}{\sqrt{1+\alpha^2_*}} + \alpha_1 z_1 \right) \zeta_2(t) ,\\
  \\
-j_{35}&=&\frac{\partial^2 \ell}{\partial \Omega_{11}\partial \Omega_{22} } =
  \frac{\lambda^4 + 6\lambda^3z_1z_2-2\lambda^2z^2_1-2\lambda^2z^2_2}{2\Omega_{11}\Omega_{22}             (1-\lambda^{2})^{2}} +
  \frac{\lambda^2}{2\Omega_{11}\Omega_{22}(1-\lambda^{2})} +
  \frac{\lambda z_1z_2}{\Omega_{11}\Omega_{22}(1-\lambda^2)} + \\
    && \phantom{ \frac{\partial^2 \ell}{\partial \Omega^2_{11} } }    
  -\frac{\lambda^4(z^2_1+z^2_2-2z_1z_2\lambda)}{\Omega_{11}\Omega_{22}(1-\lambda^{2})^{3}} + 
  \frac{\alpha_1\alpha_2\lambda\tau}{4\Omega_{11}\Omega_{22}\sqrt{1+\alpha^2_*}}    
  \left( 1 - \frac{\alpha_1\alpha_2\lambda}{1+\alpha^2_*}\right) \zeta_1(t) + \\
    && \phantom{ \frac{\partial^2 \ell}{\partial \Omega^2_{11} } }  
  +\frac{1}{4\Omega_{11}\Omega_{22}} 
  \left( \frac{\alpha_1\alpha_2\lambda\tau}{\sqrt{1+\alpha^2_*} } + \alpha_1 z_1 \right) 
  \left( \frac{\alpha_1\alpha_2\lambda\tau}{\sqrt{1+\alpha^2_*} } + \alpha_2 z_2 \right) \zeta_2(t),\\
-j_{36}&=&\frac{\partial^2 \ell}{\partial \Omega_{11}\partial \alpha_1 } =
  \frac{1}{2\Omega_{11}}\left( \frac{\alpha_1\alpha_2\lambda(\alpha_2\lambda+\alpha_1)\tau}                     {(1+\alpha^2_*)^{3/2}} -\frac{\alpha_2\lambda\tau}{\sqrt{1+\alpha^2_*}}  - z_1 \right) \zeta_1(t)+\\
    && \phantom{ \frac{\partial^2 \ell}{\partial \Omega^2_{11} } }    
  -\frac{1}{2\Omega_{11}} \left( \frac{\alpha_1\alpha_2\lambda\tau}{\sqrt{1+\alpha^2_*}} + \alpha_1 z_1 \right)         \left( \frac{(\alpha_2\lambda+\alpha_1)\tau}{\sqrt{1+\alpha^2_*}} + z_1 \right) \zeta_2(t),\\ 
-j_{37}&=&\frac{\partial^2 \ell}{\partial \Omega_{11}\partial\alpha_2 } = 
    \frac{1}{2\Omega_{11}}\left(  \frac{\alpha_1\alpha_2\lambda(\alpha_1\lambda+\alpha_2)\tau}                      {(1+\alpha^2_*)^{3/2}} -\frac{\alpha_1\lambda\tau}{\sqrt{1+\alpha^2_*}}\right) \zeta_1(t)+\\
    && \phantom{ \frac{\partial^2 \ell}{\partial \Omega^2_{11} } }    
    -\frac{1}{2\Omega_{11}} \left( \frac{\alpha_1\alpha_2\lambda\tau}{\sqrt{1+\alpha^2_*}} + \alpha_1 z_1 \right)         \left( \frac{(\alpha_2+\alpha_1\lambda)\tau}{\sqrt{1+\alpha^2_*}} + z_2 \right) \zeta_2(t),\\   
-j_{38}&=&\frac{\partial^2 \ell}{\partial \Omega_{11}\partial\tau} =
   - \frac{\alpha_1\alpha_2\lambda}{2\Omega_{11}\sqrt{1+\alpha^2_*}} \zeta_1(t) -
   \frac{\sqrt{1+\alpha^2_*} }{2\Omega_{11}} \left( \frac{\alpha_1\alpha_2\lambda\tau}{\sqrt{1+\alpha^2_*}} +         \alpha_1 z_1 \right)  \zeta_2(t),
\end{eqnarray*}
\begin{eqnarray*}
-j_{44}& = &\frac{\partial^2 \ell}{\partial \Omega^2_{12} } =
  \frac{1}{(1-\lambda^2)\Omega_{11}\Omega_{22}} + 
  \frac{6\lambda z_1 z_2 -z^2_1 - z^2_2 +2\lambda^2  }  {(1-\lambda^2)^2\Omega_{11}\Omega_{22}} 
  -\frac{4\lambda^2(z^2_1 + z^2_2 -2 z_1 z_2\lambda)}{(1-\lambda^2)^3\Omega_{11}\Omega_{22}} +\\
    && \phantom{ \frac{\partial^2 \ell}{\partial \Omega^2_{11} } }    
  +\frac{\alpha^2_1\alpha^2_2\tau}{\Omega_{11}\Omega_{22} (1+\alpha^2_*)} \left( \tau \zeta_2(t) -
  \frac{\zeta_1(t)}{\sqrt{1+\alpha^2_*}} \right) ,  \\  
-j_{45}& = &\frac{\partial^2 \ell}{\partial \Omega_{12}\partial \Omega_{22}  } =
    -\frac{\lambda+z_1z_2}{(1-\lambda^{2})\Omega^{3/2}_{22}\sqrt{\Omega_{11}}} + 
  \frac{2\lambda z^2_2 + \lambda z^2_1 - 5\lambda^2 z_1z_2 - \lambda^3}{(1-\lambda^{2})^2                     \Omega^{3/2}_{22}\sqrt{\Omega_{11}}} +\\
    && \phantom{ \frac{\partial^2 \ell}{\partial \Omega^2_{11} } }    
  +\left( \frac{\alpha^2_1\alpha^2_2\tau\lambda}{2\Omega^{3/2}_{22}\sqrt{\Omega_{11}} (1+\alpha^2_*)^{3/2} } -      \frac{\alpha_1\alpha_2\tau}{2\Omega^{3/2}_{22}\sqrt{\Omega_{11}}\sqrt{1+\alpha^2_*}}\right) \zeta_1(t)  +\\ 
    && \phantom{ \frac{\partial^2 \ell}{\partial \Omega^2_{11} } }
  -\frac{\alpha_1\alpha_2\tau}{2\Omega^{3/2}_{22}\sqrt{\Omega_{11}}\sqrt{1+\alpha^2_*}}
   \left( \frac{\alpha_1\alpha_2\lambda\tau}{\sqrt{1+\alpha^2_*}} + \alpha_2 z_2 \right) \zeta_2(t)+\\
    && \phantom{ \frac{\partial^2 \ell}{\partial \Omega^2_{11} } }
+\frac{2\lambda^3(z^2_1+ z^2_2 -2z_1z_2\lambda)} {(1-\lambda^{2})^3\Omega^{3/2}_{22}\sqrt{\Omega_{11}}},\\    
-j_{46}& = &\frac{\partial^2 \ell}{\partial \Omega_{12}\partial\alpha_1 } =
  \frac{\alpha_2\tau}{\sqrt{\Omega_{11}\Omega_{22}} \sqrt{1+\alpha^2_*}} 
  \left( 1- \frac{\alpha_1(\alpha_2\lambda+\alpha_1)}{1+\alpha^2_*}\right)\zeta_1(t)+\\
    && \phantom{ \frac{\partial^2 \ell}{\partial \Omega^2_{11} } }    
  +\frac{\alpha_1\alpha_2\tau}{\sqrt{\Omega_{11}\Omega_{22}}\sqrt{1+\alpha^2_*}}
  \left( \frac{(\alpha_1+\lambda\alpha_2)\tau}{\sqrt{1+\alpha^2_*}} +z_1 \right)\zeta_2(t),\\
-j_{47}& = &\frac{\partial^2 \ell}{\partial \Omega_{12}\partial\alpha_2 } =
  \frac{\alpha_1\tau}{\sqrt{\Omega_{11}\Omega_{22}} \sqrt{1+\alpha^2_*}} 
  \left( 1- \frac{\alpha_2(\alpha_1\lambda+\alpha_2)}{1+\alpha^2_*}\right)\zeta_1(t)+\\
    && \phantom{ \frac{\partial^2 \ell}{\partial \Omega^2_{11} } }    
  +\frac{\alpha_1\alpha_2\tau}{\sqrt{\Omega_{11}\Omega_{22}}\sqrt{1+\alpha^2_*}}
  \left( \frac{(\alpha_2+\lambda\alpha_1)\tau}{\sqrt{1+\alpha^2_*}} +z_2 \right)\zeta_2(t),\\
-j_{48}& = &\frac{\partial^2 \ell}{\partial \Omega_{12}\partial\tau} =
  \frac{\alpha_1\alpha_2}{\sqrt{\Omega_{11}\Omega_{22}}}
  \left( \frac{\zeta_1(t)}{\sqrt{1+\alpha^2_*}} +\tau\zeta_2(t) \right),
\end{eqnarray*}
\begin{eqnarray*}
-j_{55}& = &\frac{\partial^2 \ell}{\partial \Omega^2_{22} } =
  \frac{\lambda^2-z^2_2+2z_1z_2\lambda}{\Omega^2_{22}(1-\lambda^{2})} +                                   \frac{4\lambda^3z_1z_2-2\lambda^2z^2_2-\lambda^2z^2_1}{\Omega^2_{22}(1-\lambda^{2})^{2}} -
  \frac{\lambda^4(z^2_1+z^2_2-2z_1z_2\lambda)}{\Omega^2_{22}(1-\lambda^{2})^{3}} +\\
    && \phantom{ \frac{\partial^2 \ell}{\partial \Omega^2_{11} } }    
  +\frac{1}{2\Omega^2_{22}} + \frac{\lambda^4}{2\Omega^2_{22}(1-\lambda^{2})^{2}}+
  \frac{1}{4\Omega^2_{22}} 
  \left( \frac{3\alpha_1\alpha_2\tau\lambda}{\sqrt{1+\alpha^2_*}} 
  - \frac{\alpha^2_1\alpha^2_2\tau\lambda^2}{(1+\alpha^2_*)^{3/2}} 
  + 3\alpha_2 z_2 \right) \zeta_1(t) + \\
    && \phantom{ \frac{\partial^2 \ell}{\partial \Omega^2_{11} } }  
  +\frac{1}{4\Omega^2_{22}} \left( \frac{\alpha_1\alpha_2\lambda\tau}{\sqrt{1+\alpha^2_*}} +\alpha_2 z_2            \right)^2 \zeta_2(t),\\
-j_{56}& = &\frac{\partial^2 \ell}{\partial \Omega_{22}\partial\alpha_1 } =
  \frac{1}{2\Omega_{22}}\left(  \frac{\alpha_1\alpha_2\lambda(\alpha_2\lambda+\alpha_1)\tau}                      {(1+\alpha^2_*)^{3/2}} -\frac{\alpha_2\lambda\tau}{\sqrt{1+\alpha^2_*}}\right) \zeta_1(t)+\\
    && \phantom{ \frac{\partial^2 \ell}{\partial \Omega^2_{11} } }    
  -\frac{1}{2\Omega_{22}} \left( \frac{\alpha_1\alpha_2\lambda\tau}{\sqrt{1+\alpha^2_*}} + \alpha_2 z_2 \right)         \left( \frac{(\alpha_1+\alpha_2\lambda)\tau}{\sqrt{1+\alpha^2_*}} + z_1 \right) \zeta_2(t),\\     
-j_{57}& = &\frac{\partial^2 \ell}{\partial \Omega_{22}\partial \alpha_2 } =
  \frac{1}{2\Omega_{22}}\left( \frac{\alpha_1\alpha_2\lambda(\alpha_1\lambda+\alpha_2)\tau}                     {(1+\alpha^2_*)^{3/2}} -\frac{\alpha_1\lambda\tau}{\sqrt{1+\alpha^2_*}}  - z_2 \right) \zeta_1(t)+\\
    && \phantom{ \frac{\partial^2 \ell}{\partial \Omega^2_{11} } }    
  -\frac{1}{2\Omega_{22}} \left( \frac{\alpha_1\alpha_2\lambda\tau}{\sqrt{1+\alpha^2_*}} + \alpha_2 z_2 \right)         \left( \frac{(\alpha_1\lambda+\alpha_2)\tau}{\sqrt{1+\alpha^2_*}} + z_2 \right) \zeta_2(t),\\ 
-j_{58}& = &\frac{\partial^2 \ell}{\partial \Omega_{22}\partial\tau} =
   - \frac{\alpha_1\alpha_2\lambda}{2\Omega_{22}\sqrt{1+\alpha^2_*}} \zeta_1(t) -
   \frac{\sqrt{1+\alpha^2_*} }{2\Omega_{22}} \left( \frac{\alpha_1\alpha_2\lambda\tau}{\sqrt{1+\alpha^2_*}} +         \alpha_2 z_2 \right)  \zeta_2(t),
\end{eqnarray*}
\begin{eqnarray*}
-j_{66}& = &\frac{\partial^2 \ell}{\partial \alpha^2_1 } =
  \left( \frac{\tau}{\sqrt{1+\alpha^2_*}} - \frac{(\alpha_2\lambda+\alpha_1)^2\tau}{(1+\alpha^2_*)^{3/2}} \right) \zeta_1(t) + \left(\frac{(\alpha_1+\lambda\alpha_2)\tau}{\sqrt{1+\alpha^2_*}}+z_1 \right)^2 
    \zeta_2(t),\\
-j_{67}& = &\frac{\partial^2 \ell}{\partial \alpha_1\partial\alpha_2} =
  \left( \frac{\lambda\tau}{\sqrt{1+\alpha^2_*}} - \frac{(\alpha_2+\lambda\alpha_1)(\alpha_1+\lambda\alpha_2)\tau}{(1+\alpha^2_*)^{3/2}} \right) \zeta_1(t) +\\
    && \phantom{ \frac{\partial^2 \ell}{\partial \Omega^2_{11} } }    
  +\left(\frac{(\alpha_1+\lambda\alpha_2)\tau}{\sqrt{1+\alpha^2_*}}+z_1 \right) 
  \left(\frac{(\alpha_2+\lambda\alpha_1)\tau}{\sqrt{1+\alpha^2_*}}+z_2 \right) \zeta_2(t),\\
-j_{68}& = &\frac{\partial^2 \ell}{\partial \alpha_1\partial\tau} =
  \frac{\alpha_1+\lambda\alpha_2}{\sqrt{1+\alpha^2_*}} \zeta_1(t) +   \left(\frac{(\alpha_1+\lambda\alpha_2)\tau}{\sqrt{1+\alpha^2_*}}+z_1 \right) \sqrt{1+\alpha^2_*} \zeta_2(t),
\end{eqnarray*}
%
\begin{eqnarray*}
-j_{77}& = &\frac{\partial^2 \ell}{\partial \alpha^2_2 } =
  \left( \frac{\tau}{\sqrt{1+\alpha^2_*}} - \frac{(\alpha_1\lambda+\alpha_2)^2\tau}{(1+\alpha^2_*)^{3/2}} \right) \zeta_1(t) + \left(\frac{(\alpha_2+\lambda\alpha_1)\tau}{\sqrt{1+\alpha^2_*}}+z_2 \right)^2 
    \zeta_2(t),\\
-j_{78}& = &\frac{\partial^2 \ell}{\partial \alpha_2\partial\tau} =
  \frac{\alpha_2+\lambda\alpha_1}{\sqrt{1+\alpha^2_*}} \zeta_1(t) +   \left(\frac{(\alpha_2+\lambda\alpha_1)\tau}{\sqrt{1+\alpha^2_*}}+z_2 \right) \sqrt{1+\alpha^2_*} \zeta_2(t) ,
\end{eqnarray*} 
\begin{eqnarray*}
-j_{88}& = &\frac{\partial^2 \ell}{\partial \tau^2 } =
  (1+\alpha^2_*)\zeta_2(t)-\zeta_2(\tau).
\end{eqnarray*}
Evaluation of $j(\theta)$ at the MLE point $\hat\theta$ yields the 
observed information matrix, $j(\hat\theta)$.

 Some of these algebraic expressions have been initially obtained  using the
symbolic manipulation program \emph{Maxima} \citep{c.elmohamed}, and then
manually refined.
The final algebraic expressions have been checked against the numerical
computed derivatives obtained from the function \texttt{hessian} of  the
package \texttt{numDeriv} \citep{man:numderiv} of the computing
environment \Rlang. 
Similarly to the evaluation of the first derivatives, these checks have
been performed for various combinations of the parameter values and
of the vector $y$. In all the cases, the values  obtained by the two methods
were essentially coincident, in the sense that the absolute values of
the differences were always smaller than the selected tolerance level of 
the \texttt{hessian}; hence, they can be interpreted as numerical errors.


\subsection{The expected information matrix} 
\label{sec:info-exp}


We first recall some well-known results to be used later on. 
They are available from various sources, for instance 
(A) \citet[p.\,279]{azzalini:1996}, 
(B) \citet[p.\,458]{book:mardia}; 
(C) \citet[p.\,288]{azzalini:1996}.

\begin{lemma}[A: the Sherman-Morrison formula]
  \label{lemma:inversionlemma}
  For a square  $n\times n$  matrix $A$ and  $n \times 1$ column vectors  $b,\,d$, 
  we have
  \begin{equation}
    (A+bd^T)^{-1} = A^{-1}-\frac{1}{1+d^T A^{-1}b}A^{-1}bd^TA^{-1}
  \end{equation}
  provided the inverse matrices exist.
\end{lemma}

\begin{lemma}[B]
  \label{lemma:determinantA}
  If $A$ is an invertible  $n\times n$ square matrix, and $b$ is a $n \times 1$
  column vector, then
  \begin{equation}
    |A+bb^T| = |A|\: (1+b^TA^{-1}b).
  \end{equation}
\end{lemma}
\begin{lemma}[C: completing the square]
\label{lemma:compquad}
If $A$ is a $k\times k$ symmetric positive definite matrix, 
and $b$ is a $k \times 1$ vector, then
\begin{equation}
\int_{\mathds{R}^{k}} \frac{1}{(2 \pi)^{k/2}} \exp\left\{ -\frac{1}{2} \left(y^T A y - 2 b^T y\right)  \right\} \ud y 
    = \frac{\exp\left\{ \frac{1}{2} (b^T A^{-1} b) \right\}} {|A|^{1/2}}.
\end{equation}
where $\ud y$ denotes $\ud y_{1}\cdots \ud y_{k}$.
\end{lemma}

Before proceeding to the computation of the elements of the expected information matrix,
we present two preliminary lemmata; to the best of our knowledge, these results are known. 
 
\begin{lemma}
  \label{lemmastenozero}
  If $Z \sim \ESN_k(0,\bar{\Omega},\alpha,\tau)$ where $\bar\Omega$ is a full-rank correlation 
  matrix, and $\zeta_1(\cdot)$  denotes the function defined by (\ref{eq:zeta1}), then
  \begin{eqnarray}
    \label{eq:lemmastenozero}
    \mean[\zeta_1(\alpha_0 +\alpha^T Z  )]= \frac{\zeta_1(\tau)}{\sqrt{1+\alpha^T\bar{\Omega}\alpha}}   
  \end{eqnarray}  
where $\alpha_{0}=\tau(1+\alpha^T \bar{\Omega}\alpha)^{1/2}$.
\end{lemma}

\begin{proof} Denote by $I$ the expected value on the left of the equality sign in
(\ref{eq:lemmastenozero}) and let $\alpha^2_*$ be as in (\ref{eq:alphaasterisco}).
It turns out that
\begin{eqnarray*}
  I &=& \int_{\mathbbm{R^{k}}} \frac{\phi(\alpha_0 +\alpha^T  z)}{\Phi(\alpha_0 +\alpha^T z)} \times
                  \frac{1}{\Phi(\tau)}  \phi_k(z,\bar{\Omega}) \Phi\left(\alpha_0 + \alpha^T z\right) \ud z \\
  &=& \frac{1}{\Phi({\tau})} \int_{\mathbbm{R^{k}}} \phi_k(z,\bar{\Omega}) \phi(\alpha_0 +\alpha^T z) \ud z\\
    &=& \frac{1}{\Phi({\tau})} \int_{\mathbbm{R^{k}}} \frac{1}{\sqrt{2\pi}} \exp\left\{ -\frac{1}{2}
         (\alpha_0 +\alpha^T z)^2 \right\}  \frac{1}{(2\pi)^\frac{k}{2} |\bar{\Omega}|^\frac{1}{2}} \times
        \exp\left\{ -\frac{1}{2} (z^T \bar{\Omega}^{-1} z) \right\}  \ud z\\ 
    &=& \frac{1}{\Phi({\tau})} \frac{1}{|\bar{\Omega} |^\frac{1}{2}} \frac{1}{\sqrt{2\pi}} \int_{\mathbbm{R^{k}}}                 \frac{1}{(2\pi)^\frac{k}{2} } \exp\left\{ -\frac{1}{2} \left[ (\alpha_0 +\alpha^T z)^2 + (z^T                           \bar{\Omega}^{-1} z) \right] \right\} \ud z.     
\end{eqnarray*}         
Write
\[
  u=\begin{pmatrix}1 \\ z\end{pmatrix},\qquad
  \beta=\begin{pmatrix}\alpha_{0} \\ \alpha\end{pmatrix},\qquad
  \Psi=\begin{pmatrix} 0 & 0 \\ 0 &    \bar{\Omega}^{-1} \end{pmatrix},
\]
\[          
 \beta\beta^T + \Psi 
   = \begin{pmatrix} \alpha^2_{0} & \alpha_{0}\alpha^T  \\ \alpha_{0}\alpha &                       \alpha\alpha^T + \bar{\Omega}^{-1} \end{pmatrix} 
   = \begin{pmatrix} \Sigma^{11} & \Sigma^{12}   \\ 
          \Sigma^{21} & \Sigma^{22}    \end{pmatrix},           
\]
where $\Sigma^{hk}$ is defined as the matching term of the preceding matrix,
so that
\begin{eqnarray*}
  Q  &=&  \left[ (\alpha_0 +\alpha^T z)^2 + (z^T \bar{\Omega}^{-1} z) \right]\\
   &=& \left[ (\alpha_{0} \, \alpha^T) \begin{pmatrix}1\\ z\end{pmatrix}\right]^2+(z^T \bar{\Omega}^{-1}                z)\\ 
   &=& (u^T  \beta)(\beta^T  u)+ (u^T \Psi \, u)\\
   &=& u^T(\beta \beta^T + \Psi)u\\
   &=& u^T \begin{pmatrix} \Sigma^{11} & \Sigma^{12}  \\ \Sigma^{21} & \Sigma^{22} \end{pmatrix} u\\ 
   &=&    \Sigma^{11} + 2\Sigma^{12}z + z^T \Sigma^{22} z.
\end{eqnarray*} 
Further, expand 
\begin{eqnarray}  
      (\Sigma^{22})^{-1} &=& (\bar{\Omega}^{-1}+\alpha\alpha^T )^{-1} = \bar{\Omega}-\frac{1}{1+\alpha^T                          \bar{\Omega}\alpha}\bar{\Omega}\alpha\alpha^T\bar{\Omega},\\
      |\Sigma^{22}| &=& |\bar{\Omega}^{-1}+\alpha\alpha^T| =|\bar{\Omega}^{-1}|                                               (1+\alpha^T\bar{\Omega}\alpha).   
    \end{eqnarray}
From here, it follows that
\begin{eqnarray*}
  I &=&  \frac{1}{\Phi({\tau})} \frac{1}{|\bar{\Omega}|^\frac{1}{2}} 
         \frac{1}{\sqrt{2\pi}} \exp \left\{ -\frac{1}{2} \Sigma^{11} \right\}
       \int_{\mathbbm{R^{k}}} \frac{1}{(2\pi)^\frac{k}{2} } \exp\left\{ -\frac{1}{2} \left[                           z^T\Sigma^{22}z-2(-\Sigma^{12})z \right] \right\} \ud z.
\end{eqnarray*}
By using Lemma~\ref{lemma:compquad}, we obtain
\begin{eqnarray}  
  \label{eq:eqzero} 
  I  &=& \frac{1}{\Phi({\tau})} \frac{1}{|\bar{\Omega} |^{\frac{1}{2}}} \frac{1}{\sqrt{2\pi}}
        \exp \left\{ -\frac{1}{2} \Sigma^{11} \right\} 
        \frac{\exp\left\{\frac{1}{2} \left[ \Sigma^{12}                           ({\Sigma^{22}})^{-1} ({\Sigma^{12}})^T \right] \right\}}{|(\Sigma^{22})|^{\frac{1}{2}}}\\ 
  &=&\frac{1}{\Phi({\tau})} \frac{1}{\sqrt{1+\alpha^T\bar{\Omega}\alpha}} 
         \frac{1}{\sqrt{2\pi}} \exp\left\{-\frac{1}{2} \alpha^{2}_{0}\right\} \notag
      \exp {\left\{\frac{1}{2} [\alpha^{2}_{0} \alpha^{T} (\Sigma^{22})^{-1} \alpha\right\} }  
      \notag \\
  &=&\frac{1}{\Phi({\tau})} \frac{1}{\sqrt{1+\alpha^T\bar{\Omega}\alpha}} 
            \frac{1}{\sqrt{2\pi}} \exp\left\{-\frac{1}{2} \tau^{2}\right\} \notag
     \exp {\left\{\frac{1}{2} [\alpha^{2}_{0} \alpha^{T} (\Sigma^{22})^{-1} \alpha - \tau^{2}   \notag                        \alpha^{T}\bar{\Omega}\alpha]\right\} } \notag \\
   &=& \frac{\phi(\tau)}{\Phi({\tau})} \frac{1}{\sqrt{1+\alpha^T\bar{\Omega}\alpha}} \:
       \exp{\left\{\frac{1}{2} \left[ \tau(1+ \alpha^T\bar{\Omega}\alpha)\alpha^T 
           \left(\bar{\Omega} - \frac{1}{1+\alpha^T \bar{\Omega}\alpha} 
           \bar{\Omega}\alpha\alpha^T\bar{\Omega}\right) \alpha 
          - \tau^2\alpha^T\bar{\Omega}\alpha\right]\right\}}\notag \\     
   &=&\frac{ \zeta_1(\tau) }{\sqrt{1+\alpha^T\bar{\Omega}\alpha}} 
       \exp{\left\{ \frac{1}{2} 
        \left[\tau(1+\alpha^2_{*})\left(\alpha^2_{*} 
        - \frac{\alpha^4_{*}}                                         {1+\alpha^2_{*}}\right)-\tau\alpha^2_{*} \right]\right\}} \notag \\ 
  &=& \frac{ \zeta_1(\tau) }{\sqrt{1+\alpha^T\bar{\Omega}\alpha}} 
       \exp{\left\{ \frac{1}{2}
       \left[\tau\left(\alpha^2_{*} + \alpha^4_{*} - \alpha^4_{*}\right)-\tau\alpha^2_{*} \right]
         \right\}} \notag \\
  &=& \frac{ \zeta_1(\tau) }{\sqrt{1+\alpha^T\bar{\Omega}\alpha}}.\notag
\end{eqnarray}
\end{proof}


The following lemma, initially presented by \citet{franco:2011}, 
includes and generalizes the situation considered in Lemma~\ref{lemmastenozero}.
Thanks to this result, it is possible to compute each term of the expected information matrix.

\begin{lemma}
  \label{lemmasteno}
If $Z \sim \ESN_k(0,\bar{\Omega},\alpha,\tau)$ where $\bar\Omega$ is a full-rank correlation 
matrix, and $\zeta_1(\cdot)$ is the function defined by (\ref{eq:zeta1}), then
\begin{eqnarray}
    \label{eq:lemmasteno}
    \mean [ h(Z)\: \zeta_1(\alpha_0 +\alpha^T Z )] &=&
      \mean[\zeta_1(\alpha_0 +\alpha^T Z  )] ~ \mean[ h(U)]
\end{eqnarray}
where  $\alpha_0$ is given by (\ref{eq:alpha0}),
\begin{eqnarray}
  U &\sim & \N_k(-\tau\delta, \:\bar{\Omega}-\delta\delta^T), \label{eq:U}   \\
  \delta &=&  \frac{1}{(1+\alpha^T\bar{\Omega}\alpha)^{1/2}}\: \bar{\Omega}\alpha,    
  \label{eq:delta_cap3}
\end{eqnarray}
for any function $h(\cdot)$  from $\mathds{R}$ to $\mathds{R}$ such that the integrals
involved exist.
\end{lemma}

\begin{proof} 
Denote by $I$ the expected value on the left of the equality sign in (\ref{eq:lemmasteno}). 
It turns out that
    \begin{eqnarray*}
      I
      &=& \int_{\mathbbm{R^{k}}} h(z) \frac{\phi(\alpha_0 +\alpha^T  z)}{\Phi(\alpha_0 +\alpha^T z)} 
                  \frac{1}{\Phi(\tau)}  \phi_k(z,\bar{\Omega}) \Phi\left(\alpha_0 + \alpha^T z\right) \ud z \\
      &=& \frac{1}{\Phi({\tau})} \frac{1}{|\bar{\Omega}|^\frac{1}{2}} \frac{1}{\sqrt{2\pi}} \exp \left\{ -\frac{1}{2}                 \Sigma^{11} \right\}
      \int_{\mathbbm{R^{k}}} h(z) \frac{1}{(2\pi)^\frac{k}{2} } \exp\left\{ -\frac{1}{2} \underbrace{\left[                     z^T\Sigma^{22}z-2(-\Sigma^{12})z \right]}_{Q_{1}} \right\} \ud z ,  
    \end{eqnarray*}
using algebraic steps similar to those of Lemma \ref{lemmastenozero}. 

In $Q_1$,  we add and subtract the term \emph{$\Sigma^{12} (\Sigma^{22})^{-1} (\Sigma^{12})^T$} 
so to obtain an integrand represented by the function $h(x)$ multiplied by the density 
function of the multivariate normal variable $U$ with mean vector 
$-(\Sigma^{22})^{-1}(\Sigma^{12})^T$ and covariance matrix $(\Sigma^{22})^{-1}$;
hence, we arrived at
\begin{eqnarray*}   
      I     
      &=& \frac{1}{\Phi({\tau})} \frac{|(\Sigma^{22})^{-1}|^{\frac{1}{2}} }{|\bar{\Omega} |^{\frac{1}{2}}} \frac{1}{\sqrt{2\pi}}
      \exp \left\{ -\frac{1}{2} \Sigma^{11} \right\} \exp\left\{\frac{1}{2} \left[ \Sigma^{12} ({\Sigma^{22}})^{-1}             ({\Sigma^{12}})^T   \right] \right\} \times \\
      && \int_{\mathbbm{R^{k}}} h(z) \frac{1}{(2\pi)^\frac{k}{2}  |(\Sigma^{22})^{-1}|^{\frac{1}{2}} }                      \exp\left\{-\frac{1}{2} \left[ z^T\Sigma^{22}z + 2\Sigma^{12}z + \Sigma^{12} ({\Sigma^{22}})^{-1}                   ({\Sigma^{12}})^T \right] \right\}  \ud z \\
      &=& \frac{1}{\Phi({\tau})} \frac{|(\Sigma^{22})^{-1}|^{\frac{1}{2}} }{|\bar{\Omega} |^{\frac{1}{2}}} \frac{1}{\sqrt{2\pi}}
      \exp \left\{ -\frac{1}{2} \Sigma^{11} \right\} \exp\left\{\frac{1}{2} \left[ \Sigma^{12} ({\Sigma^{22}})^{-1}             ({\Sigma^{12}})^T   \right] \right\} \times \mean [h(U)].
    \end{eqnarray*} 
Since in the above expression the terms before   $\mean [h(U)]$ coincide with the result
obtained in equation  (\ref{eq:eqzero}) of the development of Lemma \ref{lemmastenozero},
then we conclude that $I$ is equal to the term on right side of  (\ref{eq:lemmasteno}).
  
Next, let us expand the mean vector  of $U$, taking into account
that $\alpha^2_{*}$ is defined at (\ref{eq:alphaasterisco}) and $\delta$ at (\ref{eq:delta_cap3}): 
    \begin{eqnarray*}
      -(\Sigma^{22})^{-1}(\Sigma^{12})^T 
      &=& -\tau\sqrt{(1+ \alpha^T\bar{\Omega}\alpha)}\left( \bar{\Omega}-\frac{1}{1+\alpha^T                          \bar{\Omega}\alpha}\bar{\Omega}\alpha\alpha^T\bar{\Omega}\right)\alpha \\
      &=& -\tau\sqrt{1+\alpha^2_{*}} \left( \bar{\Omega}\alpha-\frac{1}                                               {\sqrt{1+\alpha^2_{*}}} \bar{\Omega}\alpha\alpha^2_{*} \right) \\
      &=& -\tau \left( \frac {(1+\alpha^2_{*})\bar{\Omega}\alpha - \alpha^2_{*}\bar{\Omega}\alpha}                        {\sqrt{1+\alpha^2_{*}}}\right)\\
      &=& -\tau\frac{1}{(1+\alpha^T\bar{\Omega}\alpha)^{1/2}}\bar{\Omega}\alpha\\
      &=& -\tau\delta 
    \end{eqnarray*}
which coincides with the expression indicated in the statement of the lemma; 
moreover, for the covariance matrix, write
\begin{eqnarray*}
(\Sigma^{22})^{-1} \,=\, (\bar{\Omega}^{-1} + \alpha\alpha^T)^{-1}  \,=\, \bar{\Omega}-\frac{1}{1+\alpha^T                  \bar{\Omega} \alpha}\bar{\Omega}\alpha\alpha^T\bar{\Omega} \,=\, \bar{\Omega}-\delta\delta^T.
\end{eqnarray*}
\end{proof}


Let us confirm that the above Lemma~\ref{lemmasteno} generalizes and
simplifies Lemma~\ref{lemmatony} presented by \citet{tesi:canale},
for the scalar case.
The simplification arises from the fact that it is no longer necessary
to compute the inverse of a function of a \rv\ and then compute its
expected value; now, we only need to compute the expected value of
a \rv\ of which we know the distribution.

Consider a \rv\ $Z\sim \ESN(0,1,\alpha,\tau)$, 
recall $\zeta_1(\cdot)$ defined by (\ref{eq:zeta1}), $\alpha_0=\tau\sqrt{1+\alpha^2}$, 
and let  $f(\cdot)$, $h(\cdot)$ be functions from $\mathds{R}$ to $\mathds{R}$ 
such that the involved integrals exist. 
On using Lemma \ref{lemmastenozero} with $k=1$, we have that
\begin{eqnarray*}
    \mean[\zeta_1(\alpha_0 +\alpha^T Z  )] 
    &=&   \frac{\zeta_1(\tau)}{\sqrt{1+\alpha^T\bar{\Omega}\alpha}} \\
    &=& \frac{\zeta_1(\tau)}{\sqrt{1+\alpha^2}}\\
    &=& \mean [ \zeta_1(\alpha_0 +\alpha Z)].
  \end{eqnarray*}
Let
  \begin{eqnarray*}
    g(z) &=& \alpha \tau + \sqrt{\alpha^2+1}\, z\quad =\quad v \\
    g^{-1}(v) &=&\frac{v-\alpha\tau}{\sqrt{1+\alpha^2}}
  \end{eqnarray*}
  Se $V \sim N(0,1)$ e $\delta$ \`e data da (\ref{eq:delta_cap3}) allora
  \begin{eqnarray*}
    \quad V^*= g^{-1}(V) = \frac{V-\alpha\tau}{\sqrt{1+\alpha^2}} \sim                                            N\left(-\frac{\alpha\tau}{\sqrt{1+\alpha^2}},\frac{1}{1+\alpha^2}\right), \\
    U \sim  N_k(-\tau\delta,\bar{\Omega}-\delta\delta^T) \sim
    N_1 \left(-\frac{\alpha\tau}{\sqrt{1+\alpha^2}},\frac{1}{1+\alpha^2}\right).
  \end{eqnarray*}
Hence, if $f(\cdot)$ and $h(\cdot)$ represent the same function, we obtain
  \begin{eqnarray*}
  \mean[f(V^*)] = \mean [f(g^{-1}(V))] = \mean[h(U)]
  \end{eqnarray*}
and
\begin{eqnarray*}
    \mean [ f(Z) \zeta_1(\alpha_0 +\alpha Z )] 
      &=& \frac{\zeta_1(\tau)}{\sqrt{\alpha^2+1}} \ \mean [f(g^{-1}(V))]                       
      \quad=\quad \mean[\zeta_1(\alpha_0 +\alpha^T Z  )] \ \mean[ h(U)].
\end{eqnarray*}

We now present some results which will be useful for computing the expected values 
of the second-order partial derivatives of the log-likelihood, using
the cumulant generating function (\ref{eq:fngen}) and Lemma~\ref{lemmasteno}.
Se $Z\sim\ESN_2(0,\bar{\Omega},\alpha,\tau)$, then:
           \begin{eqnarray}
\label{eq:E_noti1}
\mean\left[ Z_1\right] &=& \zeta_1(\tau)\delta_1 \\ 
\mean\left[ Z_2\right] &=& \zeta_1(\tau)\delta_2 \\ 
\mean\left[ Z_1Z_2\right] &=& \lambda+\delta_1\delta_2[\zeta_2(\tau)+\zeta_1(\tau)^2]\\ 
\mean\left[ Z^2_1\right] &=& 1+\delta^2_1[\zeta_2(\tau)+\zeta_1(\tau)^2]\\
\mean\left[ Z^2_2\right] &=& 1+\delta^2_2[\zeta_2(\tau)+\zeta_1(\tau)^2]
\end{eqnarray}

Let $U$ be as in (\ref{eq:U}), whose mean value is the $k\times1$ vector $-\tau\delta$,
the postive definite covariance matrix is $V=\bar{\Omega}-\delta\delta^T$;
also let $\lambda$ be as in (\ref{eq:lambda}) and $\alpha^2_{*}=\alpha^T\bar{\Omega}\alpha$. 
Specifically, for $k=2$ we have:
\begin{eqnarray}
\label{eq:alphazerostar}
\alpha^2_{*} &=& \alpha^2_1 + \alpha^2_2 + 2\alpha_1\alpha_2\lambda,\notag\\
\alpha_0 &=& \tau\sqrt{1+\alpha^2_*},
\end{eqnarray}
the components of the mean vector
\begin{eqnarray*}
-\tau\delta_1 &=& -\tau \frac{1}{\sqrt{1+\alpha^2_*}} (\alpha_1+\alpha_2\lambda) ,\\
-\tau\delta_2 &=& -\tau \frac{1}{\sqrt{1+\alpha^2_*}} (\alpha_2+\alpha_1\lambda) ,
\end{eqnarray*}
and those of the covariance matrix
\[          
  V = \begin{pmatrix} v_{11} & v_{12}  \\ v_{21} & v_{22} \end{pmatrix}
\]
where 
\begin{eqnarray*}
v_{11}&=& \frac{1+\alpha_2(1-\lambda^2)}{1+\alpha_*^2},\\
v_{22}&=&\frac{1+\alpha_1(1-\lambda^2)}{1+\alpha_*^2},\\
v_{12}&=& \frac{ \lambda- \alpha_1\alpha_2(1-\lambda^2)}  {1+\alpha_*^2}.
\end{eqnarray*}

\noindent 
Define by
\begin{eqnarray}
   \label{eq:T}
   T= \alpha_0 +\alpha^T Z = \alpha_0+\alpha_1Z_1+\alpha_2Z_2,
\end{eqnarray}
the \rv\ associated to the statistic (\ref{eq:t}), and let $\alpha_0$ as in (\ref{eq:alphazerostar}).
Then
\begin{eqnarray}  
  \mean [\zeta_1(T)] &=& \frac{\zeta_1(\tau)}{\sqrt{1+\alpha_*^2}} ,\label{eq:valorinoti1}\\
  \mean [Z_1 \zeta_1(T)] &=& -\tau\delta_1 \frac{\zeta_1(\tau)}{\sqrt{1+\alpha_*^2}} ,\\
  \mean [Z_2 \zeta_1(T)] &=& -\tau\delta_2 \frac{\zeta_1(\tau)}{\sqrt{1+\alpha_*^2}} ,\\
  \mean [Z_1^2 \zeta_1(T)] &=& (\tau^2\delta^2_1+v_{11}) \frac{\zeta_1(\tau)}{\sqrt{1+\alpha_*^2}} ,\\
  \mean [Z_2^2 \zeta_1(T)] &=& (\tau^2\delta^2_2+v_{22}) \frac{\zeta_1(\tau)}{\sqrt{1+\alpha_*^2}} ,\\
  \mean [T \zeta_1(T)] &=& (\alpha_0-\tau\alpha^T\delta) 
        \frac{\zeta_1(\tau)}{\sqrt{1+\alpha_*^2}} , \\
  \mean [Z_1 T \zeta_1(T)] &=& [-\alpha_0\tau\delta_1 + \alpha_1(\tau^2\delta^2_1+v_{11}) +                                         \alpha_2(\tau^2\delta_1\delta_2+v_{12})] \frac{\zeta_1(\tau)}{\sqrt{1+\alpha_*^2}}  , \\  
  \mean [Z_2 T \zeta_1(T)] &=& [-\alpha_0\tau\delta_2 + \alpha_2(\tau^2\delta^2_2+v_{22}) +                                         \alpha_1(\tau^2\delta_1\delta_2+v_{12})] \frac{\zeta_1(\tau)}{\sqrt{1+\alpha_*^2}}  , \\
  \mean [\zeta_2(T)] &=& -(\alpha_0-\tau\alpha^T \delta) 
        \frac{\zeta_1(\tau)}{\sqrt{1+\alpha_*^2}} - \mean[\zeta_1(T)^2 ],\\
  \mean [Z_1\zeta_2(T)] &=& - [-\alpha_0\tau\delta_1 + \alpha_1(\tau^2\delta^2_1+v_{11}) +                                              \alpha_2(\tau^2\delta_1\delta_2+v_{12})] \frac{\zeta_1(\tau)}{\sqrt{1+\alpha_*^2}} - \notag \\
                      && \phantom{  [Z_1 }              
                    \mean[Z_1\zeta_1(T)^2 ],\\
  \mean [Z_2\zeta_2(T)] &=& - [-\alpha_0\tau\delta_2 + \alpha_2(\tau^2\delta^2_2+v_{22}) +                                              \alpha_1(\tau^2\delta_1\delta_2+v_{12})] \frac{\zeta_1(\tau)}{\sqrt{1+\alpha_*^2}} - \notag \\
                      && \phantom{  [Z_1 }  
                    \mean[Z_2\zeta_1(T)^2 ],\\
  \mean [Z^2_1\zeta_2(T)] &=& - \bigg[ \alpha_0(v_{11}+\tau^2\delta^2_1)-\alpha_1\tau\delta_1(\tau^2\delta^2_1+3v_{11})  +\notag\\
                              && \phantom{ [Z_1}
                              \alpha_2\tau\left( \frac{v_{12}}{v_{11}}\delta_1-\delta_2\right)(\tau^2\delta^2_1+v_{11}) - \notag\\
                              && \phantom{ [Z_1}
                              \alpha_2\tau\delta_1\frac{v_{12}}{v_{11}}(\tau^2\delta^2_1+3v_{11}) \bigg]
                              \frac{\zeta_1(\tau)}{\sqrt{1+\alpha_*^2}} - \mean[Z^2_1\zeta_1(T)^2 ],\\      
  \mean [Z^2_2\zeta_2(T)] &=& - \bigg[ \alpha_0(v_{22}+\tau^2\delta^2_2)-\alpha_2\tau\delta_2(\tau^2\delta^2_2+3v_{22})  +\notag\\
                              && \phantom{ [Z_1}
                              \alpha_1\tau\left( \frac{v_{12}}{v_{22}}\delta_2-\delta_1\right)(\tau^2\delta^2_2+v_{22}) - \notag\\
                              && \phantom{ [Z_1}
                              \alpha_1\tau\delta_2\frac{v_{12}}{v_{22}}(\tau^2\delta^2_2+3v_{22}) \bigg]
                              \frac{\zeta_1(\tau)}{\sqrt{1+\alpha_*^2}} - \mean[Z^2_2\zeta_1(T)^2 ],\\                              
\label{eq:E_noti2}                              
\mean [Z_1Z_2\zeta_2(T)] &=&  - \bigg[ \alpha_0(v_{12}+\tau^2\delta_1\delta_2)+
                            \alpha_1\tau\left( \frac{v_{12}}{v_{11}}\delta_1-\delta_2\right)(\tau^2\delta^2_1+v_{11}) - \notag \\
                            && \phantom{  [Z_1 }
                            \alpha_1\tau\delta_1\frac{v_{12}}{v_{11}}(\tau^2\delta^2_1+3v_{11}) + 
                            \alpha_2\tau\left( \frac{v_{12}}{v_{22}}\delta_2-\delta_1\right)(\tau^2\delta^2_2+v_{22}) - \notag \\                                           && \phantom{  [Z_1 }
                            \alpha_2\tau\delta_2\frac{v_{12}}{v_{22}}(\tau^2\delta^2_2+3v_{22}) \bigg] 
                            \frac{\zeta_1(\tau)}{\sqrt{1+\alpha_*^2}} - \mean[Z_1Z_2\zeta_1(T)^2 ].  
\end{eqnarray}

The above expressions involve a small number of expected values which are not
amenable to explicit evaluation, namely 
\begin{eqnarray}
a_{j,p} &=& \mean \left[Z_j^p \zeta_1(T)^2 \right],  \label{eq:akesn} \\
a_{ij} &=& \mean \left[Z_i Z_j\zeta_1(T)^2 \right], \label{eq:aijesn}
\end{eqnarray}
and
\begin{equation}
\label{eq:ak0esn}
a_{0} = \mean \left[\zeta_1(T)^2 \right]
\end{equation}
where $j=1,2$, $i=1,2$ and $p=0,1,2$.
The definition  of these quantities, as well as their notation, is analogous to those
introduced by \cite{art:azza85} and \citet{tesi:canale,canale:2011} for the 
information matrices of the univariate SN and ESN distributions, respectively.

Given the lack of explicit expressions, the computation of (\ref{eq:akesn}) 
to (\ref{eq:ak0esn})  will be  accomplished by numerical integration
using the function \texttt{cuhre} of the \Rlang\ package \texttt{cubature}
by \citet{narasimhan:etal:R-pkg-cubature}.

The expected information matrix, $i(\theta)$, turns out to be
\[
i(\theta)= \left(
\begin{tabular}{cccccccc}
$i_{11}$ & $i_{12}$ & $i_{13}$ & $i_{14}$ & $i_{15}$ & $i_{16}$ & $i_{17}$ & $i_{18}$ \\ 
 & $i_{22}$ & $i_{23}$ & $i_{24}$ & $i_{25}$ & $i_{26}$ & $i_{27}$ & $i_{28}$ \\ 
 &  & $i_{33}$ & $i_{34}$ & $i_{35}$ & $i_{36}$ & $i_{37}$ & $i_{38}$ \\ 
 &  &  & $i_{44}$ & $i_{45}$ & $i_{46}$ & $i_{47}$ & $i_{48}$ \\ 
 &  &  &  & $i_{55}$ & $i_{56}$ & $i_{57}$ & $i_{58}$ \\ 
 &  &  &  &  & $i_{66}$ & $i_{67}$ & $i_{68}$ \\ 
 &  &  &  &  &  & $i_{77}$ & $i_{78}$ \\ 
 &  &  &  &  &  &  & $i_{88}$ \\ 
\end{tabular}
\right)
\]
which is symmetric positive semi-defined when the lower triangle is completed by symmetry.
Its terms are given by:
\begin{eqnarray}
  i_{11}&=&-\mean\left[ \frac{\partial^2\ell}{\partial \xi^2_1}\right]  = 
    \frac{1}{\Omega_{11}}\left[ \frac{1}{1-\lambda^2} - \alpha^2_{1} \mean \left[ \zeta_{2}(T) \right] \right],\\
  i_{12}&=&-\mean\left[\frac{\partial^2 \ell}{\partial \xi_{1}  \partial \xi_{2} }\right] =
     -\frac{1}{\sqrt{\Omega_{11}\Omega_{22}}} \left[ \frac{\lambda}{1-\lambda^2} + \alpha_1\alpha_2 \mean[ \zeta_{2}(T)] \right] ,\\
  i_{13}&=&-\mean\left[ \frac{\partial^2 \ell}{\partial \xi_{1}  \partial \Omega_{11} }\right] 
  = -\frac{\lambda \mean\left[Z_2\right]-\mean\left[Z_1\right]}{(1-\lambda^{2})^{2}\Omega^{3/2}_{11}}   
      -\frac{\alpha_1}{2\Omega^{3/2}_{11}} \mean\left[ \zeta_1(T)\right]+\notag\\
      && \phantom{ \partial^2 \ell }          
       -\frac{\alpha_1}{2\Omega^{3/2}_{11}}   
      \left( \frac{\alpha_1\alpha_2\lambda\tau}{\sqrt{1+\alpha^2_*}}\mean\left[\zeta_2(T)\right]+\alpha_1 \mean\left[ Z_1\zeta_2(T)\right] \right),\\ 
i_{14}&=&-\mean\left[ \frac{\partial^2 \ell}{\partial \xi_{1}  \partial \Omega_{12} } \right]= 
  \frac{2\lambda(\lambda\mean\left[Z_2\right]-\mean\left[Z_1\right] )} {(1-\lambda^2)^2 \Omega_{11}\sqrt{\Omega_{22}}} +
  \frac{\mean\left[Z_2\right]}{(1-\lambda^2) \Omega_{11}\sqrt{\Omega_{22}}} +\notag\\
    && \phantom{ \partial^2 \ell }
  +\frac{\alpha^2_1\alpha_2\tau}{\Omega_{11}\sqrt{\Omega_{22}}\sqrt{1+\alpha^2_*}}  \mean\left[\zeta_2(T)\right],\\
i_{15}&=&-\mean\left[ \frac{\partial^2 \ell}{\partial \xi_{1}  \partial \Omega_{22} } \right] = 
  - \frac{\lambda(\mean\left[Z_2\right]-\mean\left[Z_1\right]\lambda)}{(1-\lambda^2)^2\Omega_{22}\sqrt{\Omega_{11}}} + \notag\\
    && \phantom{ \partial^2 \ell }      
  -\frac{\alpha_1}{2\Omega_{22}\sqrt{\Omega_{11}}}
  \left( \frac{\alpha_1\alpha_2\lambda\tau}{\sqrt{1+\alpha^2_*}}\mean\left[\zeta_2(T)\right] + \alpha_2 \mean\left[ Z_2\zeta_2(T)\right]\right),\\
i_{16}&=&-\mean\left[\frac{\partial^2 \ell}{\partial \xi_{1}  \partial \alpha_1 } \right]  =
  \frac{\alpha_1}{\sqrt{\Omega_{11}}}\left[\frac{(\alpha_1+\lambda\alpha_2)\tau}                                                            {\sqrt{1+\alpha^2_*}}\mean\left[\zeta_2(T)\right]+\mean\left[Z_1\zeta_2(T)\right]\right] +\notag\\
    && \phantom{ \partial^2 \ell }
  +\frac{1}{\sqrt{\Omega_{11}}} \mean\left[\zeta_1(T)\right]  ,\\
i_{17}&=&-\mean\left[\frac{\partial^2 \ell}{\partial \xi_{1}  \partial \alpha_2 } \right] =
  \frac{\alpha_1}{\sqrt{\Omega_{11}}}  \left[\frac{(\alpha_2+\lambda\alpha_1)\tau}{\sqrt{1+\alpha^2_*}}\mean\left[\zeta_2(T)\right]+\mean\left[Z_2\zeta_2(T)\right]\right] , \\
i_{18}&=&-\mean\left[\frac{\partial^2 \ell}{\partial \xi_{1}  \partial \tau } \right] = 
  \frac{\alpha_1\sqrt{1+\alpha^2_*}}{\sqrt{\Omega_{11}}}\mean\left[\zeta_2(T)\right],
\end{eqnarray}  
\begin{eqnarray}
i_{22}&=&-\mean\left[\frac{\partial^2\ell}{\partial \xi^2_2} \right] =
  \frac{1}{\Omega_{22}}\left[ \frac{1}{1-\lambda^2} - \alpha^2_{2} \mean\left[\zeta_2(T)\right]  \right], \\
i_{23}&=&-\mean\left[ \frac{\partial^2 \ell}{\partial \xi_{2}  \partial \Omega_{11} } \right]=
  - \frac{\lambda(\mean\left[Z_1\right]-\mean\left[Z_2\right]\lambda)}{(1-\lambda^2)^2\Omega_{11}\sqrt{\Omega_{22}}} + \notag\\
    && \phantom{ \partial^2 \ell }    
  -\frac{\alpha_2}{2\Omega_{11}\sqrt{\Omega_{22}}}
  \left( \frac{\alpha_1\alpha_2\lambda\tau}{\sqrt{1+\alpha^2_*}}\mean\left[\zeta_2(T)\right] +\alpha_1 \mean\left[ Z_1\zeta_2(T)\right]\right),\\
i_{24}&=&-\mean\left[ \frac{\partial^2 \ell}{\partial \xi_{2}  \partial \Omega_{12} }  \right]=
  \frac{2\lambda(\lambda \mean\left[Z_1\right]-\mean\left[Z_2\right])} {(1-\lambda^2)^2 \Omega_{22}\sqrt{\Omega_{11}}} + 
  \frac{\mean\left[Z_1\right]}{(1-\lambda^2) \Omega_{22}\sqrt{\Omega_{11}}} +\notag\\
    && \phantom{ \partial^2 \ell }
  +\frac{\alpha^2_2\alpha_1\tau}{\Omega_{22}\sqrt{\Omega_{11}}\sqrt{1+\alpha^2_*}} \mean\left[\zeta_2(T)\right],\\
i_{25}&=&-\mean\left[ \frac{\partial^2 \ell}{\partial \xi_{2}  \partial \Omega_{22} } \right] =
  -\frac{\lambda\mean\left[Z_1\right]-\mean\left[Z_2\right]}{(1-\lambda^{2})^{2}\Omega^{3/2}_{22}}  
  -\frac{\alpha_2}{2\Omega^{3/2}_{22}}\mean\left[ \zeta_1(T)\right]+ \notag \\
    && \phantom{ \partial^2 \ell }      
  -\frac{\alpha_2}{2\Omega^{3/2}_{22}}    
  \left( \frac{\alpha_1\alpha_2\lambda\tau}{\sqrt{1+\alpha^2_*}}\mean\left[\zeta_2(T)\right] +\alpha_2 \mean\left[ Z_2\zeta_2(T)\right]\right),\\ 
i_{26}&=&-\mean\left[\frac{\partial^2 \ell}{\partial \xi_{2}  \partial \alpha_1 }\right] =
  \frac{\alpha_2}{\sqrt{\Omega_{22}}}  \left[\frac{(\alpha_1+\lambda\alpha_2)\tau}{\sqrt{1+\alpha^2_*}} \mean\left[\zeta_2(T)\right]
  +\mean\left[Z_1\zeta_2(T)\right] \right], \\
i_{27}&=&-\mean\left[\frac{\partial^2 \ell}{\partial \xi_{2}  \partial \alpha_2 } \right] =
  \frac{\alpha_2}{\sqrt{\Omega_{22}}}\left[\frac{(\alpha_2+\lambda\alpha_1)\tau}{\sqrt{1+\alpha^2_*}}\mean\left[\zeta_2(T)\right]
  +\mean\left[Z_2\zeta_2(T)\right]\right]+ \notag\\
    && \phantom{ \partial^2 \ell }
  +\frac{1}{\sqrt{\Omega_{22}}}\mean\left[\zeta_1(T)\right],\\
i_{28}&=&-\mean\left[\frac{\partial^2 \ell}{\partial \xi_{2}  \partial \tau } \right] =
  \frac{\alpha_2\sqrt{1+\alpha^2_*}}{\sqrt{\Omega_{22}}}\mean\left[\zeta_2(T)\right],
\end{eqnarray}  
\begin{eqnarray}
i_{33}&=&-\mean\left[\frac{\partial^2 \ell}{\partial\Omega^2_{11} } \right] = 
  - \frac{\lambda^2-\mean\left[Z^2_1 \right] +2\mean\left[ Z_1Z_2\right] \lambda}{\Omega^2_{11}(1-\lambda^{2})} 
  -   \frac{4\lambda^3\mean\left[ Z_1Z_2\right] -2\lambda^2 \mean\left[ Z^2_1 \right] -\lambda^2\mean\left[  Z^2_2 \right]}{\Omega^2_{11}(1-\lambda^{2})^{2}} + \notag \\
    && \phantom{ \partial^2 \ell }    
  -\frac{1}{4\Omega^2_{11}} 
  \left( \frac{3\alpha_1\alpha_2\tau\lambda}{\sqrt{1+\alpha^2_*}} \mean\left[\zeta_1(T) \right]
  - \frac{\alpha^2_1\alpha^2_2\tau\lambda^2}{(1+\alpha^2_*)^{3/2}} \mean\left[\zeta_1(T) \right]
  + 3\alpha_1 \mean\left[Z_1\zeta_1(T) \right] \right)  + \notag \\
    && \phantom{ \partial^2 \ell }
  -\frac{1}{4\Omega^2_{11}} \left( \frac{\alpha^2_1\alpha^2_2\lambda^2\tau^2}{1+\alpha^2_*}\mean\left[\zeta_2(T)\right] +
  \alpha^2_1\mean\left[Z^2_1\zeta_2(T) \right]+ 2\frac{\alpha^2_1\alpha_2\lambda\tau}{\sqrt{1+\alpha^2_*}} \mean\left[Z_1\zeta_2(T) \right] \right)+ \notag \\
    && \phantom{ \partial^2 \ell }
  +\frac{\lambda^4(\mean\left[ Z^2_1\right]+\mean\left[  Z^2_2\right]-2\mean\left[  Z_1Z_2\right]\lambda)}{\Omega^2_{11}(1-\lambda^{2})^{3}} 
  - \frac{1}{2\Omega^2_{11}} - \frac{\lambda^4}{2\Omega^2_{11}(1-\lambda^{2})^{2}},\\
i_{34}&=& -\mean\left[\frac{\partial^2 \ell}{\partial \Omega_{11} \partial \Omega_{12} } \right] =
  \frac{\lambda+\mean\left[ Z_1Z_2 \right] }{(1-\lambda^{2}) \Omega^{3/2}_{11}\sqrt{\Omega_{22}}} - 
  \frac{2\lambda \mean\left[ Z^2_1 \right]  + \lambda\mean\left[ Z^2_2 \right] - 5\lambda^2 \mean\left[ Z_1Z_2 \right] - \lambda^3}                 {(1-\lambda^{2})^2 \Omega^{3/2}_{11}\sqrt{\Omega_{22}}} + \notag\\
  && \phantom{ \partial^2 \ell }  
  -\frac{2\lambda^3(\mean\left[ Z^2_1\right]+\mean\left[  Z^2_2\right]-2\mean\left[  Z_1Z_2\right]\lambda)}                                       {(1-\lambda^{2})^3\Omega^{3/2}_{11}\sqrt{\Omega_{22}}} +\notag\\
  && \phantom{ \partial^2 \ell }      
  -\left( \frac{\alpha^2_1\alpha^2_2\tau\lambda}{2\Omega^{3/2}_{11}\sqrt{\Omega_{22}} (1+\alpha^2_*)^{3/2} } -
  \frac{\alpha_1\alpha_2\tau}{2\Omega^{3/2}_{11}\sqrt{\Omega_{22}}\sqrt{1+\alpha^2_*}}\right) \mean\left[  \zeta_1(T) \right]   +\notag\\ 
    && \phantom{ \frac{\partial^2 \ell}{\partial \Omega^2_{11} } }
  +\frac{\alpha_1\alpha_2\tau}{2\Omega^{3/2}_{11}\sqrt{\Omega_{22}}\sqrt{1+\alpha^2_*}}
   \left( \frac{\alpha_1\alpha_2\lambda\tau}{\sqrt{1+\alpha^2_*}}\mean\left[\zeta_2(T) \right] + \alpha_1\mean\left[Z_1\zeta_2(T) \right]\right),\\
i_{35}&=&-\mean\left[\frac{\partial^2 \ell}{\partial \Omega_{11}\partial \Omega_{22} } \right]=
  - \frac{\lambda^4 + 6\lambda^3\mean\left[Z_1Z_2 \right]-2\lambda^2\mean\left[Z^2_1\right]-2\lambda^2\mean\left[ Z^2_2 \right]}                {2\Omega_{11}\Omega_{22}  (1-\lambda^{2})^{2}} + \notag\\
   && \phantom{ \partial^2 \ell }
  -\frac{\lambda^2}{2\Omega_{11}\Omega_{22}(1-\lambda^{2})} 
  -\frac{\lambda\mean\left[  Z_1Z_2\right] }{\Omega_{11}\Omega_{22}(1-\lambda^2)} + \notag \\
  && \phantom{ \partial^2 \ell }
  +\frac{\lambda^4(\mean\left[ Z^2_1\right]+\mean\left[  Z^2_2\right]-2\mean\left[  Z_1Z_2\right]\lambda)}{\Omega_{11}\Omega_{22}(1-\lambda^{2})^{3}} +\notag\\ 
  && \phantom{ \partial^2 \ell }
  -\frac{\alpha_1\alpha_2\lambda\tau}{4\Omega_{11}\Omega_{22}\sqrt{1+\alpha^2_*}}   
  \left( 1 - \frac{\alpha_1\alpha_2\lambda}{\sqrt{1+\alpha^2_*}}\right)  \mean\left[  \zeta_1(T) \right]  + \notag \\
  && \phantom{ \partial^2 \ell }  
  -\frac{\alpha_1\alpha_2}{4\Omega_{11}\Omega_{22}} 
  \Bigg( \frac{\alpha_1\alpha_2\lambda^2\tau^2}{1+\alpha^2_*} \mean\left[\zeta_2(T) \right] +
      \frac{\alpha_2\lambda\tau}{\sqrt{1+\alpha^2_*} }  \mean\left[Z_2\zeta_2(T) \right]  + \notag\\
        && \phantom{\frac{\alpha_1\alpha_2}{4\Omega_{11}\Omega_{22}} }  
      +\frac{\alpha_1\lambda\tau}{\sqrt{1+\alpha^2_*} }  \mean\left[Z_1\zeta_2(T) \right]  + 
      \mean\left[Z_1Z_2\zeta_2(T) \right] \Bigg),\\
i_{36}&=& - \mean\left[ \frac{\partial^2 \ell}{\partial \Omega_{11}\partial \alpha_1 } \right] =
  -\frac{1}{2\Omega_{11}}\Bigg( \frac{\alpha_1\alpha_2\lambda(\alpha_2\lambda+\alpha_1)\tau}  {(1+\alpha^2_*)^{3/2}}\mean\left[  \zeta_1(T) \right] 
   -\frac{\alpha_2\lambda\tau}{\sqrt{1+\alpha^2_*}}\mean\left[  \zeta_1(T) \right] + \notag \\
    && \phantom{ \frac{\partial^2 \ell}{\partial \Omega^2_{11} } }
    - \mean\left[ Z_1 \zeta_1(T) \right]  \Bigg) + \notag \\
    && \phantom{ \frac{\partial^2 \ell}{\partial \Omega^2_{11} } }    
  +\frac{1}{2\Omega_{11}} 
  \Bigg( \frac{\alpha_1\alpha_2\lambda\tau^2 (\alpha_2\lambda+\alpha_1) }{1+\alpha^2_*} \mean\left[\zeta_2(T) \right] +
      \frac{\alpha_1\alpha_2\lambda\tau }{\sqrt{1+\alpha^2_*}}\mean\left[Z_1\zeta_2(T) \right]+ \notag \\
        && \phantom{ \frac{\partial^2 \ell}{\partial \Omega^2_{11} } }  
      +\frac{\alpha_1(\alpha_2\lambda+\alpha_1)\tau}{\sqrt{1+\alpha^2_*}} \mean\left[Z_1\zeta_2(T) \right]+
      \alpha_1 \mean\left[Z_1^2\zeta_2(T) \right]
   \Bigg) ,\\
i_{37}&=& - \mean\left[ \frac{\partial^2 \ell}{\partial \Omega_{11}\partial\alpha_2 }\right]  =
    -\frac{1}{2\Omega_{11}}\left( \frac{\alpha_1\alpha_2\lambda(\alpha_1\lambda+\alpha_2)\tau}  {(1+\alpha^2_*)^{3/2}} 
    -\frac{\alpha_1\lambda\tau}{\sqrt{1+\alpha^2_*}}\right) \mean\left[  \zeta_1(T) \right] + \notag \\
    && \phantom{ \frac{\partial^2 \ell}{\partial \Omega^2_{11} } }    
    +\frac{1}{2\Omega_{11}} 
  \Bigg( \frac{\alpha_1\alpha_2\lambda\tau^2 (\alpha_2+\alpha_1\lambda) }{1+\alpha^2_*} \mean\left[\zeta_2(T) \right] +
      \frac{\alpha_1\alpha_2\lambda\tau }{\sqrt{1+\alpha^2_*}}\mean\left[Z_2\zeta_2(T) \right]+ \notag \\
        && \phantom{ \frac{\partial^2 \ell}{\partial \Omega^2_{11} } }  
      +\frac{\alpha_1(\alpha_2+\alpha_1\lambda)\tau}{\sqrt{1+\alpha^2_*}} \mean\left[Z_1\zeta_2(T) \right]+
      \alpha_1 \mean\left[Z_1Z_2\zeta_2(T) \right]
   \Bigg) ,\\
i_{38}&=&- \mean\left[ \frac{\partial^2 \ell}{\partial \Omega_{11}\partial\tau} \right]  = 
   \frac{\alpha_1\alpha_2\lambda}{2\Omega_{11}\sqrt{1+\alpha^2_*}} \mean\left[\zeta_1(T) \right] +\notag \\
      && \phantom{ \frac{\partial^2 \ell}{\partial \Omega^2_{11} } }  
   +\frac{\sqrt{1+\alpha^2_*} }{2\Omega_{11}} \left( \frac{\alpha_1\alpha_2\lambda\tau}{\sqrt{1+\alpha^2_*}} \mean\left[\zeta_2(T) \right]+ \alpha_1\mean\left[Z_1\zeta_2(T) \right] \right),
\end{eqnarray}  
\begin{eqnarray}
i_{44}&=&- \mean\left[\frac{\partial^2 \ell}{\partial \Omega^2_{12} } \right] =
  -\frac{1}{(1-\lambda^2)\Omega_{11}\Omega_{22}} - 
  \frac{6\lambda \mean\left[Z_1 Z_2\right]  -\mean\left[ Z^2_1 \right]  -\mean\left[Z^2_2 \right]  +2\lambda^2  } {(1-\lambda^2)^2\Omega_{11}\Omega_{22}} + \notag \\
  && \phantom{ \frac{\partial^2 \ell}{\partial \Omega^2_{22} } }    
  +\frac{4\lambda^2(\mean\left[Z^2_1\right] + \mean\left[Z^2_2\right]  -2 \mean\left[ Z_1 Z_2\right] \lambda)}{(1-\lambda^2)^3\Omega_{11}\Omega_{22}} +\notag \\
    && \phantom{ \frac{\partial^2 \ell}{\partial \Omega^2_{11} } }    
  -\frac{\alpha^2_1\alpha^2_2\tau}{\Omega_{11}\Omega_{22} (1+\alpha^2_*)} \left( \tau \mean\left[ \zeta_2(T) \right] -
  \frac{ \mean\left[ \zeta_1(T) \right] }{\sqrt{1+\alpha^2_*}} \right) ,  \\  
\notag
i_{45}&=&-\mean\left[\frac{\partial^2 \ell}{\partial \Omega_{12} \partial \Omega_{22}  } \right] =
  \frac{\lambda+\mean\left[ Z_1Z_2 \right] }{(1-\lambda^{2}) \Omega^{3/2}_{22}\sqrt{\Omega_{11}}} - 
  \frac{2\lambda \mean\left[ Z^2_2 \right]  + \lambda\mean\left[ Z^2_1 \right] - 5\lambda^2 \mean\left[ Z_1Z_2 \right] - \lambda^3}                 {(1-\lambda^{2})^2 \Omega^{3/2}_{22}\sqrt{\Omega_{11}}} + \notag\\
  && \phantom{ \partial^2 \ell }  
  -\frac{2\lambda^3(\mean\left[ Z^2_1\right]+\mean\left[  Z^2_2\right]-2\mean\left[  Z_1Z_2\right]\lambda)}                                       {(1-\lambda^{2})^3\Omega^{3/2}_{22}\sqrt{\Omega_{11}}} +\notag\\
  && \phantom{ \partial^2 \ell }      
  -\left( \frac{\alpha^2_1\alpha^2_2\tau\lambda}{2\Omega^{3/2}_{22}\sqrt{\Omega_{11}} (1+\alpha^2_*)^{3/2} } -
  \frac{\alpha_1\alpha_2\tau}{2\Omega^{3/2}_{22}\sqrt{\Omega_{11}}\sqrt{1+\alpha^2_*}}\right) \mean\left[  \zeta_1(T) \right]   +\notag\\ 
    && \phantom{ \frac{\partial^2 \ell}{\partial \Omega^2_{11} } }
  +\frac{\alpha_1\alpha_2\tau}{2\Omega^{3/2}_{22}\sqrt{\Omega_{11}}\sqrt{1+\alpha^2_*}}
   \left( \frac{\alpha_1\alpha_2\lambda\tau}{\sqrt{1+\alpha^2_*}}\mean\left[\zeta_2(T) \right] + \alpha_2\mean\left[Z_2\zeta_2(T) \right]\right),\\
i_{46}&=&- \mean\left[\frac{\partial^2 \ell}{\partial \Omega_{12}\partial\alpha_1 }  \right] =
  - \frac{\alpha_2\tau}{\sqrt{\Omega_{11}\Omega_{22}} \sqrt{1+\alpha^2_*}} 
  \left( 1- \frac{\alpha_1(\alpha_2\lambda+\alpha_1)}{1+\alpha^2_*}\right)\mean\left[ \zeta_1(T) \right] +\notag\\
    && \phantom{ \frac{\partial^2 \ell}{\partial \Omega^2_{11} } }    
  -\frac{\alpha_1\alpha_2\tau}{\sqrt{\Omega_{11}\Omega_{22}}\sqrt{1+\alpha^2_*}}
  \left( \frac{(\alpha_1+\lambda\alpha_2)\tau}{\sqrt{1+\alpha^2_*}} \mean\left[ \zeta_2(T) \right]+\mean\left[ Z_1\zeta_2(T) \right] \right),\\
i_{47}&=&- \mean\left[  \frac{\partial^2 \ell}{\partial \Omega_{12}\partial\alpha_2 } \right] =
  -\frac{\alpha_1\tau}{\sqrt{\Omega_{11}\Omega_{22}} \sqrt{1+\alpha^2_*}} 
  \left( 1- \frac{\alpha_2(\alpha_1\lambda+\alpha_2)}{1+\alpha^2_*}\right)\mean\left[  \zeta_1(T)\right]+\notag\\
    && \phantom{ \frac{\partial^2 \ell}{\partial \Omega^2_{11} } }    
  -\frac{\alpha_1\alpha_2\tau}{\sqrt{\Omega_{11}\Omega_{22}}\sqrt{1+\alpha^2_*}}
  \left( \frac{(\alpha_2+\lambda\alpha_1)\tau}{\sqrt{1+\alpha^2_*}} \mean\left[ \zeta_2(T) \right]+\mean\left[ Z_2\zeta_2(T) \right]  \right),\\
i_{48}&=&- \mean\left[ \frac{\partial^2 \ell}{\partial \Omega_{12}\partial\tau} \right] =
  -\frac{\alpha_1\alpha_2}{\sqrt{\Omega_{11}\Omega_{22}}}
  \left( \frac{\mean\left[ \zeta_1(T) \right]}{\sqrt{1+\alpha^2_*}} +\tau\mean\left[ \zeta_2(T) \right] \right),
\end{eqnarray}  
\begin{eqnarray}
i_{55}&=&-\mean\left[\frac{\partial^2 \ell}{\partial\Omega^2_{22} } \right] =
  - \frac{\lambda^2-\mean\left[Z^2_2 \right] +2\mean\left[ Z_1Z_2\right] \lambda}{\Omega^2_{22}(1-\lambda^{2})} 
  -   \frac{4\lambda^3\mean\left[ Z_1Z_2\right] -2\lambda^2 \mean\left[ Z^2_2 \right] -\lambda^2\mean\left[  Z^2_1 \right]}{\Omega^2_{22}(1-\lambda^{2})^{2}} + \notag \\
    && \phantom{ \partial^2 \ell }
  +\frac{\lambda^4(\mean\left[ Z^2_1\right]+\mean\left[  Z^2_2\right]-2\mean\left[  Z_1Z_2\right]\lambda)}{\Omega^2_{22}(1-\lambda^{2})^{3}} 
  - \frac{1}{2\Omega^2_{22}} - \frac{\lambda^4}{2\Omega^2_{22}(1-\lambda^{2})^{2}} + \notag \\
    && \phantom{ \partial^2 \ell }
  -\frac{1}{4\Omega^2_{22}} 
  \left( \frac{3\alpha_1\alpha_2\tau\lambda}{\sqrt{1+\alpha^2_*}} \mean\left[\zeta_1(T) \right]
  - \frac{\alpha^2_1\alpha^2_2\tau\lambda^2}{(1+\alpha^2_*)^{3/2}} \mean\left[\zeta_1(T) \right]
  + 3\alpha_2 \mean\left[Z_2\zeta_1(T) \right] \right)  + \notag \\
    && \phantom{ \partial^2 \ell }
  -\frac{1}{4\Omega^2_{22}} \left( \frac{\alpha^2_1\alpha^2_2\lambda^2\tau^2}{1+\alpha^2_*}\mean\left[\zeta_2(T)\right] +
  \alpha^2_2\mean\left[Z^2_2\zeta_2(T) \right]+ 2\frac{\alpha_1\alpha^2_2\lambda\tau}{\sqrt{1+\alpha^2_*}} \mean\left[Z_2\zeta_2(T) \right] \right),\\
i_{56}&=&- \mean\left[ \frac{\partial^2 \ell}{\partial \Omega_{22}\partial\alpha_1 }\right]  =  
    -\frac{1}{2\Omega_{22}}\left( \frac{\alpha_1\alpha_2\lambda(\alpha_2\lambda+\alpha_1)\tau}  {(1+\alpha^2_*)^{3/2}} 
    -\frac{\alpha_2\lambda\tau}{\sqrt{1+\alpha^2_*}}\right) \mean\left[  \zeta_1(T) \right] + \notag \\
    && \phantom{ \frac{\partial^2 \ell}{\partial \Omega^2_{22} } }    
    +\frac{1}{2\Omega_{22}} 
  \Bigg( \frac{\alpha_1\alpha_2\lambda\tau^2 (\alpha_1+\alpha_2\lambda) }{1+\alpha^2_*} \mean\left[\zeta_2(T) \right] +
      \frac{\alpha_1\alpha_2\lambda\tau }{\sqrt{1+\alpha^2_*}}\mean\left[Z_1\zeta_2(T) \right]+ \notag \\
        && \phantom{ \frac{\partial^2 \ell}{\partial \Omega^2_{22} } }  
      +\frac{\alpha_2(\alpha_1+\alpha_2\lambda)\tau}{\sqrt{1+\alpha^2_*}} \mean\left[Z_2\zeta_2(T) \right]+
      \alpha_2 \mean\left[Z_1Z_2\zeta_2(T) \right]
   \Bigg),\\
i_{57}&=&- \mean\left[ \frac{\partial^2 \ell}{\partial \Omega_{22}\partial \alpha_2 } \right] = 
  -\frac{1}{2\Omega_{22}}\Bigg( \frac{\alpha_1\alpha_2\lambda(\alpha_1\lambda+\alpha_2)\tau}  {(1+\alpha^2_*)^{3/2}}\mean\left[  \zeta_1(T) \right] 
   -\frac{\alpha_1\lambda\tau}{\sqrt{1+\alpha^2_*}}\mean\left[  \zeta_1(T) \right]+  \notag \\
    && \phantom{ \frac{\partial^2 \ell}{\partial \Omega^2_{22} } }
    - \mean\left[ Z_2 \zeta_1(T) \right]  \Bigg) + \notag \\
    && \phantom{ \frac{\partial^2 \ell}{\partial \Omega^2_{22} } }    
  +\frac{1}{2\Omega_{22}} 
  \Bigg( \frac{\alpha_1\alpha_2\lambda\tau^2 (\alpha_1\lambda+\alpha_2) }{1+\alpha^2_*} \mean\left[\zeta_2(T) \right] +
      \frac{\alpha_1\alpha_2\lambda\tau }{\sqrt{1+\alpha^2_*}}\mean\left[Z_2\zeta_2(T) \right]+ \notag \\
        && \phantom{ \frac{\partial^2 \ell}{\partial \Omega^2_{22} } }  
      +\frac{\alpha_2(\alpha_1\lambda+\alpha_2)\tau}{\sqrt{1+\alpha^2_*}} \mean\left[Z_2\zeta_2(T) \right]+
      \alpha_2 \mean\left[Z_2^2\zeta_2(T) \right]
   \Bigg),\\
i_{58}&=&- \mean\left[ \frac{\partial^2 \ell}{\partial \Omega_{22}\partial\tau} \right] =
   \frac{\alpha_1\alpha_2\lambda}{2\Omega_{22}\sqrt{1+\alpha^2_*}} \mean\left[\zeta_1(T) \right] +\notag \\
      && \phantom{ \frac{\partial^2 \ell}{\partial \Omega^2_{22} } }  
   +\frac{\sqrt{1+\alpha^2_*} }{2\Omega_{22}} \left( \frac{\alpha_1\alpha_2\lambda\tau}{\sqrt{1+\alpha^2_*}} \mean\left[\zeta_2(T) \right]+ \alpha_2\mean\left[Z_2\zeta_2(T) \right] \right),
\end{eqnarray}  
\begin{eqnarray}
i_{66}&=&- \mean\left[\frac{\partial^2 \ell}{\partial \alpha^2_1 }\right] =
  -\left( \frac{\tau}{\sqrt{1+\alpha^2_*}} - \frac{(\alpha_2\lambda+\alpha_1)^2\tau}{(1+\alpha^2_*)^{3/2}} \right) \mean\left[\zeta_1(T)\right]  -
  \frac{(\alpha_1+\lambda\alpha_2)^2\tau^2}{1+\alpha^2_*} \mean\left[ \zeta_2(T) \right]+ \notag\\
  && \phantom{ \frac{\partial^2 \ell}{\partial \Omega^2_{11} } }    
  - \mean\left[ Z^2_1\zeta_2(T) \right] - 2\frac{(\alpha_1+\lambda\alpha_2)\tau}{\sqrt{1+\alpha^2_*}}\mean\left[ Z_1\zeta_2(T) \right],\\
i_{67}&=&- \mean\left[\frac{\partial^2 \ell}{\partial \alpha_1\partial\alpha_2} \right] =
  - \left( \frac{\lambda\tau}{\sqrt{1+\alpha^2_*}} - \frac{(\alpha_2+\lambda\alpha_1)(\alpha_1+\lambda\alpha_2)\tau}{(1+\alpha^2_*)^{3/2}} \right) \mean\left[\zeta_1(T)\right] +\notag\\
    && \phantom{ \frac{\partial^2 \ell}{\partial \Omega^2_{11} } }    
  -\frac{(\alpha_1+\lambda\alpha_2)(\alpha_2+\lambda\alpha_1)\tau^2}{1+\alpha^2_*} \mean\left[\zeta_2(T) \right] -
  \frac{(\alpha_1+\lambda\alpha_2)\tau}{\sqrt{1+\alpha^2_*}} \mean\left[ Z_2\zeta_2(T) \right] + \notag \\
    && \phantom{ \frac{\partial^2 \ell}{\partial \Omega^2_{11} } }    
  -\frac{(\alpha_2+\lambda\alpha_1)\tau}{\sqrt{1+\alpha^2_*}} \mean\left[ Z_1\zeta_2(T) \right] -
  \mean\left[Z_1Z_2\zeta_2(T) \right],\\
i_{68}&=&- \mean\left[\frac{\partial^2 \ell}{\partial \alpha_1\partial\tau} \right]  = 
  -\frac{\alpha_1+\lambda\alpha_2}{\sqrt{1+\alpha^2_*}} \mean\left[ \zeta_1(T)\right]  +\notag \\
  && \phantom{ \frac{\partial^2 \ell}{\partial \Omega^2_{11} } }  
    -\left(\frac{(\alpha_1+\lambda\alpha_2)\tau}{\sqrt{1+\alpha^2_*}}\mean\left[ \zeta_2(T)\right]+\mean\left[ Z_1\zeta_2(T)\right]\right) \sqrt{1+\alpha^2_*} ,  
\end{eqnarray}  
\begin{eqnarray}
i_{77}&=&- \mean\left[\frac{\partial^2 \ell}{\partial \alpha^2_2 } \right] =
  - \left( \frac{\tau}{\sqrt{1+\alpha^2_*}} - \frac{(\alpha_1\lambda+\alpha_2)^2\tau}{(1+\alpha^2_*)^{3/2}} \right) \mean\left[\zeta_1(T)\right]  -     
      \frac{(\alpha_2+\lambda\alpha_1)^2\tau^2}{1+\alpha^2_*} \mean\left[ \zeta_2(T) \right]+\notag\\
  && \phantom{ \frac{\partial^2 \ell}{\partial \Omega^2_{11} } }    
  -\mean\left[ Z^2_2\zeta_2(T) \right] - 2\frac{(\alpha_2+\lambda\alpha_1)\tau}{\sqrt{1+\alpha^2_*}}\mean\left[ Z_2\zeta_2(T) \right],\\
i_{78}&=&- \mean\left[\frac{\partial^2 \ell}{\partial \alpha_2\partial\tau} \right]  = 
  -\frac{\alpha_2+\lambda\alpha_1}{\sqrt{1+\alpha^2_*}} \mean\left[ \zeta_1(T)\right]  +\notag \\
      && \phantom{ \frac{\partial^2 \ell}{\partial \Omega^2_{11} } }  
    -\left(\frac{(\alpha_2+\lambda\alpha_1)\tau}{\sqrt{1+\alpha^2_*}}\mean\left[ \zeta_2(T)\right]+\mean\left[ Z_2\zeta_2(T)\right]\right) \sqrt{1+\alpha^2_*} ,
\end{eqnarray}  
\begin{eqnarray}
i_{88}&=&- \mean\left[\frac{\partial^2 \ell}{\partial \tau^2 } \right] =
  -(1+\alpha^2_*)  \mean\left[\zeta_2(T)\right]+\zeta_2(\tau).
\end{eqnarray}  

The  numerical computations of the quantities (\ref{eq:akesn}))--(\ref{eq:ak0esn})
had originally been performed by \citet{franco:2011} using
the function \texttt{cuhre} of the \Rlang\ package \texttt{R2Cuba} \citep{man:R2Cuba}. 
Since that package is no longer available on the CRAN site, it has been replaced 
in the present exposition by its successor \texttt{cubature}.
\citet{franco:2011} has also checked that each  element $i_{hk}$
of the expected information coincided with the matching
value obtained by numerical integration of the second derivative.

\subsection{Some peculiar aspects of the information matrix} \label{s:peculiar}

\subsubsection{Formal preliminaries} \label{s:formal-prelim}

Denote by $i(\nu)$ a generic entry on the main diagonal of the expected information matrix.
If we reparameterize the corresponding parameter from $\nu$ to $\psi(\nu)$, with a monotonic
differentiable functions, the expected Fisher information changes according to a simple rule,
that is, the information for $\psi$ is  
\begin{eqnarray}
  \left\{\psi'(\nu)\right\}^{-2}i(\nu)  \Bigg\vert_{\nu=\nu(\psi)}
\end{eqnarray}

For the univariate $\ESN$, the scale parameter has been parameterized by $\omega$,
while in the bivariate case by $\Omega_{ii}$, for $i=1,2$, where $\Omega_{ii}$ 
plays the role of the square on $\omega$. 
Therefore, when we transition from bivariate to univariate, the scale parameter is reparameterized as follows:
\begin{eqnarray}
  \left\{\psi'(\omega)\right\}^{-2}i(\omega)  \Bigg\vert_{\omega=\omega(\psi)}
\end{eqnarray}
where
\begin{eqnarray}
  \psi(\omega)=\Omega_{ii}=\omega^2,  
\end{eqnarray}
so that
\begin{eqnarray}
  \left\{\psi'(\omega)\right\}^{-2} = \frac{1}{4 \omega^2}.
\end{eqnarray}


As already mentioned, setting $\tau=0$ leads back to the skew-normal
distribution, and clearly, this fact must also hold for the information matrix. 
On setting $\tau=0$, we have checked numerically for various
combinations of the other seven parameters, the equality of the expected information
$\SN_{2}$ matrix as studied by \citet{art:azzarei2008} and the one for
the $\ESN_{2}$ studied here (clearly, after removal of the final row and
column which pertain to $\tau$).
This fact provides an indication towards the correctness of the results.
For a given combination of the parameters, it has been checked numerically
that the introduction of the eighth parameter $\tau=0$ leads to a decrease
of the determinant of the matrix with respect to the $SN_{2}$ distribution.


We recall the following result on pairwise conditional independence for
a $\ESN_k$ \rv\ \citep{art:capiteal03}.

\begin{prop} [Pairwise conditional independence]
  \label{prop:indipendenza}
  If $Y\sim \ESN_k(\xi,\Omega,\alpha,\tau)$, then
  \[ Y_i \indep Y_j \:|\: (\hbox{all other variables}) \]
  if and only if the following conditions simultaneously hold:\\
  \hspace*{1em}(a) $\Omega^{ij} = 0$,\\
  \hspace*{1em}(b) $\alpha_i \alpha_j = 0$\\
  where $\Omega^{ij}$ denotes the $(i,j)$-th entry of $\Omega^{-1}$.
\end{prop}

If $Y\sim \ESN_2(\xi,\Omega,\alpha,\tau)$, we obtain the following results
for specific configurations of the parameters:
\begin{itemize}
\item
when $\Omega_{12}=0$, $\alpha_1=0$ and $\tau=0$,  we have that
\begin{equation}
   f_{\ESN_2}(y;\xi,\Omega, \alpha, \tau) 
     = f_{\N}(y_1;\xi_1,\Omega_{11}) \times f_{\SN}(y_2;\xi_2,\Omega_{22},\alpha_2)
\end{equation}
in an obvious notation;
\item
when $\Omega_{12}=0$ and $\alpha_1=0$, we have that
\begin{equation}
  f_{\ESN_2}(y;\xi,\Omega, \alpha, \tau) 
    =f_{\N}(y_1;\xi_1,\Omega_{11}) \times f_{\ESN}(y_2;\xi_2,\Omega_{22},\alpha_2,\tau).
\end{equation}
\end{itemize}

On suitably rearranging the components, the expected information matrix of 
a $\ESN_{2}$ variable in the two special cases indicated above
can be partitioned as follows:
\[
i(\theta_A)= \left( \begin{tabular}{c|c}
 $i_N(\theta_1)$ &  0 \\ 
\hline  0 & $i_{SN}(\theta_2)$ \\ 
 \end{tabular} 
 \right), \quad
i(\theta_B)=\left(
\begin{tabular}{c|c}
 $i_N(\theta_1)$ &  0 \\ 
\hline 0  & $i_{ESN}(\theta_3)$ \\ 
 \end{tabular} \right)
 \]
where $\theta_A=(\xi_1,\xi_2,\Omega_{11},\Omega_{22},\alpha_2 )$, $\theta_B=(\xi_1,\xi_2,\Omega_{11},\Omega_{22},\alpha_2,\tau)$, $\theta_1=(\xi_1,\Omega_{11})$, $\theta_2=(\xi_2,\Omega_{22},\alpha_2)$ and $\theta_3=(\xi_2,\Omega_{22},\alpha_2,\tau)$.

\subsubsection{Numerical exploration}

Next, we shall examine some features of the expected information matrix on the 
basis of numerical and graphical evidence. Among these aspects, a relevant one is
the singularity of the matrix as  $\alpha \to 0$ and as $\tau \to \pm \infty$.

\begin{figure}[h]
  \centering
    \includegraphics[scale=0.8]{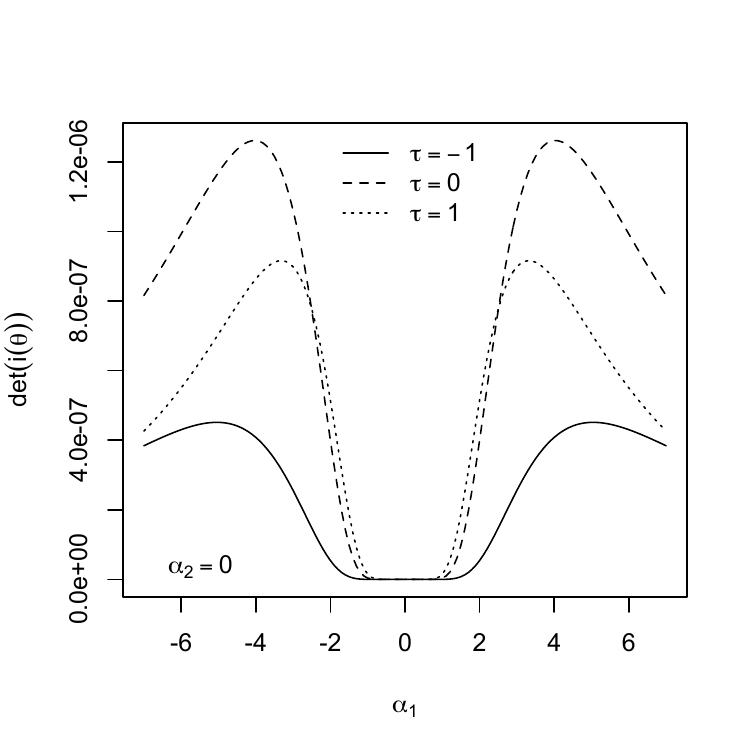}
  \caption{Determinant of the expected information matrix as a function of
  $\alpha_1$ for three values of $\tau$,
  when $\xi_1=\xi_2=0, \alpha_2=0$, $\Omega_{11}=\Omega_{22}=1, \Omega_{12}=0$.}
  \label{fig:alpha_1_det_alpha2_0_indip}
\end{figure}

On setting $\alpha_{2}=0$, the expected information matrix tends to singularity
as  $\alpha_1 \to 0$. On this respect, see Figure~\ref{fig:alpha_1_det_alpha2_0_indip} 
which displays three curves for different values of $\tau$; on the abscissa, 
we have $\alpha_{1}$ and on the ordinate we have the determinant
of the information matrix; 
the other parameters are as fixed at $\xi=(0,0)^T$, $\Omega=\mathrm{diag}(1,1)$.

The determinant of the matrix decreases to $0$ as $\alpha_1 \to 0$. 
This fact persists for each choice of $\tau$, whose value only affects
the speed of the descent to $0$ of the determinant.
If we swap the roles of $\alpha_1$ and $\alpha_2$, the behaviour 
of the determinant would remain the same.
 
\begin{figure}[h]
  \centering
    \includegraphics[scale=0.7]{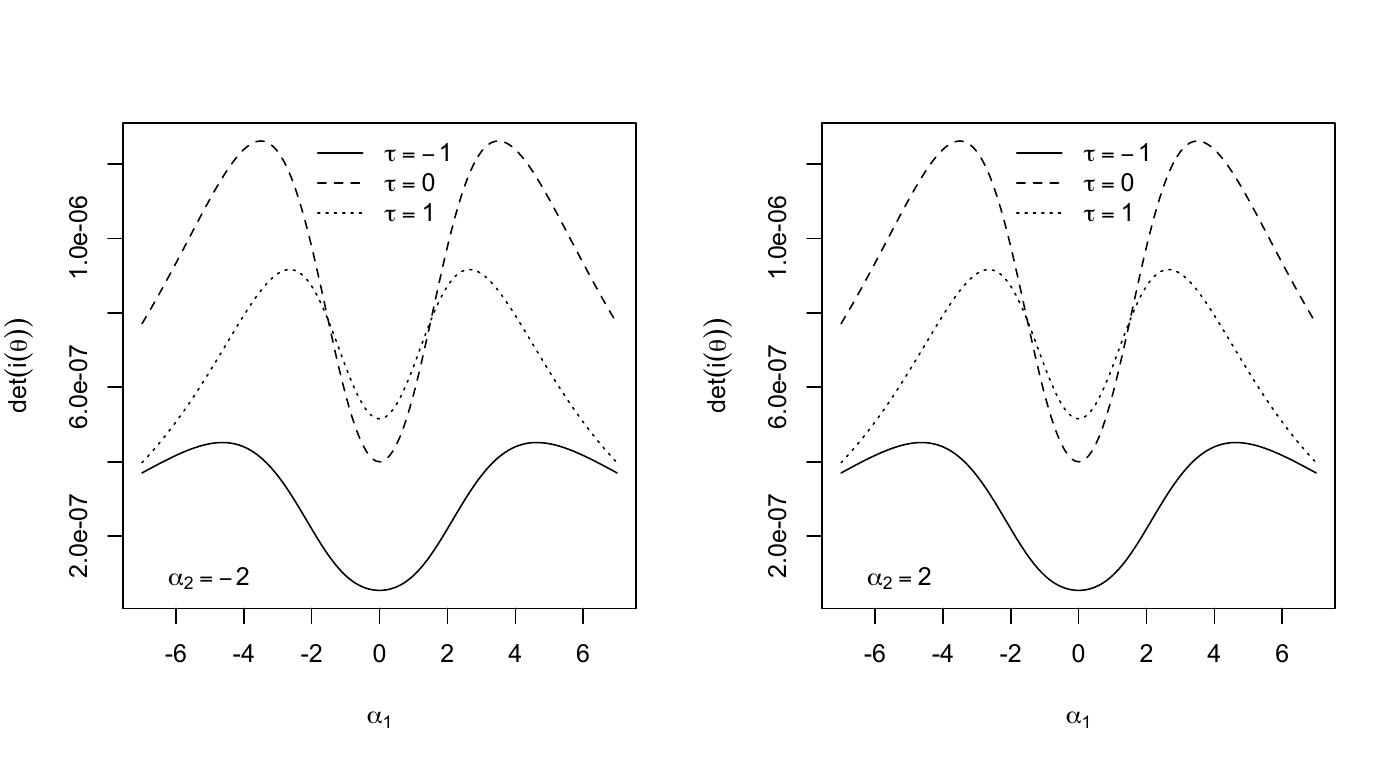}
  \caption{
  Determinant of the expected information matrix as a function of $\alpha_1$  
  for three values of $\tau$ and two values of $\alpha_2$, 
  when $\xi_1=\xi_2=0, \Omega_{11}=\Omega_{22}=1, \Omega_{12}=0$.}
  \label{fig:alpha_1_det_alpha2_2_vs_m2_indip_per2}
\end{figure}

From Figure~\ref{fig:alpha_1_det_alpha2_0_indip}, 
we see that a set of  $\alpha_1$ values for which the determinant tends
to $0$ are those between $-2$ and $2$ or about so. 
In Figure~\ref{fig:alpha_1_det_alpha2_2_vs_m2_indip_per2} it has then
been decided to fix $\alpha_2=-2$ and $\alpha_2=2$, for the two plots,
while keeping the other parameters unchanged. 
Notice that the two plots of Figure~\ref{fig:alpha_1_det_alpha2_2_vs_m2_indip_per2}
are identical; in fact, the sign of $\alpha$ does not affect the behaviour.
The new choice of the $\alpha_2$ values avoids that the determinant vanishes
as $\alpha_1 \to0$,  like in Figure~\ref{fig:alpha_1_det_alpha2_0_indip}, 
but it does vanish when $\alpha_1 \to \pm \infty$.
 
\begin{figure}
  \centering
    \includegraphics[scale=0.8]{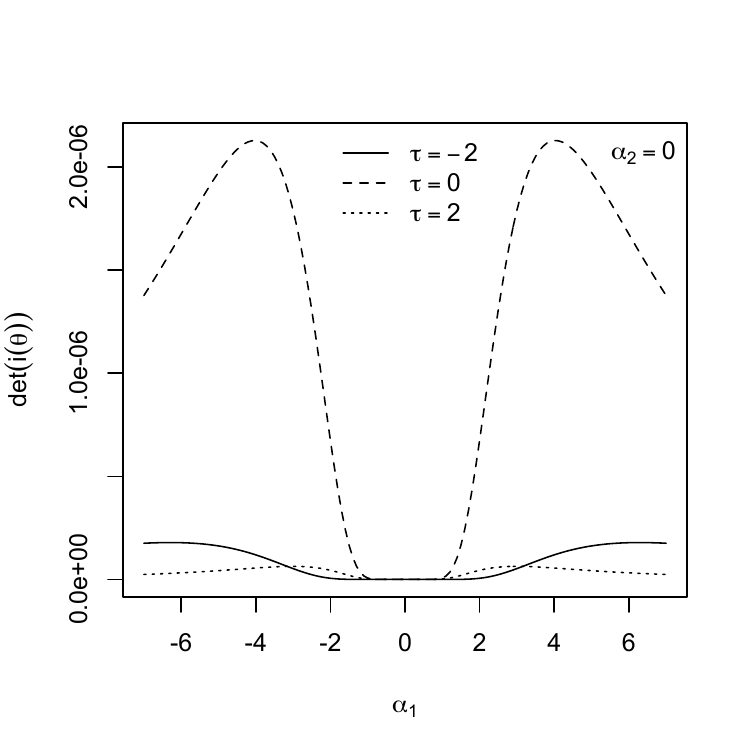}
  \caption{
  Determinant of the expected information matrix as a function of $\alpha_1$  
  for three values of $\tau$ when $\xi_1=\xi_2=0, \alpha_2=0$, $\Omega_{11}=\Omega_{22}=1, \Omega_{12}=0.4$.}
  \label{fig:alpha_1_det_alpha2_0_04_diversitau}
\end{figure}

Next, consider Figure~\ref{fig:alpha_1_det_alpha2_0_04_diversitau}, 
which refers to another  $\tau$ triplet, namely $(-2,0,2)$, and see
that the determinant vanishes as $\alpha_1\to0$ when $\alpha_2=0$.
A qualitatively simular situation had occurred in Figure~\ref{fig:alpha_1_det_alpha2_0_indip};
in this case the values of $\tau$ did not lead to any improvement
as for invertibility of the matrix when  $\alpha=(0,0)^T$.

\begin{figure}
  \centering
    \includegraphics[scale=0.7]{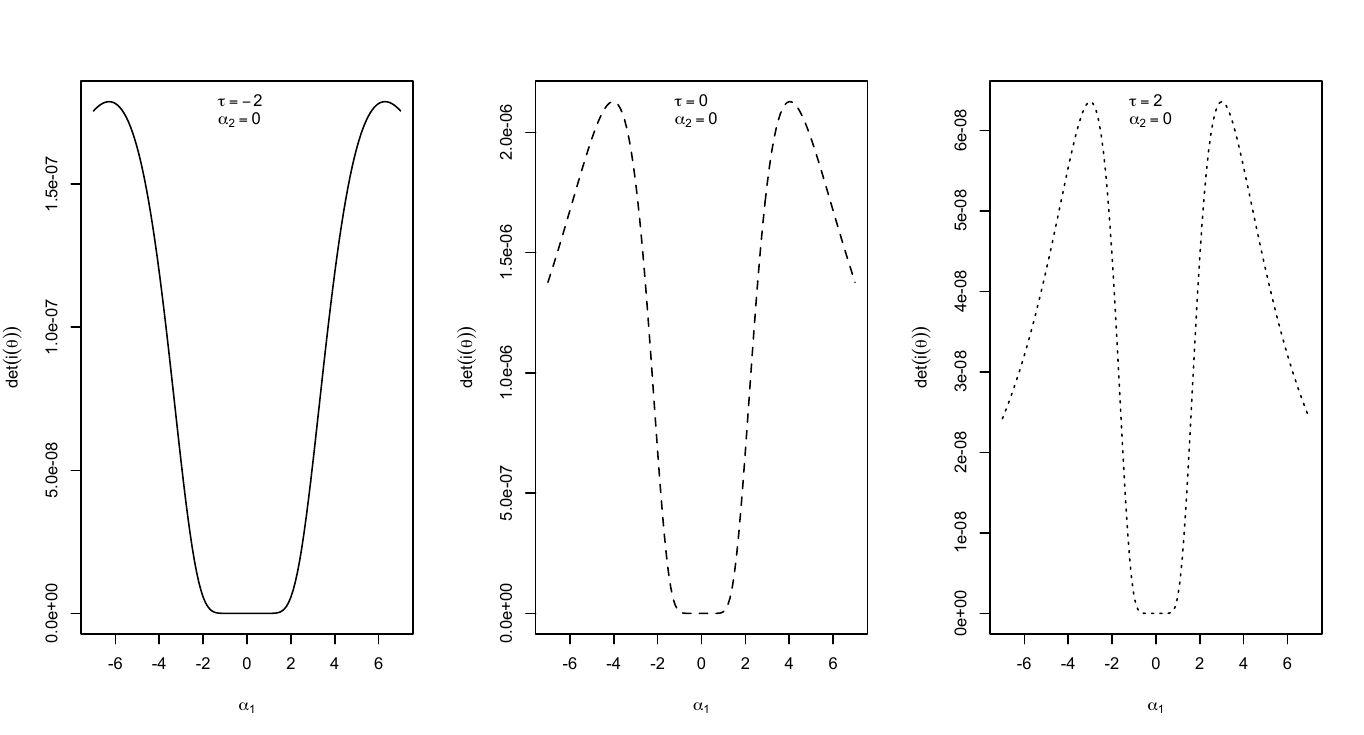}
  \caption{
  Determinant of the expected information matrix as a function of $\alpha_1$  
  for three values of $\tau$ when $\xi_1=\xi_2=0, \alpha_2=0$, 
  $\Omega_{11}=\Omega_{22}=1, \Omega_{12}=0.4$.}
\end{figure}

\begin{figure}
  \centering
    \includegraphics[scale=0.8]{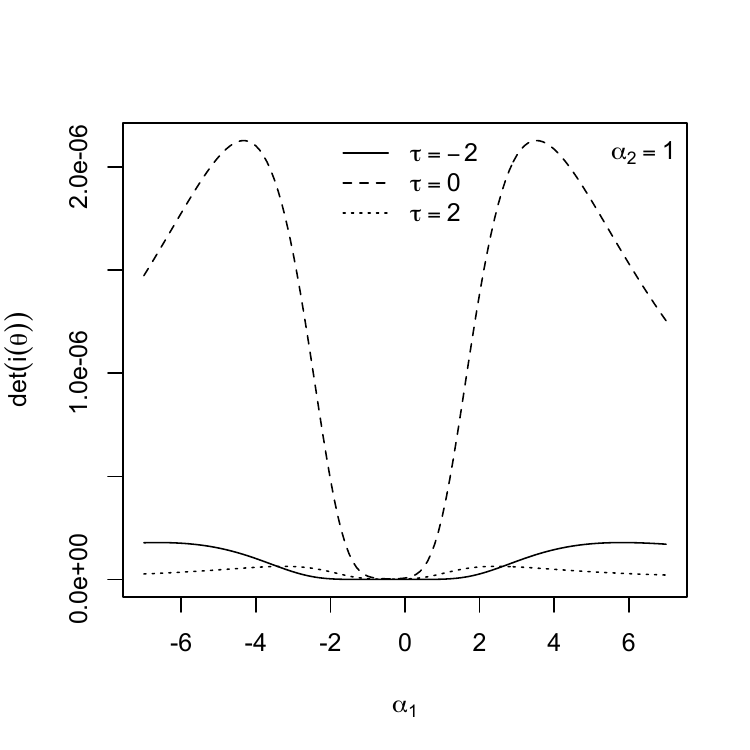}
  \caption{
  Determinant of the expected information matrix as a function of $\alpha_1$  
  for three values of $\tau$ when  
  $\xi_1=\xi_2=0, \alpha_2=1$, $\Omega_{11}=\Omega_{22}=1, \Omega_{12}=0.4$.}
  \label{fig:alpha_1_det_alpha2_1_04_diversitau}
\end{figure}

\begin{figure}
  \centering
    \includegraphics[scale=0.7]{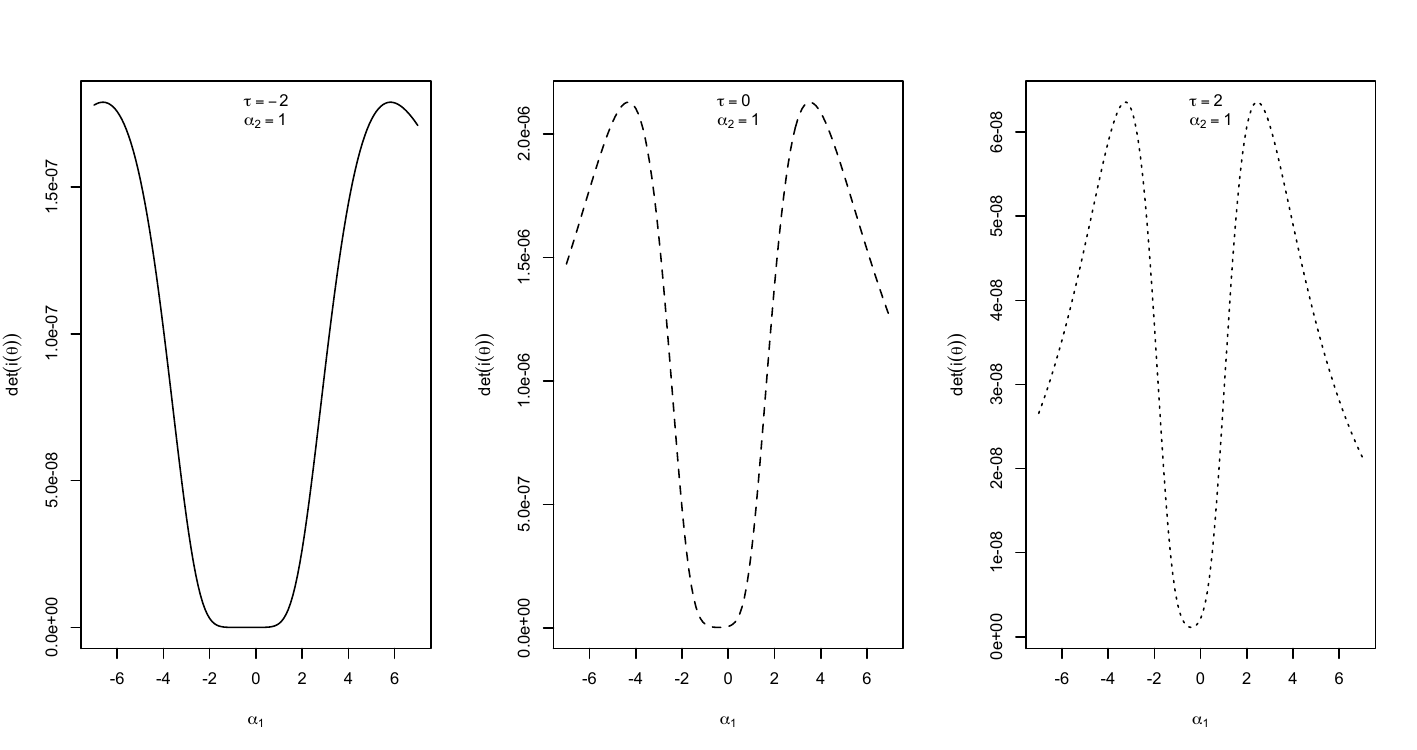}
  \caption{
  Determinant of the expected information matrix as a function of $\alpha_1$  
  for three values of $\tau$ when  
  $\xi_1=\xi_2=0, \alpha_2=1$, $\Omega_{11}=\Omega_{22}=1, \Omega_{12}=0.4$.}
  \label{fig:alpha_1_det_alpha2_1_04_diversitau_per3}
\end{figure}

\begin{figure}
  \centering
    \includegraphics[scale=0.8]{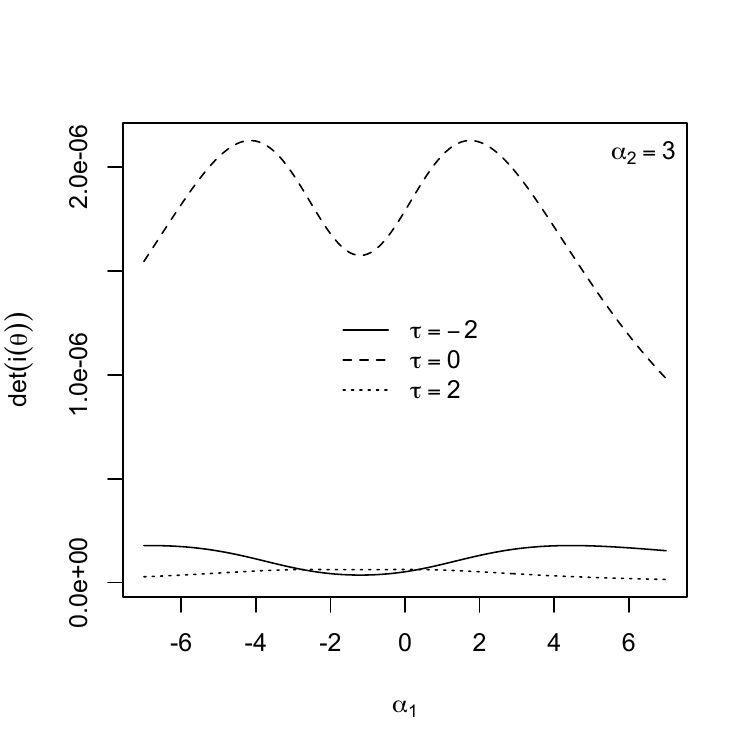}
  \caption{Determinant of the expected information matrix as a function of $\alpha_1$  
  for three values of $\tau$ when $\xi_1=\xi_2=0, \alpha_2=3$, $\Omega_{11}=\Omega_{22}=1, \Omega_{12}=0.4$.}
  \label{fig:alpha_1_det_alpha2_3_04_diversitau}
\end{figure}

\begin{figure}
  \centering
    \includegraphics[scale=0.7]{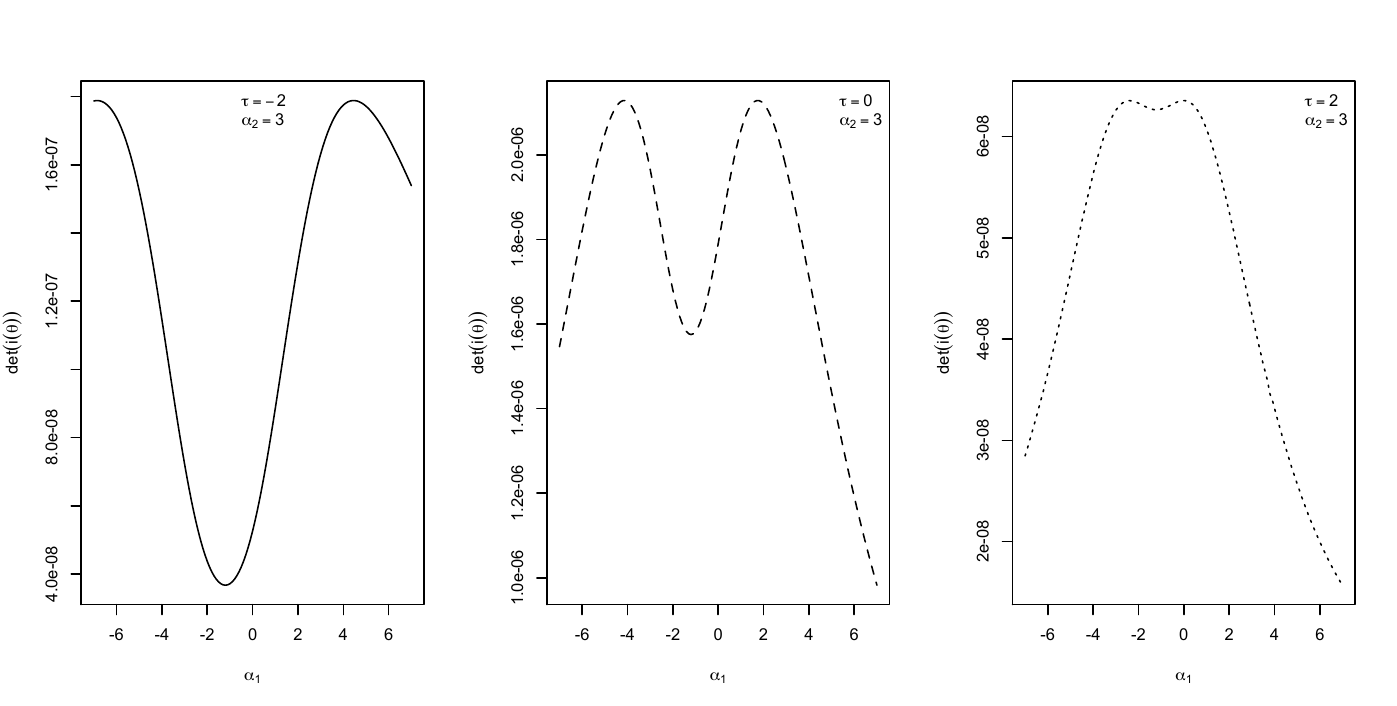}
  \caption{Determinant of the expected information matrix as a function of $\alpha_1$  
     for three values of $\tau$ when $\xi_1=\xi_2=0, \alpha_2=3$, $\Omega_{11}=\Omega_{22}=1, \Omega_{12}=0.4$.}
    \label{fig:alpha_1_det_alpha2_3_04_diversitau_per3}
\end{figure}

\begin{figure}
  \centering
    \includegraphics[scale=0.8]{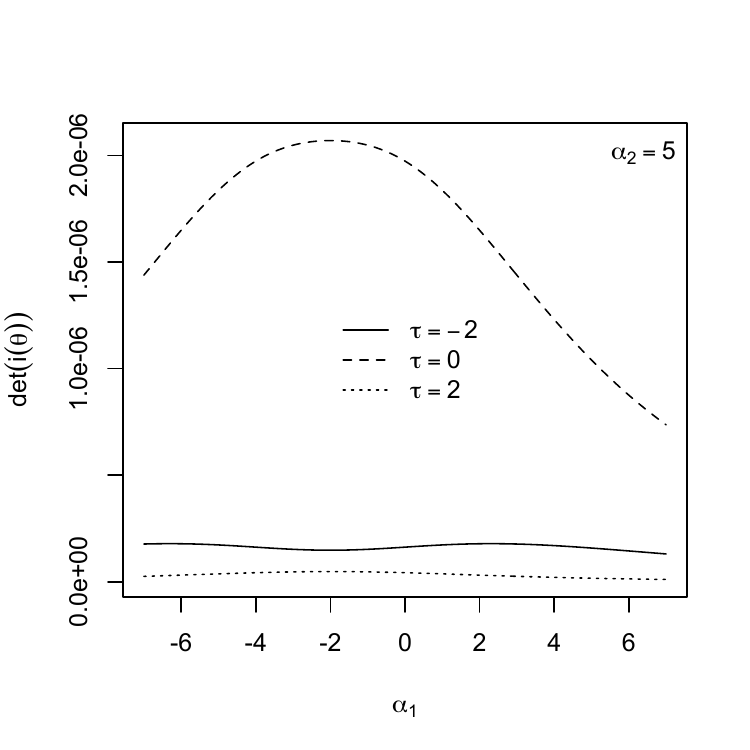}
  \caption{Determinant of the expected information matrix as a function of $\alpha_1$  
  for three values of $\tau$ when $\xi_1=\xi_2=0, \alpha_2=5$, $\Omega_{11}=\Omega_{22}=1, \Omega_{12}=0.4$.}
  \label{fig:alpha_1_det_alpha2_5_04_diversitau}
\end{figure}

\begin{figure}
  \centering
    \includegraphics[scale=0.7]{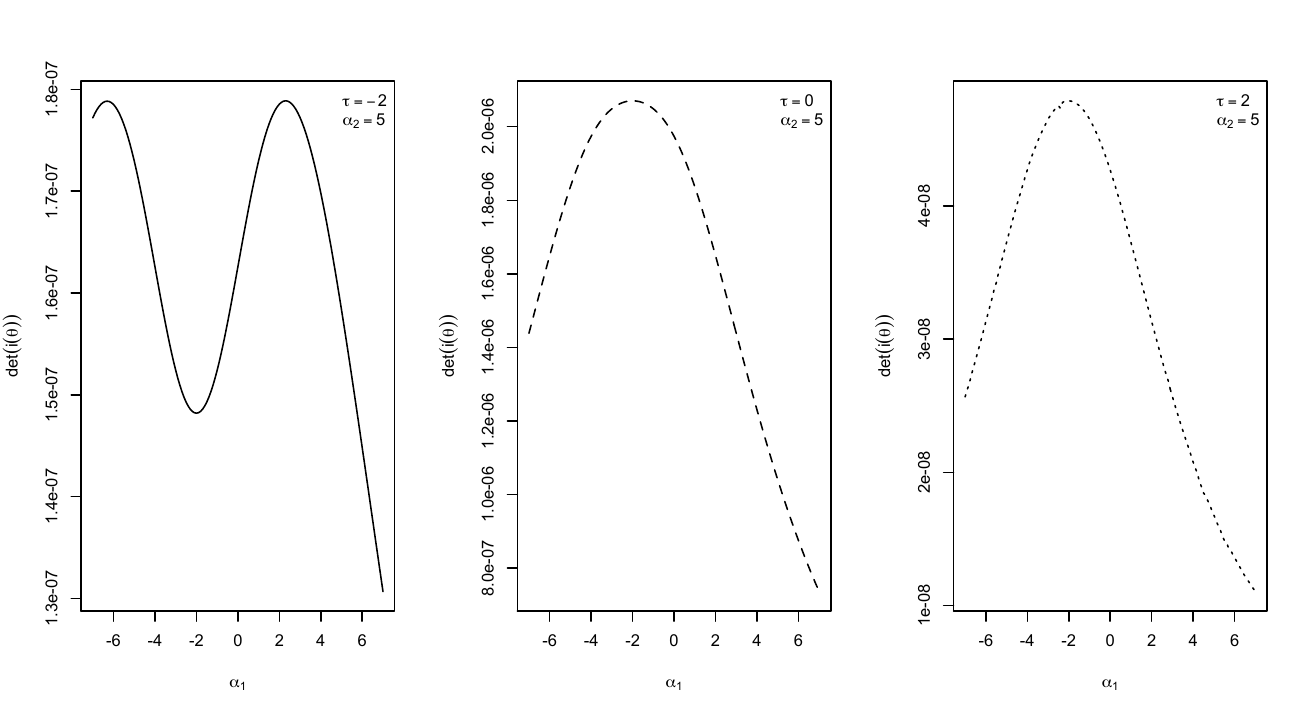}
  \caption{Determinant of the expected information matrix as a function of $\alpha_1$  
  for three values of $\tau$ when $\xi_1=\xi_2=0, \alpha_2=5$, 
  $\Omega_{11}=\Omega_{22}=1, \Omega_{12}=0.4$.}
  \label{fig:alpha_1_det_alpha2_5_04_diversitau_per3}
\end{figure}

In Figures~\ref{fig:alpha_1_det_alpha2_1_04_diversitau} to \ref{fig:alpha_1_det_alpha2_5_04_diversitau_per3},
the determinant goes to $0$ as  $\alpha_1 \to \pm \infty$, when $\alpha_2=1$, $3$ e $5$.
It appears that the matrix is singular only when $\alpha =(0,0)$, but not with other values
of the asymmetry parameter.
 
\begin{figure}
  \centering
    \includegraphics[scale=0.8]{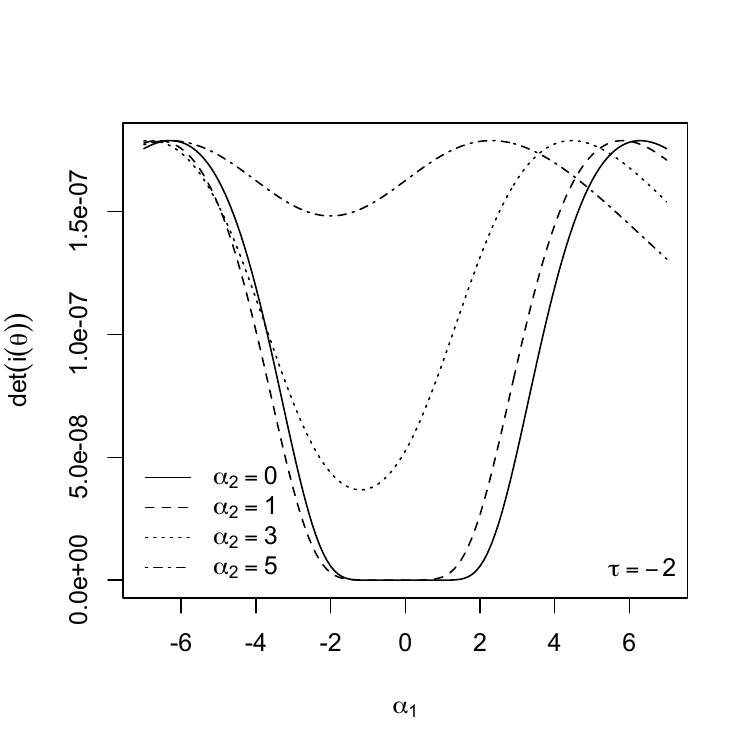}
  \caption{Determinante della matrice d'informazione attesa in funzione di $\alpha_1$ \\
  per quattro valori di $\alpha_2$, dove $\xi_1=\xi_2=0$, $\tau=-2$, $\Omega_{11}=\Omega_{22}=1, \Omega_{12}=0.4$}
  \label{fig:alpha_1_det_tau_m2_vari_alpha2}
\end{figure}

\begin{figure} 
  \centering
    \includegraphics[scale=0.8]{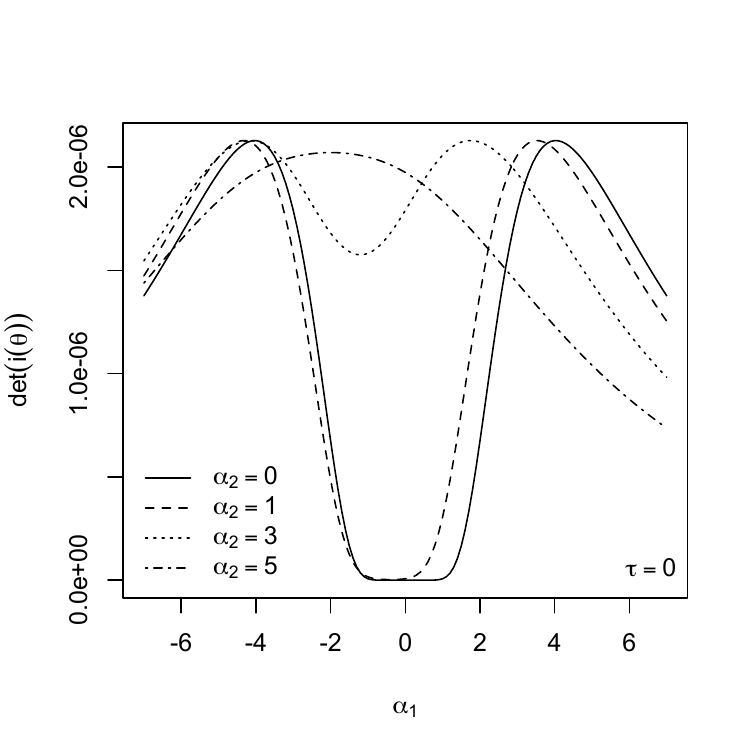}
  \caption{Determinant of the expected information matrix as a function of $\alpha_1$ 
         for four values of $\alpha_2$, when $\xi_1=\xi_2=0$, 
         $\tau=0$, $\Omega_{11}=\Omega_{22}=1, \Omega_{12}=0.4$.}
  \label{fig:alpha_1_det_tau_0_vari_alpha2}
\end{figure}

\begin{figure} 
  \centering
    \includegraphics[scale=0.8]{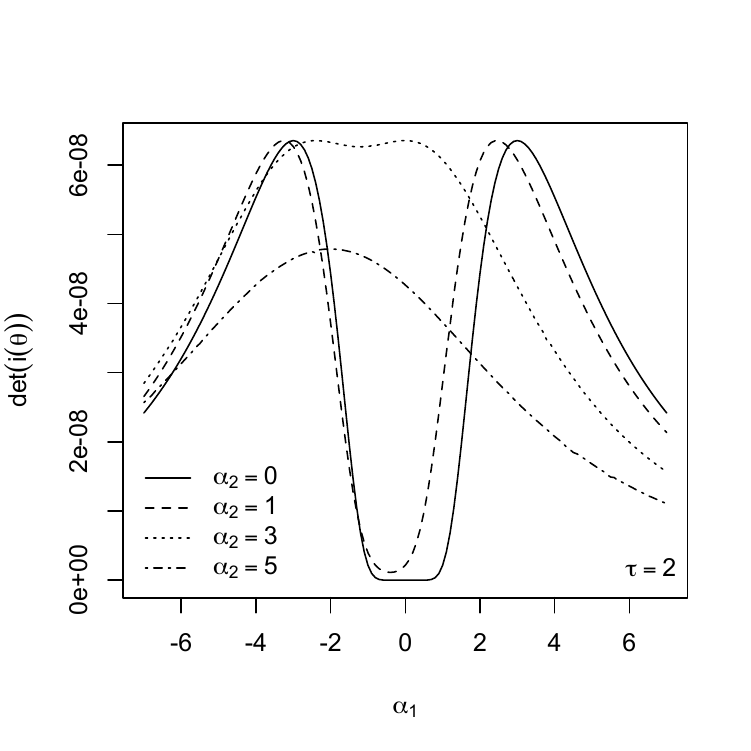}
  \caption{Determinant of the expected information matrix as a function of $\alpha_1$ 
  for four values of $\alpha_2$, when $\xi_1=\xi_2=0$, $\tau=-2$, $\Omega_{11}=\Omega_{22}=1, \Omega_{12}=0.4$.}
  \label{fig:alpha_1_det_tau_2_vari_alpha2}
\end{figure}

\begin{figure} 
  \centering
    \includegraphics[scale=0.75]{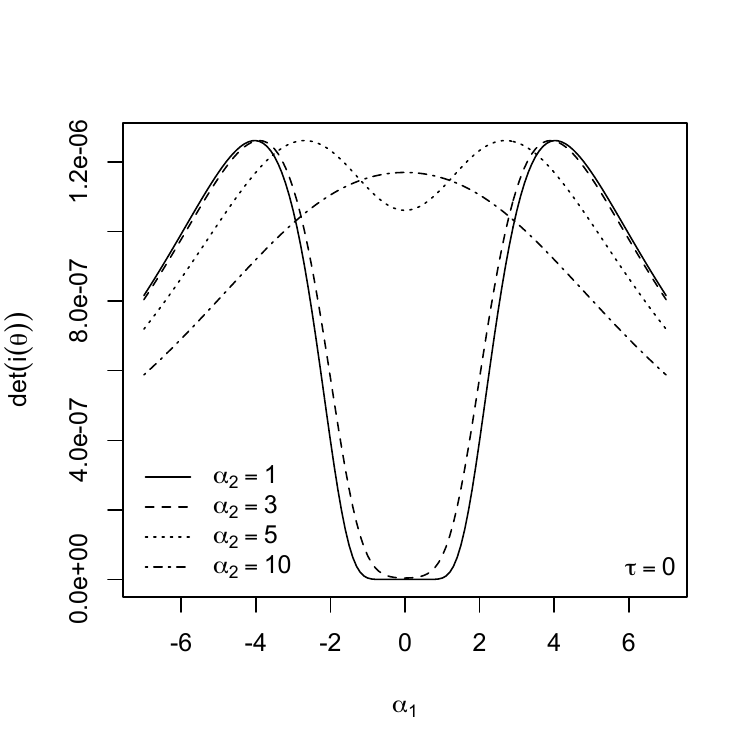} 
  \caption{Determinant of the expected information matrix as a function of $\alpha_1$ 
  for four values of $\alpha_2$, when $\xi_1=\xi_2=0$, $\tau=0$, $\Omega_{11}=\Omega_{22}=1, \Omega_{12}=0$.}
  \label{fig:alpha_1_det_varialpha2_Omega12_0_tau0}
\end{figure}

This pattern is even more clear from  Figures \ref{fig:alpha_1_det_tau_m2_vari_alpha2},
\ref{fig:alpha_1_det_tau_0_vari_alpha2}, \ref{fig:alpha_1_det_tau_2_vari_alpha2} and
\ref{fig:alpha_1_det_varialpha2_Omega12_0_tau0} where, for a given $\tau$,
the matrix appears to be singular only when  $\alpha =(0,0)^T$ and, 
as $\alpha_2$ increases, the determinant tends to vanish as $\alpha_1 \to \pm \infty$ .

Osserviamo ora gli stessi tipi di grafici invertendo il ruolo di $\alpha_1$ e
$\tau$, dove il resto dei parametri sono fissati come in precedenza.

\begin{figure}
  \centering
    \includegraphics[scale=0.75]{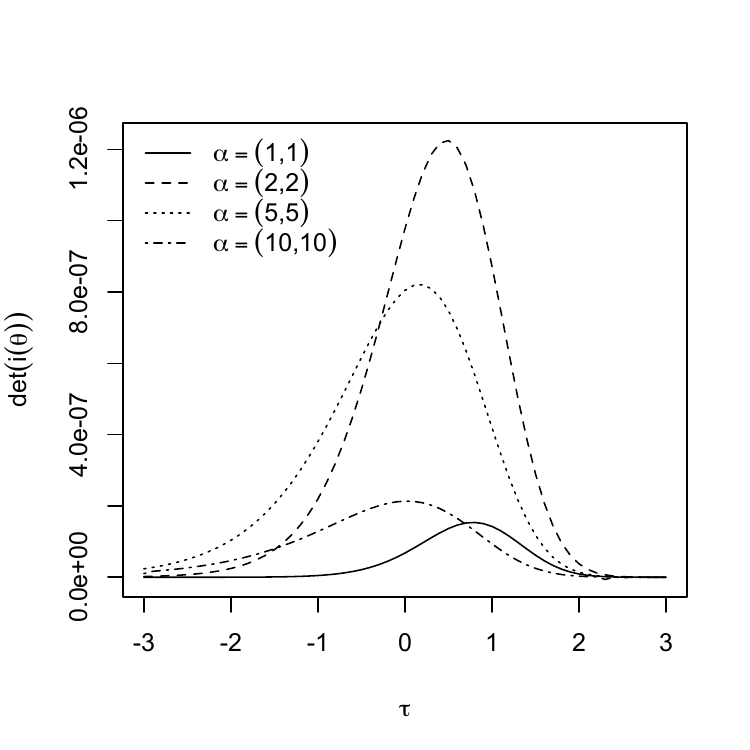}
  \caption{Determinant of the expected information matrix as a function of $\tau$ 
  for four choices of $\alpha$, when $\xi_1=\xi_2=0$, $\Omega_{11}=\Omega_{22}=1, \Omega_{12}=0$.}
  \label{fig:0_0_1_1_0_alpha_1_0_tau_1}
\end{figure}

\begin{figure}
  \centering
    \includegraphics[scale=0.75]{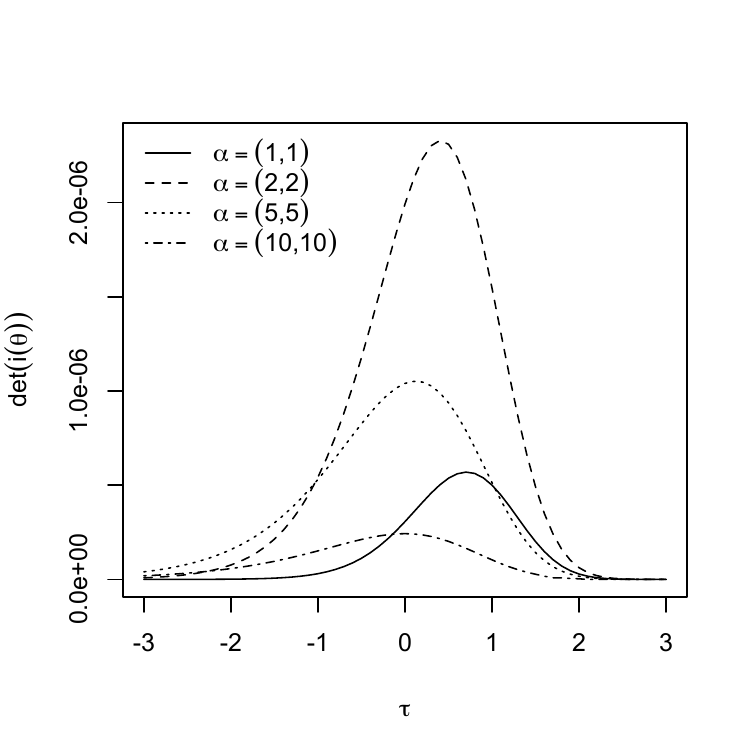}
  \caption{Determinant of the expected information matrix as a function of  $\tau$ 
  for four choices of $\alpha$, when $\xi_1=\xi_2=0$, $\Omega_{11}=\Omega_{22}=1, \Omega_{12}=0.4$.}
  \label{fig:0_0_1_1_0_4_alpha_1_0_tau_1}
\end{figure}

Next, consider the same sort of plot except that we swap the roles 
of $\alpha_{1}$ and $\tau$, while the remaining parameters stay the same as before.
In Fig.~\ref{fig:0_0_1_1_0_alpha_1_0_tau_1} and Fig.~\ref{fig:0_0_1_1_0_4_alpha_1_0_tau_1}, 
we see that the determinant of the information matrix tends to vanish when
$\tau \to \pm \infty$ and that it decreases very rapidly as $\alpha$ increases.

\subsection{Summary and final remarks}  
 
Starting from the density function of the bivariate ESN distribution, 
its observed information matrix has been developed in Section~\ref{cap:t3}.
Thanks to Lemma~\ref{lemmasteno}, which holds for the $k$-dimensional case,
it has been possible to compute all the quantities,  
from (\ref{eq:valorinoti1}) to (\ref{eq:E_noti2}), 
required for expressing the expected information matrix.

On setting $\tau=0$, a numerical check has confirmed the equality of the expected information matrix 
obtained here with the one of the $\SN_{2}$ studied by \citet{art:azzarei2008}; 
clearly, the final row and column corresponding to parameter $\tau$ have been omitted.
This check provides supporting evidence for the correctness of our results.
For a given combination of parameters, it has also been checked that the introduction of
the eight parameters $\tau=0$ leads to a decrease of the determinant with respect to 
the case $\SN_{2}$.
 
For some special cases indicated in Subsection~\ref{s:formal-prelim}, the $\ESN_{2}$ information matrix
can be partitioned in suitable blocks, corresponding to independent marginal distribution.

From the inferential viewpoint, the $\ESN_{2}$ distribution exhibits problems similar to those
of the classical $\SN$ distribution when the direct parameterization is adopted;
in fact, numerical and graphical outcomes indicate that the expected information matrix becomes
singular as $\alpha \to 0$ and as $\tau \to \pm \infty$. 
For the classical SN distribution, however, the singularity problem has been overcome
by adopting the so-called centred parameterization, which has definite advantages.
First of all, the new parameters have a clearer and more familiar meaning since they
are linked to the first three cumulants of the distribution.
In addition, for inferential work, the centred parameterization avoids issues
otherwise occurring with the SN distribution.

In \citet{tesi:canale}, it has emerged that  the centred parameterization for the univariate ESN,
following the well-known approach for the first three parameters while keeping $\tau$ unchanged,
has led to some improvement by removing the singularity of the information matrix at $\alpha=0$
and improving the distribution of the MLEs. However, this reparameterization also exhibits
some limitations; specifically, the lack of modification of $\tau$ causes a singularity problem
of the information matrix when $\tau \to \pm \infty$, and numerical difficulties in the maximization
of the log-likelihood.
Therefore, an alternative reparameterization of $\tau$ seems an interesting research 
target to consider.

%
\clearpage

\appendix                          
\section*{Appendix: \Rlang\ code}

\subsection*{A.1: The log-likelihood and related functions}

The \Rlang\ code below computes the $\ESN_{2}$ log-likelihood function, its score function
and the observed information matrix. 
The arbitrary constant of the log-likelihood function is taken to be $-\log 2\pi$,
so that the returned value equals the log of the density function.
The term `observed information matrix' is used here in a somewhat loose terminology,
since it actually denotes the hessian matrix of the log-likelihood function evaluated
at any DP choice, not necessarily the maximum likelihood point.

It is assumed that the \Rlang\ package \texttt{sn} has been loaded \citep{azzalini:R-pkg-sn}.

\begin{footnotesize}
\begin{verbatim}
#-----------------------------------------------------------------------------------------------
# ESN_2 log-likelihood function  
#-----------------------------------------------------------------------------------------------
# Description
# Computes the log-likelihood function for the extended bivariate skew-normal
# distribution in the DP parametrizations.
#
# Usage
# esn2.loglik(dp, y1, y2)
#
# Arguments:
#   dp      vector of direct parameters
#   y1, y2  vectors of observations, must have equal length
#
# Example:
#   y1 <- c(0.7, -0.4, 1.1);   
#   y2 <- c(-1.2, 0.5, 0.9)   
#   loglik <- esn2.loglik(dp=c(0, 0, 1, 0.6, 1, 2, 3, 1), y1, y2)
# where
#   location:  xi1 <- dp[1];   xi2 <- dp[2]
#   scale:     Omega11 <- dp[3]; Omega12 <- dp[4]; Omega22 <- dp[5]
#   trunc:     tau <- dp[8]
#------------------------------------------------------------------------------ 
esn2.loglik <- function(dp, y1, y2) 
{  
  if(length(y1) != length(y2)) stop("y1, y2 must have equal length")  
  dp.ok  <- (length(dp) == 8) & (dp[3]>0) & (dp[5]>0) & (dp[3]*dp[5]- dp[4]^2>0) 
  if(!dp.ok) stop("inadmissible dp vector")  
  
  xi1 <- dp[1]
  xi2 <- dp[2]
  Omega11 <- dp[3]
  Omega12 <- dp[4]
  Omega22 <- dp[5]
  alpha1 <- dp[6]
  alpha2 <- dp[7]
  tau <- dp[8]
  z1 <- (y1 - xi1)/sqrt(Omega11)
  z2 <- (y2 - xi2)/sqrt(Omega22)
  lambda <- Omega12/sqrt(Omega11 * Omega22)
  alpha0 <- tau * sqrt(1 + (alpha1^2) + (alpha2^2) + 2 * alpha1 * alpha2 * lambda)
  t <- (alpha0 + alpha1 * z1 + alpha2 * z2)
  zeta0tau <- (sn::zeta(0, tau))
  zeta0t <- (sn::zeta(0, t))
  sum(-log(2*pi) 
    -0.5 * log(Omega11) - 0.5 * log(Omega22) - 0.5 * log(1 - lambda^2) 
    - (0.5/(1 - lambda^2)) * (z1^2 + z2^2 - 2 * lambda * z1 * z2) + zeta0t 
    - zeta0tau)
}


#------------------------------------------------------------------------------ 
#  ESN_2 score function
#------------------------------------------------------------------------------ 
# Description
# Computes the score function for the extended bivariate skew-normal
# distribution in the DP parametrizations.
#
# Usage
# esn2.score(dp, y1, y2)
#
# Arguments:
#   dp      vector of direct parameters
#   y1, y2  vectors of observations, must have equal length
# 
# Example:
#   y1 <- c(0.7, -0.4, 1.1) 
#   y2 <- c(-1.2, 0.5, 0.9)   
#   info <- esn2.score(dp=c(0, 0, 1, 0.6, 1, 2, 3, 1), y1, y2)
# where
#   location:  xi1 <- dp[1];   xi2 <- dp[2]
#   scale:     Omega11 <- dp[3]; Omega12 <- dp[4]; Omega22 <- dp[5]
#   trunc:     tau <- dp[8]
#------------------------------------------------------------------------------
esn2.score <- function(dp, y1, y2)
{
  if(length(y1) != length(y2)) stop("y1, y2 must have equal length")
  dp.ok  <- (length(dp) == 8) & (dp[3]>0) & (dp[5]>0) & (dp[3]*dp[5]- dp[4]^2>0) 
  if(!dp.ok) stop("inadmissible dp vector") 
  
  xi1 <- dp[1]
  xi2 <- dp[2]
  Omega11 <- dp[3]
  Omega12 <- dp[4]
  Omega22 <- dp[5]
  alpha1 <- dp[6]
  alpha2 <- dp[7]
  tau <- dp[8]
  
  omega1 <- sqrt(Omega11)
  omega2 <- sqrt(Omega22)
  z1 <- (y1 - xi1)/omega1
  z2 <- (y2 - xi2)/omega2
  lambda <- Omega12/(omega1 * omega2)
  u <- 1/(1 - lambda^2)
  r.1pasq <- sqrt(1 + alpha1^2 + 2 * alpha1 * alpha2 * lambda + alpha2^2)
  alpha0 <- tau * r.1pasq
  t <- alpha0 + alpha1 * z1 + alpha2 * z2
  zeta1 <- (sn::zeta(1, t))
  
  s.xi1 <- sum((z1 - lambda*z2)*u - alpha1 *zeta1)/omega1
  s.xi2 <- sum((z2 - lambda*z1)*u - alpha2 *zeta1)/omega2
  w <- (z1^2 + z2^2 - 2*z1*z2*lambda) * lambda * u^2
  s.O11 <- sum(w*lambda + (z1^2- 2*z1*z2*lambda -1)*u 
              -(alpha1*alpha2*lambda*tau/r.1pasq + alpha1*z1)*zeta1)/(2*Omega11)
  s.O12 <- sum((lambda + z1*z2)*u - w + alpha1*alpha2*tau*zeta1/r.1pasq)/(omega1*omega2)
  s.O22 <- sum(w*lambda + (z2^2- 2*z1*z2*lambda -1)*u
              -(alpha1*alpha2*lambda*tau/r.1pasq + alpha2*z2)*zeta1)/(2*Omega22)
  s.a1 <- sum(((alpha1 + alpha2*lambda)*tau/r.1pasq + z1)*zeta1)
  s.a2 <- sum(((alpha2 + alpha1*lambda)*tau/r.1pasq + z2)*zeta1)           
  s.tau <- sum(r.1pasq*zeta1 - sn::zeta(1, tau))    
  scores <- c(s.xi1, s.xi2, s.O11, s.O12, s.O22, s.a1, s.a2, s.tau)
  names(scores) <- names(dp)
  return(scores)    
}


#------------------------------------------------------------------------------ 
#  ESN_2 observed information matrix
#------------------------------------------------------------------------------ 
# Description
# Computes the observed information matrix for the extended bivariate 
# skew-normal distribution in the DP parametrizations.
#
# Usage
# esn2.Oinfo(dp, y1, y2)
#
# Arguments:
#   dp      vector of direct parameters
#   y1, y2  vectors of observations, must have equal length
# 
# Example:
#   y1 <- c(0.7, -0.4, 1.1) 
#   y2 <- c(-1.2, 0.5, 0.9)   
#   info <- esn2.Oinfo(dp=c(0, 0, 1, 0.6, 1, 2, 3, 1), y1, y2)
# where
#   location:  xi1 <- dp[1];   xi2 <- dp[2]
#   scale:     Omega11 <- dp[3]; Omega12 <- dp[4]; Omega22 <- dp[5]
#   trunc:     tau <- dp[8]
#------------------------------------------------------------------------------
esn2.Oinfo <- function(dp, y1, y2) 
{
  if(length(y1) != length(y2)) stop("y1, y2 must have equal length") 
  dp.ok  <- (length(dp) == 8) & (dp[3]>0) & (dp[5]>0) & (dp[3]*dp[5]- dp[4]^2>0) 
  if(!dp.ok) stop("inadmissible dp vector") 
  
  xi1 <- dp[1]
  xi2 <- dp[2]
  xi <- c(xi1,xi2) 
  Omega11 <- dp[3]
  Omega12 <- dp[4]
  Omega22 <- dp[5]
  Omega <- matrix(c(Omega11,Omega12,Omega12,Omega22),2,2)
  alpha1 <- dp[6]
  alpha2 <- dp[7]
  alpha <- c(alpha1,alpha2)
  tau <- dp[8]
  
  z1 <- (y1-xi1)/sqrt(Omega11)
  z2 <- (y2-xi2)/sqrt(Omega22)
  lambda <-  Omega12/sqrt(Omega11*Omega22)
  alpha0 <- tau*sqrt(1 + (alpha1^2) + (alpha2^2) + 2*alpha1*alpha2*lambda)
  t <- (alpha0 + alpha1*z1 + alpha2*z2)
  alphastar <- (alpha1^2) + 2*alpha1*alpha2*lambda + (alpha2^2)
  den <- sqrt(1+alphastar)

  zeta1 <- (sn::zeta(1,t))
  zeta2 <- (sn::zeta(2,t))
  zeta1tau <- (sn::zeta(1,tau))
  zeta2tau <- (sn::zeta(2,tau))
  
  j11 <- sum( (-1/Omega11)*( (1-lambda^2)^(-1) - (alpha1^2)*zeta2)) 
  
  j12 <- sum((1/sqrt(Omega11*Omega22))*((lambda/(1-lambda^2))+alpha1*alpha2*zeta2))
    
  j13 <-  sum((lambda*z2-z1)/((1-lambda^2)^2*(Omega11^(3/2))) + 
      (alpha1/(2*Omega11^(3/2))) * ((alpha1*alpha2*lambda*tau)/den + alpha1*z1)*zeta2 + 
      (alpha1/(2*Omega11^(3/2)))*zeta1)
  
  j14 <- sum(-(2*lambda*(lambda*z2-z1))/((1-lambda^2)^2 * Omega11*sqrt(Omega22)) - 
        (z2)/((1-lambda^2) * Omega11*sqrt(Omega22)) -
        (alpha1^2*alpha2*tau)/(Omega11*sqrt(Omega22)*den)*zeta2)
  
  j15 <- sum((lambda*(z2-z1*lambda))/((1-lambda^2)^2*Omega22*sqrt(Omega11))
         +(alpha1/(2*Omega22*sqrt(Omega11)))*((alpha1*alpha2*lambda*tau)/den+alpha2*z2)*zeta2)

        
  j16 <- sum(-(alpha1/sqrt(Omega11))*( ((alpha1+lambda*alpha2)*tau)/(den) + z1 )*zeta2 
         - (1/sqrt(Omega11)) * zeta1)
  
  j17 <- sum(-(alpha1/sqrt(Omega11))*( ((alpha2+lambda*alpha1)*tau)/(den) + z2 )*zeta2)
  
  j18 <- sum(-((alpha1*den)/sqrt(Omega11))*zeta2)
    
  j22 <- sum((-1/Omega22)*( (1-lambda^2)^(-1) - (alpha2^2)*zeta2))
  
  j23 <- sum((lambda*(z1-z2*lambda)) / ((1-lambda^2)^2 * Omega11*sqrt(Omega22)) +
       (alpha2/(2*Omega11*sqrt(Omega22))) * ((alpha1*alpha2*lambda*tau)/den + alpha1*z1) *zeta2)
  
  j24 <- sum(-(2*lambda*(lambda*z1-z2))/((1-lambda^2)^2 * Omega22*sqrt(Omega11)) - 
        (z1)/((1-lambda^2) * Omega22*sqrt(Omega11)) -
        (alpha2^2*alpha1*tau)/(Omega22*sqrt(Omega11)*den)*zeta2)
  
  j25 <- sum((lambda*z1-z2)/((1-lambda^2)^2*(Omega22^(3/2))) + 
        (alpha2/(2*Omega22^(3/2)))*((alpha1*alpha2*lambda*tau)/den + alpha2*z2)*zeta2 + 
        (alpha2/(2*Omega22^(3/2)))*zeta1)
  
  j26 <- sum(-(alpha2/sqrt(Omega22))*( ((alpha1+lambda*alpha2)*tau)/(den) + z1 )*zeta2)
  
  j27 <- sum(-(alpha2/sqrt(Omega22))*( ((alpha2+lambda*alpha1)*tau)/(den) + z2 )*zeta2 
         - (1/sqrt(Omega22)) * zeta1)  
  
  j28 <- sum(-((alpha2*den)/sqrt(Omega22))*zeta2)
  

  j33 <- sum(  (lambda^2-z1^2+2*z1*z2*lambda)/((Omega11^2)*(1-lambda^2)) + 
      (4*lambda^3*z1*z2-2*lambda^2*z1^2-lambda^2*z2^2)/((Omega11^2)*(1-lambda^2)^2) -
      (lambda^4*(z1^2+z2^2-2*z1*z2*lambda))/((Omega11^2)*(1-lambda^2)^3) + 1/(2*Omega11^2) +
      lambda^4/((2*Omega11^2)*(1-lambda^2)^2) + 
      (1/(4*Omega11^2))*
      ( (3*alpha1*alpha2*tau*lambda)/den - (alpha1^2*alpha2^2*tau*lambda^2)/den^(3) 
       + 3*alpha1*z1 )* zeta1 +
      (1/(4*Omega11^2))*( (alpha1*alpha2*tau*lambda)/den + alpha1*z1 )^2*zeta2 )
        
  j34 <- sum(  (-1*(lambda+z1*z2)/((1-lambda^2)*(Omega11^(3/2))*sqrt(Omega22)))+
      (2*lambda*z1^2+lambda*z2^2-5*lambda^2*z1*z2-lambda^3)
          /(((1-lambda^2)^2)*(Omega11^(3/2))*sqrt(Omega22))+
      (2*lambda^3*(z1^2+z2^2-2*z1*z2*lambda))/(((1-lambda^2)^3)*(Omega11^(3/2))*sqrt(Omega22)) + 
      (((alpha1^2)*(alpha2^2)*tau*lambda)/(2*(Omega11^(3/2))*sqrt(Omega22)*den^3)
        -(alpha1*alpha2*tau)/((2*Omega11^(3/2))*sqrt(Omega22)*den))*zeta1-
      (((alpha1*alpha2*tau)/(2*(Omega11^(3/2))*sqrt(Omega22)*den))
        *((alpha1*alpha2*lambda*tau)/den + alpha1*z1)*zeta2) )
      
  j35 <- sum( lambda^2*(6*lambda*z1*z2-2*z1^2-2*z2^2+lambda^2)/(2*Omega11*Omega22*(1-lambda^2)^2)+
        (2*z1*z2*lambda+lambda^2)/(2*Omega11*Omega22*(1-lambda^2)) 
        - lambda^4*(z1^2+z2^2-2*z1*z2*lambda)/(Omega11*Omega22*(1-lambda^2)^3) 
        + ((alpha1*alpha2*lambda*tau)/(4*Omega11*Omega22*den)) 
        * (1-(alpha1*alpha2*lambda)/(1+alphastar))*zeta1 
        + (1/(4*Omega11*Omega22)) * ((alpha1*alpha2*lambda*tau)/den + alpha1*z1) 
        * ((alpha1*alpha2*lambda*tau)/den + alpha2*z2)*zeta2   ) 

  j36 <- sum((1/(2*Omega11))*( (alpha1*alpha2*lambda*(alpha2*lambda+alpha1)*tau)/(den^(3)) 
        - (alpha2*lambda*tau)/(den) -z1 )*zeta1 
        - (1/(2*Omega11)) * ((alpha1*alpha2*lambda*tau)/den + alpha1*z1) 
        * (((alpha2*lambda+alpha1)*tau)/den + z1)*zeta2)
    
  j37 <- sum( (1/(2*Omega11))*( (alpha1*alpha2*lambda*(alpha1*lambda+alpha2)*tau)/(den^(3))
         -(alpha1*lambda*tau)/(den))*zeta1 - (1/(2*Omega11)) * ((alpha1*alpha2*lambda*tau)/(den) 
          + alpha1*z1) *(((alpha2+alpha1*lambda)*tau)/(den) + z2)*zeta2)
  
  j38 <- sum( -((alpha1*alpha2*lambda)/(2*Omega11*den))*zeta1 -
       (den/(2*Omega11))*((alpha1*alpha2*lambda*tau)/den + alpha1*z1 )*zeta2 )

  j44 <- sum( (1/((1-lambda^2)*Omega11*Omega22))+(6*lambda*z1*z2-z1^2-z2^2+2*lambda^2)
         /(((1-lambda^2)^2)*Omega11*Omega22)
         -(4*lambda^2* (z1^2+z2^2-2*z1*z2*lambda))/(((1-lambda^2)^3)*Omega11*Omega22)+
         ((alpha1^2*alpha2^2*tau)/(Omega11*Omega22*den^2))*(tau*zeta2-zeta1/den)  )           

  j45 <- sum(-(lambda+z1*z2)/((1-lambda^2)*Omega22^(3/2)*sqrt(Omega11)) + 
           (2*lambda*z2^2+lambda*z1^2-5*lambda^2*z1*z2-lambda^3)
           /((1-lambda^2)^2*Omega22^(3/2)*sqrt(Omega11)) + 
          (2*lambda^3*(z1^2+z2^2-2*z1*z2*lambda))/((1-lambda^2)^3*Omega22^(3/2)*sqrt(Omega11)) +
          ( (alpha1^2*alpha2^2*tau*lambda)/(2*Omega22^(3/2)*sqrt(Omega11)*den^3) - 
          (alpha1*alpha2*tau)/(2*Omega22^(3/2)*sqrt(Omega11)*den) )*zeta1 - 
          ((alpha1*alpha2*tau)/(2*Omega22^(3/2)*sqrt(Omega11)*den)) 
          * ((alpha1*alpha2*lambda*tau)/den + alpha2*z2)*zeta2)
  
  j46 <- sum((alpha2*tau)/(sqrt(Omega11*Omega22)*den)*(1-(alpha1*(alpha2*lambda+alpha1))/den^2)*zeta1
         +(alpha1*alpha2*tau)/(sqrt(Omega11*Omega22)*den) * (((alpha1+lambda*alpha2)*tau)/den 
         + z1)*zeta2 )
  
  j47 <- sum( 
        (alpha1*tau)/(sqrt(Omega11*Omega22)*den)*(1-(alpha2*(alpha1*lambda+alpha2))/den^2)*zeta1 +
        (alpha1*alpha2*tau)/(sqrt(Omega11*Omega22)*den) * (((alpha2+lambda*alpha1)*tau)/den 
        + z2)*zeta2 )

  j48 <- sum(((alpha1*alpha2)/sqrt(Omega11*Omega22))*(zeta1/den + tau*zeta2))
   
  j55 <- sum( (lambda^2-z2^2+2*z1*z2*lambda)/((Omega22^2)*(1-lambda^2)) + 
      (4*lambda^3*z1*z2-2*lambda^2*z2^2-lambda^2*z1^2)/((Omega22^2)*(1-lambda^2)^2) -
      (lambda^4*(z1^2+z2^2-2*z1*z2*lambda))/((Omega22^2)*(1-lambda^2)^3) + 1/(2*Omega22^2) +
      lambda^4/((2*Omega22^2)*(1-lambda^2)^2) + 
      (1/(4*Omega22^2))*( (3*alpha1*alpha2*tau*lambda)/den - 
      (alpha1^2*alpha2^2*tau*lambda^2)/den^(3) + 3*alpha2*z2 )* zeta1 +
      (1/(4*Omega22^2))*( (alpha1*alpha2*tau*lambda)/den + alpha2*z2 )^2*zeta2 )
  
  j56 <- sum( (1/(2*Omega22))*( (alpha1*alpha2*lambda*(alpha2*lambda+alpha1)*tau)/(den^(3))-
      (alpha2*lambda*tau)/(den))*zeta1 -
      (1/(2*Omega22)) * ((alpha1*alpha2*lambda*tau)/(den) + alpha2*z2) * 
      (((alpha1+alpha2*lambda)*tau)/(den) + z1)*zeta2 )
  
  j57 <- sum( (1/(2*Omega22))*( (alpha1*alpha2*lambda*(alpha1*lambda+alpha2)*tau)/(den^(3)) - 
      (alpha1*lambda*tau)/(den) -z2 )*zeta1 -
      (1/(2*Omega22)) * ((alpha1*alpha2*lambda*tau)/den + alpha2*z2) * 
      (((alpha1*lambda+alpha2)*tau)/den + z2)*zeta2 )
  
  j58 <- sum( -((alpha1*alpha2*lambda)/(2*Omega22*den))*zeta1 -
        (den/(2*Omega22))*((alpha1*alpha2*lambda*tau)/den + alpha2*z2 )*zeta2)
  
  j66 <- sum( (tau/den-(((alpha2*lambda+alpha1)^2*tau)/den^(3)))*zeta1+ 
          (((alpha1+lambda*alpha2)*tau)/den + z1)^2*zeta2 )
      
  j67 <- sum( ((lambda*tau)/den - ((alpha2+lambda*alpha1)*(alpha1+lambda*alpha2)*tau)/den^3)*zeta1 
        + (((alpha1+lambda*alpha2)*tau)/den + z1) *(((alpha2+lambda*alpha1)*tau)/den + z2)*zeta2 )
  
  j68 <- sum( ((alpha1+lambda*alpha2)/den)*zeta1+(((alpha1+lambda*alpha2)*tau)/den + z1)*den*zeta2 )

  j77<- sum( (tau/den - (((alpha1*lambda+alpha2)^2*tau)/den^(3)) )*zeta1 + 
        (((alpha2+lambda*alpha1)*tau)/den + z2 )^2 * zeta2) 
  
  j78<- sum(((alpha2+lambda*alpha1)/den)*zeta1 + (((alpha2+lambda*alpha1)*tau)/den + z2)*den*zeta2)

  j88<- sum(den^2*zeta2-zeta2tau)
  
  r1 <- c(j11,j12,j13,j14,j15,j16,j17,j18)
  r2 <- c(j12,j22,j23,j24,j25,j26,j27,j28)
  r3 <- c(j13,j23,j33,j34,j35,j36,j37,j38)
  r4 <- c(j14,j24,j34,j44,j45,j46,j47,j48)
  r5 <- c(j15,j25,j35,j45,j55,j56,j57,j58)
  r6 <- c(j16,j26,j36,j46,j56,j66,j67,j68)
  r7 <- c(j17,j27,j37,j47,j57,j67,j77,j78)
  r8 <- c(j18,j28,j38,j48,j58,j68,j78,j88)
  
  j.dp <- matrix(rbind(r1,r2,r3,r4,r5,r6,r7,r8),8,8)
  return(j.dp)
}
\end{verbatim}
\end{footnotesize}


\subsection*{A.2: The expected information matrix}

The \Rlang\ code below computes the $\ESN_{2}$ expected information matrix.
It assumes that the \Rlang\ packages  \texttt{sn} \citep{azzalini:R-pkg-sn} 
and \texttt{cubature} \citep{narasimhan:etal:R-pkg-cubature} have been loaded.

\begin{footnotesize}
\begin{verbatim}
 #-----------------------------------------------------------------------------------------------------
# ESN_2 expected information matrix
#-----------------------------------------------------------------------------------------------------
# Description
# Computes the expected information matrix for the extended bivariate skew normal distribution 
# in the DP parametrizations.
#
# Usage
# esn2.Einfo(dp)
#
# Arguments
# dp   vector of direct parameters
#-----------------------------------------------------------------------------------------------------
# Example: 
#   info <- esn2.Einfo(dp=c(0, 0, 1, 0.6, 1, 2, 3, 1))
#
# where the 'dp' components are intended as follows:
# location:  xi1 <- dp[1] ; xi2 <- dp[2]
# scale:     Omega11 <- dp[3]; Omega12 <- dp[4];  Omega22 <- dp[5]
# trunc:     tau <- dp[8]
#-----------------------------------------------------------------------------------------------------
esn2.Einfo <- function(dp) 
{    
    dp.ok  <- (length(dp) == 8) & (dp[3]>0) & (dp[5]>0) & (dp[3]*dp[5]- dp[4]^2>0) 
    if(!dp.ok) stop("inadmissible dp vector") 
  
    xi1 <- dp[1]
    xi2 <- dp[2]
    xi <- c(xi1, xi2)
    
    Omega11 <- dp[3]
    Omega12 <- dp[4]
    Omega22 <- dp[5]
    Omega <- matrix(c(Omega11, Omega12, Omega12, Omega22), 2, 2)
    
    alpha1 <- dp[6]
    alpha2 <- dp[7]
    alpha <- c(alpha1, alpha2)
    
    tau <- dp[8]
    
    lambda <- Omega12/sqrt(Omega11 * Omega22)
    Omegab <- matrix(c(1, lambda, lambda, 1), 2, 2)
    alphastar <- (alpha1^2) + (alpha2^2) + 2 * alpha1 * alpha2 * lambda
    alpha0 <- tau * (sqrt(1 + alphastar))
    den <- sqrt(1 + alphastar)
    
    zeta1tau <- (sn::zeta(1, tau))
    zeta2tau <- (sn::zeta(2, tau))
    
    delta <- Omegab %*% alpha %*% (1/(sqrt(1 + t(alpha) %*% Omegab %*% alpha)))
    delta1 <- (alpha1 + lambda * alpha2)/den
    delta2 <- (alpha2 + lambda * alpha1)/den
    
    v11 <- (1 + alpha2^2 * (1 - lambda^2))/(1 + alphastar)
    v22 <- (1 + alpha1^2 * (1 - lambda^2))/(1 + alphastar)
    v12 <- (lambda - alpha1 * alpha2 * (1 - lambda^2))/(1 + alphastar)
    
    E_Z <- zeta1tau * delta
    E_Z1 <- zeta1tau * delta1
    E_Z2 <- zeta1tau * delta2
    E_Z1Z2 <- lambda + delta1 * delta2 * (zeta1tau^2 + zeta2tau)
    E_Z1q <- 1 + delta1^2 * (zeta2tau + zeta1tau^2)
    E_Z2q <- 1 + delta2^2 * (zeta2tau + zeta1tau^2)
    
    zeta1q <- function(z) {
        (sn::zeta(1, (alpha0 + alpha1 * z[1] + alpha2 * z[2])))^2 * dmsn(z, xi = rep(0, 
            2), Omega = Omegab, alpha = alpha, tau = tau)
    }
    a0 <- cuhre(zeta1q, 1, lower = rep(-5, 2), upper = rep(5, 2), relTol = 1e-3, 
        absTol = 1e-15, flags = list(verbose = 0, final = 0))$integral
    
    # a1_1
    
    z1_zeta1q <- function(z) {
        z[1] * (sn::zeta(1, (alpha0 + alpha1 * z[1] + alpha2 * z[2])))^2 * dmsn(z, 
            xi = rep(0, 2), Omega = Omegab, alpha = alpha, tau = tau)
    }
    a1_1 <- cuhre(z1_zeta1q, 1, lower = rep(-5, 2), upper = rep(5, 2), relTol = 1e-3, 
        absTol = 1e-15, flags = list(verbose = 0, final = 0))$integral
    
    # a2_1
    
    z2_zeta1q <- function(z) {
        z[2] * (sn::zeta(1, (alpha0 + alpha1 * z[1] + alpha2 * z[2])))^2 * dmsn(z, 
            xi = rep(0, 2), Omega = Omegab, alpha = alpha, tau = tau)
    }
    a2_1 <- cuhre(z2_zeta1q, 1, lower = rep(-5, 2), upper = rep(5, 2), relTol = 1e-3, 
        absTol = 1e-15, flags = list(verbose = 0, final = 0))$integral
    
    # a1_2
    
    Z1q_zeta1q <- function(z) {
        z[1]^2 * (sn::zeta(1, (alpha0 + alpha1 * z[1] + alpha2 * z[2])))^2 * dmsn(z, 
            xi = rep(0, 2), Omega = Omegab, alpha = alpha, tau = tau)
    }
    a1_2 <- cuhre(Z1q_zeta1q, 1, lower = rep(-5, 2), upper = rep(5, 2), relTol = 1e-3, 
        absTol = 1e-15, flags = list(verbose = 0, final = 0))$integral
        
    # a2_2
    
    Z2q_zeta1 <- function(z) {
        z[2]^2 * (sn::zeta(1, (alpha0 + alpha1 * z[1] + alpha2 * z[2])))^2 * dmsn(z, 
            xi = rep(0, 2), Omega = Omegab, alpha = alpha, tau = tau)
    }
    a2_2 <- cuhre(Z2q_zeta1, 1, lower = rep(-5, 2), upper = rep(5, 2), relTol = 1e-3, 
        absTol = 1e-15, flags = list(verbose = 0, final = 0))$integral
    
    # a12
    
    Z1Z2_zeta1q <- function(z) {
        z[1] * z[2] * (sn::zeta(1, (alpha0 + alpha1 * z[1] + alpha2 * z[2])))^2 * 
            dmsn(z, xi = rep(0, 2), Omega = Omegab, alpha = alpha, tau = tau)
    }
    a12 <- cuhre(Z1Z2_zeta1q, 1, lower = rep(-5, 2), upper = rep(5, 2), relTol = 1e-3, 
        absTol = 1e-15, flags = list(verbose = 0, final = 0))$integral
    
    
    Ezeta_1 <- (sn::zeta(1, tau))/den
    
    E_Z1_zeta1 <- -tau * delta1 * Ezeta_1
    
    E_Z2_zeta1 <- -tau * delta2 * Ezeta_1
    
    Ezeta_2 <- -(alpha0 - tau * t(alpha) %*% delta) * Ezeta_1 - a0
    
    E_Z1_zeta2 <- -(-alpha0 * tau * delta1 + alpha1 * (tau^2 * delta1^2 + v11) + 
        alpha2 * (tau^2 * delta1 * delta2 + v12)) * Ezeta_1 - a1_1
    
    E_Z2_zeta2 <- -(-alpha0 * tau * delta2 + alpha2 * (tau^2 * delta2^2 + v22) + 
        alpha1 * (tau^2 * delta1 * delta2 + v12)) * Ezeta_1 - a2_1
    
    EZ1q_zeta2 <- -(alpha0 * (v11 + tau^2 * delta1^2) - alpha1 * tau * delta1 * (tau^2 * 
        delta1^2 + 3 * v11) + alpha2 * tau * ((v12/v11) * delta1 - delta2) * (tau^2 * 
        delta1^2 + v11) - alpha2 * tau * delta1 * (v12/v11) * (tau^2 * delta1^2 + 
        3 * v11)) * Ezeta_1 - a1_2
    
    EZ2q_zeta2 <- -(alpha0 * (v22 + tau^2 * delta2^2) - alpha2 * tau * delta2 * (tau^2 * 
        delta2^2 + 3 * v22) + alpha1 * tau * ((v12/v22) * delta2 - delta1) * (tau^2 * 
        delta2^2 + v22) - alpha1 * tau * delta2 * (v12/v22) * (tau^2 * delta2^2 + 
        3 * v22)) * Ezeta_1 - a2_2
    
    E_Z1Z2_zeta2 <- -(alpha0 * (v12 + tau^2 * delta1 * delta2) + alpha1 * tau * ((v12/v11) * 
        delta1 - delta2) * (tau^2 * delta1^2 + v11) - alpha1 * tau * delta1 * (v12/v11) * 
        (tau^2 * delta1^2 + 3 * v11) + alpha2 * tau * ((v12/v22) * delta2 - delta1) * 
        (tau^2 * delta2^2 + v22) - alpha2 * tau * delta2 * (v12/v22) * (tau^2 * delta2^2 + 
        3 * v22)) * Ezeta_1 - a12
    
    
    i11 <- (1/Omega11) * ((1 - lambda^2)^(-1) - (alpha1^2) * Ezeta_2)
    
    i12 <- (-1/sqrt(Omega11 * Omega22)) * ((lambda/(1 - lambda^2)) + alpha1 * alpha2 * 
        Ezeta_2)
    
    i13 <- -(lambda * E_Z2 - E_Z1)/((1 - lambda^2)^2 * (Omega11^(3/2))) - (alpha1/(2 * 
        Omega11^(3/2))) * (((alpha1 * alpha2 * lambda * tau)/den) * Ezeta_2 + alpha1 * 
        E_Z1_zeta2) - (alpha1/(2 * Omega11^(3/2))) * Ezeta_1
    
    i14 <- (2 * lambda * (lambda * E_Z2 - E_Z1))/((1 - lambda^2)^2 * Omega11 * sqrt(Omega22)) + 
        (E_Z2)/((1 - lambda^2) * Omega11 * sqrt(Omega22)) + ((alpha1^2 * alpha2 * 
        tau)/(Omega11 * sqrt(Omega22) * den)) * Ezeta_2
    
    i15 <- -(lambda * (E_Z2 - E_Z1 * lambda))/((1 - lambda^2)^2 * Omega22 * sqrt(Omega11)) - 
        (alpha1/(2 * Omega22 * sqrt(Omega11))) * (((alpha1 * alpha2 * lambda * tau)/den) * 
            Ezeta_2 + alpha2 * E_Z2_zeta2)
    
    i16 <- (alpha1/sqrt(Omega11)) * ((((alpha1 + lambda * alpha2) * tau)/(den)) * 
        Ezeta_2 + E_Z1_zeta2) + (1/sqrt(Omega11)) * Ezeta_1
    
    i17 <- (alpha1/sqrt(Omega11)) * ((((alpha2 + lambda * alpha1) * tau)/(den)) * 
        Ezeta_2 + E_Z2_zeta2)
    
    i18 <- ((alpha1 * den)/sqrt(Omega11)) * Ezeta_2
    
    
    i22 <- (1/Omega22) * ((1 - lambda^2)^(-1) - (alpha2^2) * Ezeta_2)
    
    i23 <- -(lambda * (E_Z1 - E_Z2 * lambda))/((1 - lambda^2)^2 * Omega11 * sqrt(Omega22)) - 
        (alpha2/(2 * Omega11 * sqrt(Omega22))) * (((alpha1 * alpha2 * lambda * tau)/den) * 
            Ezeta_2 + alpha1 * E_Z1_zeta2)
    
    i24 <- (2 * lambda * (lambda * E_Z1 - E_Z2))/((1 - lambda^2)^2 * Omega22 * sqrt(Omega11)) + 
        (E_Z1)/((1 - lambda^2) * Omega22 * sqrt(Omega11)) + ((alpha2^2 * alpha1 * 
        tau)/(Omega22 * sqrt(Omega11) * den)) * Ezeta_2
    
    i25 <- -(lambda * E_Z1 - E_Z2)/((1 - lambda^2)^2 * (Omega22^(3/2))) - (alpha2/(2 * 
        Omega22^(3/2))) * (((alpha1 * alpha2 * lambda * tau)/den) * Ezeta_2 + alpha2 * 
        E_Z2_zeta2) - (alpha2/(2 * Omega22^(3/2))) * Ezeta_1
    
    i26 <- (alpha2/sqrt(Omega22)) * ((((alpha1 + lambda * alpha2) * tau)/(den)) * 
        Ezeta_2 + E_Z1_zeta2)
    
    i27 <- (alpha2/sqrt(Omega22)) * ((((alpha2 + lambda * alpha1) * tau)/(den)) * 
        Ezeta_2 + E_Z2_zeta2) + (1/sqrt(Omega22)) * Ezeta_1
    
    i28 <- ((alpha2 * den)/sqrt(Omega22)) * Ezeta_2
       
    i33 <- -(lambda^2 - E_Z1q + 2 * E_Z1Z2 * lambda)/((Omega11^2) * (1 - lambda^2)) - 
        (4 * lambda^3 * E_Z1Z2 - 2 * lambda^2 * E_Z1q - lambda^2 * E_Z2q)/((Omega11^2) * 
            (1 - lambda^2)^2) + (lambda^4 * (E_Z1q + E_Z2q - 2 * E_Z1Z2 * lambda))/((Omega11^2) * 
        (1 - lambda^2)^3) - 1/(2 * Omega11^2) - lambda^4/((2 * Omega11^2) * (1 - 
        lambda^2)^2) - (1/(4 * Omega11^2)) * (((3 * alpha1 * alpha2 * tau * lambda)/den) * 
        Ezeta_1 - ((alpha1^2 * alpha2^2 * tau * lambda^2)/den^(3)) * Ezeta_1 + 3 * 
        alpha1 * E_Z1_zeta1) - (1/(4 * Omega11^2)) * (((alpha1^2 * alpha2^2 * tau^2 * 
        lambda^2)/den^2) * Ezeta_2 + alpha1^2 * EZ1q_zeta2 + ((2 * alpha1^2 * alpha2 * 
        tau * lambda)/den) * E_Z1_zeta2)
    
    i34 <- ((lambda + E_Z1Z2)/((1 - lambda^2) * (Omega11^(3/2)) * sqrt(Omega22))) - 
        (2 * lambda * E_Z1q + lambda * E_Z2q - 5 * lambda^2 * E_Z1Z2 - lambda^3)/(((1 - 
            lambda^2)^2) * (Omega11^(3/2)) * sqrt(Omega22)) - (2 * lambda^3 * (E_Z1q + 
        E_Z2q - 2 * E_Z1Z2 * lambda))/(((1 - lambda^2)^3) * (Omega11^(3/2)) * sqrt(Omega22)) - 
        (((alpha1^2) * (alpha2^2) * tau * lambda)/(2 * (Omega11^(3/2)) * sqrt(Omega22) * 
            den^3) - (alpha1 * alpha2 * tau)/((2 * Omega11^(3/2)) * sqrt(Omega22) * 
            den)) * Ezeta_1 + (((alpha1 * alpha2 * tau)/(2 * (Omega11^(3/2)) * sqrt(Omega22) * 
        den)) * (((alpha1 * alpha2 * lambda * tau)/den) * Ezeta_2 + alpha1 * E_Z1_zeta2))
    
    i35 <- -(lambda^4 + 6 * lambda^3 * E_Z1Z2 - 2 * lambda^2 * E_Z1q - 2 * lambda^2 * 
        E_Z2q)/(2 * Omega11 * Omega22 * (1 - lambda^2)^2) - (lambda^2)/(2 * Omega11 * 
        Omega22 * (1 - lambda^2)) - (lambda * E_Z1Z2)/(Omega11 * Omega22 * (1 - lambda^2)) + 
        (lambda^4 * (E_Z1q + E_Z2q - 2 * E_Z1Z2 * lambda))/(Omega11 * Omega22 * (1 - 
            lambda^2)^3) - ((alpha1 * alpha2 * lambda * tau)/(4 * Omega11 * Omega22 * 
        den)) * (1 - (alpha1 * alpha2 * lambda)/(1 + alphastar)) * Ezeta_1 - ((alpha1 * 
        alpha2)/(4 * Omega11 * Omega22)) * (((alpha1 * alpha2 * lambda^2 * tau^2)/den^2) * 
        Ezeta_2 + ((alpha2 * lambda * tau)/den) * E_Z2_zeta2 + ((alpha1 * lambda * 
        tau)/den) * E_Z1_zeta2 + E_Z1Z2_zeta2)
    
    i36 <- -(1/(2 * Omega11)) * (((alpha1 * alpha2 * lambda * (alpha2 * lambda + 
        alpha1) * tau)/(den^(3))) * Ezeta_1 - ((alpha2 * lambda * tau)/(den)) * Ezeta_1 - 
        E_Z1_zeta1) + (1/(2 * Omega11)) * (((alpha1 * alpha2 * lambda * tau^2 * (alpha2 * 
        lambda + alpha1))/den^2) * Ezeta_2 + ((alpha1 * alpha2 * lambda * tau)/den) * 
        E_Z1_zeta2 + ((alpha1 * (alpha2 * lambda + alpha1) * tau)/den) * E_Z1_zeta2 + 
        alpha1 * EZ1q_zeta2)
    
    i37 <- -(1/(2 * Omega11)) * ((alpha1 * alpha2 * lambda * (alpha1 * lambda + alpha2) * 
        tau)/(den^(3)) - (alpha1 * lambda * tau)/(den)) * Ezeta_1 + (1/(2 * Omega11)) * 
        (((alpha1 * alpha2 * lambda * tau^2 * (alpha2 + alpha1 * lambda))/(den)^2) * 
            Ezeta_2 + ((alpha1 * alpha2 * lambda * tau)/(den)) * E_Z2_zeta2 + ((alpha1 * 
            (alpha2 + alpha1 * lambda) * tau)/(den)) * E_Z1_zeta2 + alpha1 * E_Z1Z2_zeta2)
    
    i38 <- ((alpha1 * alpha2 * lambda)/(2 * Omega11 * den)) * Ezeta_1 + (den/(2 * 
        Omega11)) * (((alpha1 * alpha2 * lambda * tau)/den) * Ezeta_2 + alpha1 * 
        E_Z1_zeta2)
    
    i44 <- -(1/((1 - lambda^2) * Omega11 * Omega22)) - (6 * lambda * E_Z1Z2 - E_Z1q - 
        E_Z2q + 2 * lambda^2)/(((1 - lambda^2)^2) * Omega11 * Omega22) + (4 * lambda^2 * 
        (E_Z1q + E_Z2q - 2 * E_Z1Z2 * lambda))/(((1 - lambda^2)^3) * Omega11 * Omega22) - 
        ((alpha1^2 * alpha2^2 * tau)/(Omega11 * Omega22 * den^2)) * (tau * Ezeta_2 - 
            Ezeta_1/den)
    
    i45 <- (lambda + E_Z1Z2)/((1 - lambda^2) * Omega22^(3/2) * sqrt(Omega11)) - (2 * 
        lambda * E_Z2q + lambda * E_Z1q - 5 * lambda^2 * E_Z1Z2 - lambda^3)/((1 - 
        lambda^2)^2 * Omega22^(3/2) * sqrt(Omega11)) - (2 * lambda^3 * (E_Z1q + E_Z2q - 
        2 * E_Z1Z2 * lambda))/((1 - lambda^2)^3 * Omega22^(3/2) * sqrt(Omega11)) - 
        ((alpha1^2 * alpha2^2 * tau * lambda)/(2 * Omega22^(3/2) * sqrt(Omega11) * 
            den^3) - (alpha1 * alpha2 * tau)/(2 * Omega22^(3/2) * sqrt(Omega11) * 
            den)) * Ezeta_1 + ((alpha1 * alpha2 * tau)/(2 * Omega22^(3/2) * sqrt(Omega11) * 
        den)) * (((alpha1 * alpha2 * lambda * tau)/den) * Ezeta_2 + alpha2 * E_Z2_zeta2)
    
    i46 <- -(alpha2 * tau)/(sqrt(Omega11 * Omega22) * den) * (1 - (alpha1 * (alpha2 * 
        lambda + alpha1))/den^2) * Ezeta_1 - (alpha1 * alpha2 * tau)/(sqrt(Omega11 * 
        Omega22) * den) * ((((alpha1 + lambda * alpha2) * tau)/den) * Ezeta_2 + E_Z1_zeta2)
    
    i47 <- -(alpha1 * tau)/(sqrt(Omega11 * Omega22) * den) * (1 - (alpha2 * (alpha1 * 
        lambda + alpha2))/den^2) * Ezeta_1 - ((alpha1 * alpha2 * tau)/(sqrt(Omega11 * 
        Omega22) * den)) * ((((alpha2 + lambda * alpha1) * tau)/den) * Ezeta_2 + 
        E_Z2_zeta2)
    
    i48 <- -((alpha1 * alpha2)/sqrt(Omega11 * Omega22)) * (Ezeta_1/den + tau * Ezeta_2)
    
    i55 <- -(lambda^2 - E_Z2q + 2 * E_Z1Z2 * lambda)/((Omega22^2) * (1 - lambda^2)) - 
        (4 * lambda^3 * E_Z1Z2 - 2 * lambda^2 * E_Z2q - lambda^2 * E_Z1q)/((Omega22^2) * 
            (1 - lambda^2)^2) + (lambda^4 * (E_Z1q + E_Z2q - 2 * E_Z1Z2 * lambda))/((Omega22^2) * 
        (1 - lambda^2)^3) - 1/(2 * Omega22^2) - lambda^4/((2 * Omega22^2) * (1 - 
        lambda^2)^2) - (1/(4 * Omega22^2)) * (((3 * alpha1 * alpha2 * tau * lambda)/den) * 
        Ezeta_1 - ((alpha1^2 * alpha2^2 * tau * lambda^2)/den^(3)) * Ezeta_1 + 3 * 
        alpha2 * E_Z2_zeta1) - (1/(4 * Omega22^2)) * (((alpha1^2 * alpha2^2 * tau^2 * 
        lambda^2)/den^2) * Ezeta_2 + alpha2^2 * EZ2q_zeta2 + ((2 * alpha1 * alpha2^2 * 
        tau * lambda)/den) * E_Z2_zeta2)
    
    i56 <- -(1/(2 * Omega22)) * ((alpha1 * alpha2 * lambda * (alpha2 * lambda + alpha1) * 
        tau)/(den^(3)) - (alpha2 * lambda * tau)/(den)) * Ezeta_1 + (1/(2 * Omega22)) * 
        (((alpha1 * alpha2 * lambda * tau^2 * (alpha1 + alpha2 * lambda))/(den)^2) * 
            Ezeta_2 + ((alpha1 * alpha2 * lambda * tau)/(den)) * E_Z1_zeta2 + ((alpha2 * 
            (alpha1 + alpha2 * lambda) * tau)/(den)) * E_Z2_zeta2 + alpha2 * E_Z1Z2_zeta2)
    
    i57 <- -(1/(2 * Omega22)) * (((alpha1 * alpha2 * lambda * (alpha1 * lambda + 
        alpha2) * tau)/(den^(3))) * Ezeta_1 - ((alpha1 * lambda * tau)/(den)) * Ezeta_1 - 
        E_Z2_zeta1) + (1/(2 * Omega22)) * (((alpha1 * alpha2 * lambda * tau^2 * (alpha1 * 
        lambda + alpha2))/den^2) * Ezeta_2 + ((alpha1 * alpha2 * lambda * tau)/den) * 
        E_Z2_zeta2 + ((alpha2 * (alpha1 * lambda + alpha2) * tau)/den) * E_Z2_zeta2 + 
        alpha2 * EZ2q_zeta2)
    
    i58 <- ((alpha1 * alpha2 * lambda)/(2 * Omega22 * den)) * Ezeta_1 + (den/(2 * 
        Omega22)) * (((alpha1 * alpha2 * lambda * tau)/den) * Ezeta_2 + alpha2 * 
        E_Z2_zeta2)
    
    i66 <- -(tau/den - (((alpha2 * lambda + alpha1)^2 * tau)/den^(3))) * Ezeta_1 - 
        (((alpha1 + lambda * alpha2)^2 * tau^2)/den^2) * Ezeta_2 - EZ1q_zeta2 - ((2 * 
        (alpha1 + lambda * alpha2) * tau)/den) * E_Z1_zeta2
    
    i67 <- -((lambda * tau)/den - ((alpha2 + lambda * alpha1) * (alpha1 + lambda * 
        alpha2) * tau)/den^3) * Ezeta_1 - (((alpha1 + lambda * alpha2) * (alpha2 + 
        lambda * alpha1) * tau^2)/den^2) * Ezeta_2 - (((alpha1 + lambda * alpha2) * 
        tau)/den) * E_Z2_zeta2 - (((alpha2 + lambda * alpha1) * tau)/den) * E_Z1_zeta2 - 
        E_Z1Z2_zeta2
    
    i68 <- -((alpha1 + lambda * alpha2)/den) * Ezeta_1 - ((((alpha1 + lambda * alpha2) * 
        tau)/den) * Ezeta_2 + E_Z1_zeta2) * den
    
    
    i77 <- -(tau/den - (((alpha1 * lambda + alpha2)^2 * tau)/den^(3))) * Ezeta_1 - 
        (((alpha2 + lambda * alpha1)^2 * tau^2)/den^2) * Ezeta_2 - EZ2q_zeta2 - ((2 * 
        (alpha2 + lambda * alpha1) * tau)/den) * E_Z2_zeta2
    
    i78 <- -((alpha2 + lambda * alpha1)/den) * Ezeta_1 - ((((alpha2 + lambda * alpha1) * 
        tau)/den) * Ezeta_2 + E_Z2_zeta2) * den
    
    i88 <- -(1 + alphastar) * Ezeta_2 + zeta2tau
    
    
    l1 <- c(i11, i12, i13, i14, i15, i16, i17, i18)
    l2 <- c(i12, i22, i23, i24, i25, i26, i27, i28)
    l3 <- c(i13, i23, i33, i34, i35, i36, i37, i38)
    l4 <- c(i14, i24, i34, i44, i45, i46, i47, i48)
    l5 <- c(i15, i25, i35, i45, i55, i56, i57, i58)
    l6 <- c(i16, i26, i36, i46, i56, i66, i67, i68)
    l7 <- c(i17, i27, i37, i47, i57, i67, i77, i78)
    l8 <- c(i18, i28, i38, i48, i58, i68, i78, i88)
        
    I.dp <- matrix(rbind(l1, l2, l3, l4, l5, l6, l7, l8), 8, 8)
    return(I.dp)
}  
\end{verbatim}
\end{footnotesize}
%
\clearpage
\bibliographystyle{apalike2}

\end{document}